\documentclass[12pt]{amsart}

\usepackage{amssymb}
\usepackage{amscd}

\newtheorem{theorem}{Theorem}[section]
\newtheorem{lemma}[theorem]{Lemma}
\newtheorem{proposition}[theorem]{Proposition}
\newtheorem{corollary}[theorem]{Corollary}

\theoremstyle{definition}
\newtheorem{definition}[theorem]{Definition}
\newtheorem{assumption}[theorem]{Assumption}
\newtheorem{notation}[theorem]{Notation}

\theoremstyle{remark}
\newtheorem{remark}[theorem]{Remark}
\newtheorem{example}[theorem]{Example}

\numberwithin{equation}{part}

\font\tenscr=rsfs10 at 12pt

\newcommand{\lan}{\langle}
\newcommand{\ran}{\rangle}
\newcommand{\vtr}{\vartriangleright}
\newcommand{\coinv}{\mathop{\mathrm{co}}\nolimits}

\newcommand{\grad}{\mathop{\mathrm{grad}}\nolimits}
\newcommand{\diag}{\mathop{\mathrm{diag}}\nolimits}
\newcommand{\Alt}{\mathop{\mathrm{Alt}}\nolimits}
\newcommand{\Sym}{\mathop{\mathrm{Sym}}\nolimits}
\newcommand{\Pf}{\mathop{\mathrm{Pf}}\nolimits}
\newcommand{\id}{\mathop{\mathrm{id}}\nolimits}
\newcommand{\tw}{\mathop{\mathrm{tw}}\nolimits}

\newcommand{\rank}{\mathop{\mathrm{rank}}\nolimits}
\newcommand{\rk}{\mathop{\mathrm{rk}}\nolimits}
\newcommand{\ch}{\mathop{\mathrm{ch}}\nolimits}
\newcommand{\re}{\mathop{\mathrm{Re}}\nolimits}
\newcommand{\im}{\mathop{\mathrm{Im}}\nolimits}
\newcommand{\vep}{\varepsilon}
\newcommand{\Ker}{\mathop{\mathrm{Ker}}\nolimits}
\newcommand{\Hom}{\mathop{\mathrm{Hom}}\nolimits}
\newcommand{\End}{\mathop{\mathrm{End}}\nolimits}
\newcommand{\Alg}{\mathrm{Alg}}
\newcommand{\Coalg}{\mathrm{Coalg}}
\newcommand{\Aut}{\mathop{\mathrm{Aut}}\nolimits}
\newcommand{\AutDlinKalg}{\mathrm{Aut}_{D,\,\mbox{\rm\scriptsize $K$-alg}}}
\newcommand{\AutDlinKTalg}{
\mathrm{Aut}_{D,\,\mbox{\rm\scriptsize $K\otimes_{k}T$-alg}}}
\newcommand{\Spec}{\mathrm{Spec}}
\newcommand{\Pic}{\mathrm{Pic}}
\newcommand{\Gal}{\mathop{\mathrm{Gal}}\nolimits}
\newcommand{\Kdim}{\mathop{\mathrm{Kdim}}\nolimits}
\newcommand{\supp}{\mathop{\mathrm{supp}}\nolimits}
\newcommand{\bF}{\mathbf{F}}
\newcommand{\bG}{\mathbf{G}}
\newcommand{\bH}{\mathbf{H}}
\newcommand{\bT}{\mathbf{T}}
\newcommand{\Ga}{\mathbf{G}_{\mathrm{a}}}
\newcommand{\Gm}{\mathbf{G}_{\mathrm{m}}}
\newcommand{\Gs}{\mathbf{G}_{\mathrm{s}}}
\newcommand{\Gu}{\mathbf{G}_{\mathrm{u}}}
\newcommand{\bAut}{\mathbf{Aut}}
\newcommand{\bAutDlinKalg}{\mathbf{Aut}_{D,\mbox{\rm\scriptsize $K$-alg}}}
\newcommand{\bGL}{\mathbf{GL}}

\newcommand{\bu}{\boldsymbol{u}}
\newcommand{\bv}{\boldsymbol{v}}
\newcommand{\bw}{\boldsymbol{w}}
\newcommand{\by}{\boldsymbol{y}}
\newcommand{\bbc}{\mathbb{C}}
\newcommand{\bbq}{\mathbb{Q}}
\newcommand{\bbr}{\mathbb{R}}
\newcommand{\bbz}{\mathbb{Z}}
\newcommand{\cN}{\mathcal{N}}
\newcommand{\ga}{\mathfrak{a}}
\newcommand{\g}{\mathfrak{g}}
\newcommand{\gp}{\mathfrak{p}}
\newcommand{\gA}{\mathfrak{A}}
\newcommand{\gS}{\mathfrak{S}}
\newcommand{\cA}{\mathcal{A}}
\newcommand{\cC}{\mathcal{C}}
\newcommand{\cH}{\mathcal{H}}
\newcommand{\cI}{\mathcal{I}}
\newcommand{\cM}{\mathcal{M}}

\newcommand{\cR}{\mathcal{R}}
\newcommand{\cS}{\mathcal{S}}
\newcommand{\sD}{\hbox{\tenscr D}}
\newcommand{\sG}{\hbox{\tenscr G}}

\begin{document}

\begin{titlepage}
\vspace*{4cm}
\begin{center}
\Huge Relative invariants, difference equations, and the Picard-Vessiot theory
(revised version 5)
\\
\vspace{3cm}
\LARGE
Katsutoshi Amano \\
\vspace{1cm}

\normalsize
Address: 8-29, Ida-Sammai-cho, Nakahara, Kawasaki, Kanagawa, 211-0037, Japan \\
E-mail: ma-pfybdb-612019@agate.dti.ne.jp \\
URL: http://www.green.dti.ne.jp/amano/index-eng.html
\vspace{\fill}

\large
Revised version of the thesis at the University of Tsukuba, March 2005 \\
The last revision: December 2006
\end{center}
\end{titlepage}

\renewcommand{\thepage}{\roman{page}}

\normalsize

\baselineskip 18pt

\addcontentsline{toc}{part}{\protect\numberline{}{Acknowledgements}}

\begin{center}
\Huge Acknowledgements
\end{center}
\vspace{20pt}

First of all I would like to thank my supervisor Professor T.\ Kimura.
He taught me how to learn mathematics from the beginning when I just started
afresh my life. About four years ago he suggested to study on archimedean
local zeta functions of several variables as my first research task for the
Master's thesis. I can not give enough thanks to his heart-warming
encouragement all over the period of this program.

I thank Professor T.\ Kogiso and Professor K.\ Sugiyama for being close
advisers on prehomogeneous vector spaces.

I also would like to thank Professor M.\ Takeuchi and Professor A.\ Masuoka.
In June or July 2003, I had a chance to know the existence of Professor
Takeuchi's paper \cite{Intro-Takeuchi1989} when I was trying to
understand \cite{Intro-vanderPut-Singer}. Since I once attended to an
introductory part of his lecture on Hopf algebras at 2002,
the paper seemed like just what I wanted. Then I was absorbed in it and
came to think that the results can be extended to involve
\cite{Intro-vanderPut-Singer}. But I could not have developed it in the
presented form without the collaboration with Professor Masuoka.
He suggested the viewpoint from relative Hopf modules as is described in
Section \ref{sec-relative-Hopf-modules}, added many interesting results, and
helped me to complete several proofs.

Finally I would like to thank my parents for supporting my decision
to study mathematics as the lifework.


\bibliographystyle{amsalpha}

\newpage

This is the revised version of my thesis at the University of Tsukuba and is
a compilation of three separate studies.
I shall make an excuse for the awful title ``Relative invariants, difference
equations, and the Picard-Vessiot theory". Though usually a doctoral thesis is
expected to be written in one theme, I could not go this way for institutional
reason. To be accepted, I needed to include at least one topic of a paper
which had been accepted to publication in a refereed journal. But, when
applying for the degree, the paper ``Picard-Vessiot extensions of artinian
simple module algebras" (see Part 3) was not accepted yet (it was accepted
just after I submitted the application for the degree).
Hence I needed to include contents of Part 1, and finally I compiled 
the thesis into three parts to show how my capricious interest had been
changed.

\tableofcontents

\newpage

\renewcommand{\thepage}{\arabic{page}}
\setcounter{page}{1}

\part{An equivariant map from $(SL_{5}\times GL_{4},(\wedge^{2}\bbc^{5})\otimes
\bbc^{4})$ to $(GL_{4},\Sym^{2}(\bbc^{4}))$}
\label{part-concomitant}

\section*{Introduction of Part \ref{part-concomitant}}

In this part, we construct an equivariant polynomial map from $(SL_{5}\times
GL_{4},(\wedge^{2}\bbc^{5})\otimes\bbc^{4})$, the prehomogeneous vector space
of quadruples of quinary alternating forms, to another prehomogeneous
vector space $(GL_{4},\Sym^{2}(\bbc^{4}))$ of quaternary quadratic forms.
The presented result was obtained by Kogiso, Fujigami, and the author
\cite{AFK} to get an expression of the irreducible relative invariant of
the former space explicitly. (Though it is known that the irreducible relative
invariant is a homogeneous polynomial in degree $40$ and it can be computed 
as the determinant of a certain $40\times 40$ matrix by
\cite[\S 4, Proposition 16]{Sato-Kimura}, our purpose was to give a more
effective calculation.)
Then I heard that the equivariant map had further importance as follows.
For a field $k$ in general, the structure of the space
$(SL_{5}\times GL_{4},(\wedge^{2}k^{5})\otimes k^{4})$ has an arithmetic
significance by reason of the correspondence between its non-singular orbits
and isomorphism classes of separable quintic $k$-algebras
(see \cite{Wright-Yukie}).
Kable \cite[Theorem 5.7]{Kable} listed all equivariant polynomial maps
from this space to any other prehomogeneous vector space and showed those maps
could be obtained from two maps; the map presented here is one of the two.
Such equivariant maps seem to be used in \cite{Kable-Yukie1,Kable-Yukie2}
for arithmetic purposes.

Let $\Alt_{n}$ be the set of all skew-symmetric $n\times n$ complex matrices
(i.e. $\Alt_{n} = \{ X\in M_{n}(\bbc)\ |\ ^{t}X=-X\}$).
One sees that the $\bbc$-vector space $(\wedge^{2}\bbc^{5})\otimes
\bbc^{4}$ is isomorphic to $\Alt_{5}^{\oplus 4}$. The space $(SL_{5}\times
GL_{4},(\wedge^{2}\bbc^{5})\otimes\bbc^{4})$ is identified with $(SL_{5}\times
GL_{4},\rho =\Lambda_{2}\otimes\Lambda_{1},\Alt_{5}^{\oplus 4})$ in which the
representation $\rho$ is defined by
\[ \rho (A,B):(X_{1},X_{2},X_{3},X_{4})\longmapsto (AX_{1}{}^{t}A,
AX_{2}{}^{t}A,AX_{3}{}^{t}A,AX_{4}{}^{t}A){}^{t}B \]
for $(X_{1},X_{2},X_{3},X_{4})\in\Alt_{5}^{\oplus 4}$ and $(A,B)\in
SL_{5}\times GL_{4}$.
Our construction of the equivariant map is inspired by the method treated in
\cite[\S 3]{Ochiai}, constructing an equivariant map from $(SL_{5}\times
GL_{3},(\wedge^{2}\bbc^{5})\otimes\bbc^{3})$ to $(GL_{3},\Sym^{2}(\bbc^{3}))$.
Especially a certain $SL_{5}$-equivariant bilinear form
$\beta : \Alt_{5}\times\Alt_{5}\rightarrow\bbc^{5}$, which is
introduced originally in \cite{Gyoja1}, plays an important role.
As in \cite{Gyoja1,Ochiai}, we define $SL_{5}$-invariant polynomials on
$\Alt_{5}^{\oplus 4}$ by
\[ [ijklm](X_{1},X_{2},X_{3},X_{4}):={}^{t}\beta (X_{i},X_{j})X_{k}
\beta (X_{l},X_{m}) \]
for $X_{1},X_{2},X_{3},X_{4}\in\Alt_{5}$ and $i,j,k,l,m\in\{1,2,3,4\}$.
Here each image of $\beta$ is considered as a column vector.
We identify $\Sym^{2}(\bbc^{4})$ with the space of $4\times 4$ symmetric
matrices. Then the equivariant map $\Phi : \Alt_{5}^{\oplus 4}\rightarrow
\Sym^{2}(\bbc^{4})$, $X\mapsto (\varphi_{st}(X))$ will be defined like
\[ \varphi_{st}=\sum_{i,j,k,l,m,i',j',k',l',m'}
c_{stijklmi'j'k'l'm'}[ijklm][i'j'k'l'm'], \]
where the coefficients $c_{stijklmi'j'k'l'm'}$ are determined suitably. 

In \S \ref{sec-beta}, we define the map $\beta$ and the polynomials $[ijklm]$,
and describe some properties of them which are used to obtain the result.
The equivariant map $\Phi$ will be defined in \S \ref{sec-construct} and we
will show the equivariance and the surjectivity of $\Phi$
(Proposition \ref{equivariance} and Theorem \ref{thm-main}). \\

\noindent
{\bf Notations.}
For $\alpha_{1},\alpha_{2},\alpha_{3},\alpha_{4},\varepsilon\in\bbc$,
let $\diag (\alpha_{1},\alpha_{2},\alpha_{3},\alpha_{4})$ and
$E_{\varepsilon}$ be the following matrices:
\[ \diag (\alpha_{1},\alpha_{2},\alpha_{3},\alpha_{4}):=\left(
\begin{array}{llll} \alpha_{1} & 0 & 0 & 0 \\ 0 & \alpha_{2} & 0 & 0 \\
0 & 0 & \alpha_{3} & 0 \\ 0 & 0 & 0 & \alpha_{4} \end{array}\right),\ 
E_{\varepsilon}:=\left(
\begin{array}{llll} 1 & \varepsilon & 0 & 0 \\ 0 & 1 & 0 & 0 \\
0 & 0 & 1 & 0 \\ 0 & 0 & 0 & 1 \end{array}\right). \]
Let $\gS_{4}$ be the fourth symmetric group. In $\gS_{4}$, a
transposition between $i$ and $j$ is denoted by $(i\ j)$. One sees each
permutation $\sigma\in\gS_{4}$ is considered as the $4\times 4$
matrix such that its $(i,j)$-element is $1$ or $0$ with respect to
$i=\sigma (j)$ or not. So we may apply for regarding one as the other.

\section{$SL_{5}$-invariant polynomials on $\Alt_{5}^{\oplus 4}$}
\label{sec-beta}

In the beginning, we define a certain $SL_{5}$-equivariant map $\beta :
\Alt_{5}\times\Alt_{5}\rightarrow\bbc^{5}$ which is used in
\cite{Gyoja1,Ochiai}. Let $\Pf$ be the Pfaffian on $\Alt_{4}$.
For $X\in\Alt_{5}$ and $i=1,\cdots ,5$, let $X^{(i)}$ denote the matrix in
$\Alt_{4}$ which is obtained by deleting $i$-th row and $i$-th column from $X$.
For $X=(x_{ij}),Y=(y_{ij})\in\Alt_{5}$, $\beta (X,Y)$ is defined by
\begin{eqnarray*} \beta (X,Y) & := & \left(
\begin{array}{l}
\Pf (X^{(1)}+Y^{(1)})-\Pf (X^{(1)})-\Pf (Y^{(1)}) \\
-(\Pf (X^{(2)}+Y^{(2)})-\Pf (X^{(2)})-\Pf (Y^{(2)})) \\
\Pf (X^{(3)}+Y^{(3)})-\Pf (X^{(3)})-\Pf (Y^{(3)}) \\
-(\Pf (X^{(4)}+Y^{(4)})-\Pf (X^{(4)})-\Pf (Y^{(4)})) \\
\Pf (X^{(5)}+Y^{(5)})-\Pf (X^{(5)})-\Pf (Y^{(5)})
\end{array}\right) \\
 & = & \left(
\begin{array}{l}
x_{23}y_{45}-x_{24}y_{35}+x_{25}y_{34}+y_{23}x_{45}-y_{24}x_{35}+y_{25}x_{34}\\
x_{34}y_{51}-x_{35}y_{41}+x_{31}y_{45}+y_{34}x_{51}-y_{35}x_{41}+y_{31}x_{45}\\
x_{45}y_{12}-x_{41}y_{52}+x_{42}y_{51}+y_{45}x_{12}-y_{41}x_{52}+y_{42}x_{51}\\
x_{51}y_{23}-x_{52}y_{13}+x_{53}y_{12}+y_{51}x_{23}-y_{52}x_{13}+y_{53}x_{12}\\
x_{12}y_{34}-x_{13}y_{24}+x_{14}y_{23}+x_{12}y_{34}-x_{13}y_{24}+x_{14}y_{23}
\end{array}\right).
\end{eqnarray*}
Then, for $i,j,k,l,m\in\{1,2,3,4\}$, we define a polynomial
$[ijklm]$ on $\Alt_{5}^{\oplus 4}$ by
\[ [ijklm](X_{1},X_{2},X_{3},X_{4}):={}^{t}\beta (X_{i},X_{j})X_{k}
\beta (X_{l},X_{m}) \]
for $X_{1},X_{2},X_{3},X_{4}\in\Alt_{5}$.
They satisfy the following lemmas:

\begin{lemma}[{\cite[\S 2, Lemma]{Gyoja1}}]
\label{sl5invariant}
For all $i,j,k,l,m\in\{1,2,3,4\}$, the polynomial $[ijklm]$ is invariant
with respect to $SL_{5}$, i.e.
\[ [ijklm](AX_{1}{}^{t}A,AX_{2}{}^{t}A,AX_{3}{}^{t}A,AX_{4}{}^{t}A)
=[ijklm](X_{1},X_{2},X_{3},X_{4}) \]
for all $A\in SL_{5}$.
\end{lemma}

\begin{lemma}[{\cite[\S 2, (4)]{Gyoja1}}]
\label{threekinds}
If there are only one or two kinds of numbers among $\{i,j,k,l,m\}$, then
$[ijklm]=0$.
\end{lemma}

\begin{lemma}[{\cite[Lemma 3.1]{Ochiai}}]
\label{relations}
For each $i,j,k,l,m\in\{1,2,3,4\}$,
\[ \begin{array}{ll}
\mbox{\rm(i)}\ [ijklm]=[jiklm],\ [ijklm]=[ijkml], &
\mbox{\rm(ii)}\ [ijklm]=-[lmkij], \\
\mbox{\rm(iii)}\ [ijklm]+[jkilm]+[kijlm]=0, &
\mbox{\rm(iv)}\ [iiklm]=-2[kiilm], \\
\mbox{\rm(v)}\ [iikli]=-[iilki]=[iklii]=-[ilkii], &
\mbox{\rm(vi)}\ [iiilm]=0,\ [ijkij]=0. 
\end{array} \]
\end{lemma}

Finally in this section, we consider the action of $GL_{4}$ on $[ijklm]$.
$GL_{4}$ is generated by the following three types of matrices:
$\diag (\alpha_{1},\alpha_{2},\alpha_{3},\alpha_{4})$, permutation matrices,
and $E_{\varepsilon}$. Thus we only need to think on these types.
For $B\in GL_{4}$ and $P$ a polynomial on $\Alt_{5}^{\oplus 4}$, let
$P^{B}$ denote the polynomial such that $P^{B}(X)=P(X{}^{t}B)$.
Diagonal matrices $D=\diag(\alpha_{1},\alpha_{2},\alpha_{3},\alpha_{4})$ and
$\sigma\in\gS_{4}$ act on $[ijklm]$ by
\[ \begin{array}{lll}
[ijklm]^{D} & = & \alpha_{i}\alpha_{j}\alpha_{k}\alpha_{l}\alpha_{m}[jiklm], \\
{[ijklm]^{\sigma}} & = &
[\sigma^{-1}(i)\sigma^{-1}(j)\sigma^{-1}(k)\sigma^{-1}(l)\sigma^{-1}(m)].
\end{array} \]
Since $[ijklm]$ are multilinear forms, we see, for $i,j,k,l,m\in\{2,3,4\}$,
\[ \begin{array}{lll}
{[ijklm]^{E_{\varepsilon}}} & = & [ijklm], \\
{[1ijkl]^{E_{\varepsilon}}} & = & [1ijkl]+\varepsilon [2ijkl], \\
{[11ijk]^{E_{\varepsilon}}} & = & [11ijk]+2\varepsilon [12ijk]
+\varepsilon^{2}[22ijk], \\
{[11ij1]^{E_{\varepsilon}}} & = & [11ij1]+\varepsilon(2[12ij1]+[11ij2])\\
 & & \ +\varepsilon^{2}(2[12ij2]+[22ij1])+\varepsilon^{3}[22ij2],
\ \mbox{etc.}
\end{array} \]

\section{Construction of the equivariant map}
\label{sec-construct}

Our first objective is to define a map $\Phi : \Alt_{5}^{\oplus 4}\rightarrow
\Sym^{2}(\bbc^{4})$, $X\mapsto (\varphi_{st}(X))$, where each
$\varphi_{st}$ is written like 
\[ \varphi_{st}=\sum c_{stijklmi'j'k'l'm'}[ijklm][i'j'k'l'm']. \]
\begin{proposition}
\label{equivariance}
There exists a polynomial map $\Phi : \Alt_{5}^{\oplus 4}\rightarrow
\Sym^{2}(\bbc^{4})$ such that
\[ \Phi (\rho(A,B)X)=(\det B)^{2}B\Phi(X)\,^{t}B \]
for $X\in\Alt_{5}^{\oplus 4}$ and $(A,B)\in SL_{5}\times GL_{4}$.
\end{proposition}

First we observe that $\Phi$ should be determined uniquely from $\varphi_{11}$
and $\varphi_{12}$ so that $\Phi(X)$ is equivariant
with respect to the action of $\gS_{4}$:
\begin{equation}
\label{others}
\varphi_{ss}=\varphi_{11}^{(1\ s)}\quad (s=1,\dotsc,4),\qquad
\varphi_{\sigma^{-1}(1)\sigma^{-1}(2)}=\varphi_{12}^{\sigma}\quad
(\sigma\in\gS_{4}).
\end{equation}
Furthermore, $\varphi_{12}$ should also be determined from $\varphi_{11}$.
To obtain $\Phi(\rho(A,E_{\varepsilon})X)
=E_{\varepsilon}\Phi(X)\,^{t}E_{\varepsilon}$ ($A\in SL_{5}$), the polynomials
$\varphi_{st}$ should satisfy at least the following:
\begin{eqnarray}
\label{11}
\varphi_{11}^{E_{\varepsilon}} & = & \varphi_{11}
+2\varepsilon\varphi_{12}+\varepsilon^{2}\varphi_{22}, \\
\label{22}
\varphi_{22}^{E_{\varepsilon}} & = & \varphi_{22}, \\
\label{33}
\varphi_{33}^{E_{\varepsilon}} & = & \varphi_{33}, \\
\label{13}
\varphi_{13}^{E_{\varepsilon}} & = & \varphi_{13}+\varepsilon\varphi_{23}, \\
\label{34}
\varphi_{34}^{E_{\varepsilon}} & = & \varphi_{34}.
\end{eqnarray}
If we obtain $\varphi_{11}$, then $\varphi_{22}=\varphi_{11}^{(1\ 2)}$ and
$\varphi_{12}$ will be determined from (\ref{11}).

Considering the action of diagonal matrices for $\Phi(X)$, we start by
assuming that each term $[ijklm][i'j'k'l'm']$ in $\varphi_{11}$ is constructed
by the following numbers:
\[ \{i,j,k,l,m,i',j',k',l',m'\}=\{1,1,1,1,2,2,3,3,4,4\}. \]
But from Lemma \ref{threekinds} and Lemma \ref{relations}, we need not to
think on the all combinations of the above numbers.
By choosing combinations and using the method of indeterminate
coefficients, it is possible to determine the polynomial $\varphi_{11}$ so that
the equations from (\ref{11}) to (\ref{34}) are satisfied. Indeed,
we conclude that the following definitions are suitable:
\begin{eqnarray*}
\varphi_{11} & := & 160[31114](3[24132]-2[21342]-2[23412]) \\
 & & +160[41112](3[32143]-2[34213]-2[31423]) \\
 & & +160[21113](3[43124]-2[41234]-2[42314]) \\
 & & +50([11233][11244]+[11322][11344]+[11422][11433]) \\
 & & -288([13241]^{2}+[14321]^{2}+[12431]^{2}) \\
 & & +224([13241][14321]+[14321][12431]+[12431][13241]),
\end{eqnarray*}
\begin{eqnarray*}
\varphi_{12} & := & 400[31114][32224] \\
 & & -100([21113][22344]+[21114][22433]) \\
 & & -100([12223][11344]+[12224][11433]) \\
 & & +20[11422](4[31423]-[34213]-[32143]) \\
 & & +20[11322](4[41324]-[43214]-[42134]) \\
 & & -25([22144][11233]+[11244][22133]) \\
 & & +368[13241][23142] \\
 & & +112([13241]([21342]+[23412])+[23142]([12341]+[13421])) \\
 & & +192([14321][23412]+[13421][24312]) \\
 & & -208([14321][21342]+[12431][23412]).
\end{eqnarray*}
These polynomials satisfy the following properties:
\begin{enumerate}
\renewcommand{\labelenumi}{(\roman{enumi})}
\item If $\sigma\in\gS_{4}$ and $\sigma(1)=1$, then
$\varphi_{11}^{\sigma}=\varphi_{11}$,
\item If $\sigma\in\gS_{4}$ and $\{\sigma(1),\sigma(2)\}=\{1,2\}$,
then $\varphi_{12}^{\sigma}=\varphi_{12}$.
\end{enumerate}
Then we define the map $\Phi : \Alt_{5}^{\oplus 4}\rightarrow
\Sym^{2}(\bbc^{4})$, $X\mapsto (\varphi_{st}(X))$ so that (\ref{others}) is
satisfied; the well-definedness follows from (i), (ii).
It is easily seen that $\varphi_{st}=\varphi_{ts}$ and $\varphi_{st}^{\sigma}
=\varphi_{\sigma^{-1}(s)\sigma^{-1}(t)}$ for all $\sigma\in\gS_{4}$.

\begin{proof}[{Proof of Proposition \ref{equivariance}}]
Let $D=\diag(\alpha_{1},\alpha_{2},\alpha_{3},\alpha_{4})$ and let $A$ be an
arbitrary element of $SL_{5}$. Since $\varphi_{st}^{D}
=(\alpha_{1}\alpha_{2}\alpha_{3}\alpha_{4})^{2}
\alpha_{s}\alpha_{t}\varphi_{st}$ for all $s,t\in\{1,2,3,4\}$, and each
$\varphi_{st}$ is invariant with respect to $SL_{5}$, we have
\[ \Phi(\rho(A,D)X)=(\det D)^{2}D\Phi(X)\,^{t}D. \]
By the definition, it follows
\[ \Phi(\rho(A,\sigma)X)=(\varphi_{\sigma^{-1}(s)\sigma^{-1}(t)}(X))
=\sigma\Phi(X)\,^{t}\sigma \]
for all $\sigma\in\gS_{4}$.

The rest of the proof is to show $\Phi(\rho(A,E_{\varepsilon})X)
=E_{\varepsilon}\Phi(X)\,^{t}E_{\varepsilon}$, i.e.
\begin{itemize}
\item $\varphi_{11}^{E_{\varepsilon}}=\varphi_{11}+2\varepsilon\varphi_{12}
{}+\varepsilon^{2}\varphi_{22}$,
\item $\varphi_{1t}^{E_{\varepsilon}}=\varphi_{t1}^{E_{\varepsilon}}
=\varphi_{1t}+\varepsilon \varphi_{2t}$ for $t=2,3,4$,
\item $\varphi_{st}^{E_{\varepsilon}}
=\varphi_{ts}^{E_{\varepsilon}}=\varphi_{st}$ for $s,t=2,3,4$.
\end{itemize}
Recall that we defined $\varphi_{st}$ to satisfy the equations from (\ref{11})
to (\ref{34}) (in fact, they are shown directly).
By $E_{\varepsilon}^{2}=E_{2\varepsilon}$ and (\ref{11}), we have
\[ \varphi_{11}^{E_{\varepsilon}^{2}}=\varphi_{11}+4\varepsilon\varphi_{12}
+4\varepsilon^{2}\varphi_{22}. \]
On the other hand, by (\ref{11}) and (\ref{22}),
\begin{eqnarray*}
\varphi_{11}^{E_{\varepsilon}^{2}} & = & (\varphi_{11}+2\varepsilon\varphi_{12}
{}+\varepsilon^{2}\varphi_{22})^{E_{\varepsilon}} \\ & = &
\varphi_{11}^{E_{\varepsilon}}+2\varepsilon\varphi_{12}^{E_{\varepsilon}}
{}+\varepsilon^{2}\varphi_{22}^{E_{\varepsilon}} \\
 & = & \varphi_{11}+2\varepsilon\varphi_{12}
+2\varepsilon\varphi_{12}^{E_{\varepsilon}}+2\varepsilon^{2}\varphi_{22}.
\end{eqnarray*}
Therefore $\varphi_{12}^{E_{\varepsilon}}
=\varphi_{12}+\varepsilon\varphi_{22}$. Similarly by (\ref{13}),
\begin{eqnarray*}
\varphi_{13}^{E_{\varepsilon}^{2}} & = & \varphi_{13}+2\varepsilon\varphi_{23}
\\
& = & \varphi_{13}+\varepsilon\varphi_{23}
{}+\varepsilon{\varphi_{23}^{E_{\varepsilon}}}.
\end{eqnarray*}
Hence $\varphi_{23}^{E_{\varepsilon}}=\varphi_{23}$. By (\ref{13})
and $E_{\varepsilon}(3\ 4)=(3\ 4)E_{\varepsilon}$, we have
\[ \varphi_{14}^{E_{\varepsilon}}=\varphi_{13}^{(3\ 4)E_{\varepsilon}}
=\varphi_{13}^{E_{\varepsilon}(3\ 4)}
=(\varphi_{13}^{E_{\varepsilon}})^{(3\ 4)}=\varphi_{14}
+\varepsilon\varphi_{24}. \]
Similarly by (\ref{33}), we have 
\[ \varphi_{44}^{E_{\varepsilon}}=\varphi_{33}^{E_{\varepsilon}(3\ 4)}
=\varphi_{33}^{(3\ 4)}=\varphi_{44}. \]
Now the proof is completed.
\end{proof}

To prove that $\Phi$ is surjective, we only need to find five points in
$\Alt_{5}^{\oplus 4}$ such that each image has rank $0,1,2,3,4$. For
\begin{center}
$\begin{array}{llllll}
X_{01} & = & \left(
\begin{array}{lllll}
0 & 1 & 0 & 0 & 0 \\
-1 & 0 & 0 & 0 & 0 \\
0 & 0 & 0 & 1 & 0 \\
0 & 0 & -1 & 0 & 0 \\
0 & 0 & 0 & 0 & 0
\end{array}\right), &
X_{02} & = & \left(
\begin{array}{lllll}
0 & 0 & 0 & 0 & 0 \\
0 & 0 & 1 & 0 & 0 \\
0 & -1 & 0 & 0 & 0 \\
0 & 0 & 0 & 0 & 1 \\
0 & 0 & 0 & -1 & 0
\end{array}\right), \vspace{5pt} \\
X_{03} & = & \left(
\begin{array}{lllll}
0 & 0 & 1 & 0 & 0 \\
0 & 0 & 0 & 0 & 1 \\
-1 & 0 & 0 & 0 & 0 \\
0 & 0 & 0 & 0 & 0 \\
0 & -1 & 0 & 0 & 0
\end{array}\right), &
X_{04} & = & \left(
\begin{array}{lllll}
0 & 0 & 0 & 0 & 0 \\
0 & 0 & 0 & 1 & 0 \\
0 & 0 & 0 & 0 & 1 \\
0 & -1 & 0 & 0 & 0 \\
0 & 0 & -1 & 0 & 0
\end{array}\right), \vspace{5pt} \\
Y_{01} & = & \left(
\begin{array}{lllll}
0 & 1 & 1 & 0 & 0 \\
-1 & 0 & 0 & 0 & 0 \\
-1 & 0 & 0 & 0 & 0 \\
0 & 0 & 0 & 0 & 0 \\
0 & 0 & 0 & 0 & 0
\end{array}\right), &
Y_{02} & = & \left(
\begin{array}{lllll}
0 & 0 & 0 & 0 & 0 \\
0 & 0 & 0 & 0 & 0 \\
0 & 0 & 0 & 0 & 0 \\
0 & 0 & 0 & 0 & 1 \\
0 & 0 & 0 & -1 & 0
\end{array}\right), \vspace{5pt} \\
Y_{03} & = & \left(
\begin{array}{lllll}
0 & 0 & 0 & 0 & 0 \\
0 & 0 & 0 & 0 & 1 \\
0 & 0 & 0 & 0 & 0 \\
0 & 0 & 0 & 0 & 0 \\
0 & -1 & 0 & 0 & 0
\end{array}\right), & & &
\end{array}$
\end{center}
we have
\begin{center}
$\begin{array}{llll}
\Phi(X_{01},X_{02},X_{03},X_{04}) & = & \left(
\begin{array}{llll}
0 & 0 & -720 & 0 \\
0 & -480 & 0 & 0 \\
-720 & 0 & 0 & 0 \\
0 & 0 & 0 & -288
\end{array}\right) & (\mbox{rank $4$}), \vspace{5pt} \\
\Phi(Y_{01},X_{02},X_{03},X_{04}) & = & \left(
\begin{array}{llll}
-192 & 0 & -192 & -96 \\
0 & -480 & 0 & 0 \\
-192 & 0 & -192 & -96 \\
-96 & 0 & -96 & -288
\end{array}\right) & (\mbox{rank $3$}), \vspace{5pt} \\
\Phi(Y_{01},Y_{02},X_{03},X_{04}) & = & \left(
\begin{array}{llll}
-192 & 0 & -192 & -96 \\
0 & 0 & 0 & 0 \\
-192 & 0 & -192 & -96 \\
-96 & 0 & -96 & -288
\end{array}\right) & (\mbox{rank $2$}), \vspace{5pt} \\
\Phi(Y_{01},Y_{02},Y_{03},X_{04}) & = & \left(
\begin{array}{llll}
-192 & 0 & 0 & 0 \\
0 & 0 & 0 & 0 \\
0 & 0 & 0 & 0 \\
0 & 0 & 0 & 0
\end{array}\right) & (\mbox{rank $1$}), \vspace{5pt} \\
\Phi(0,0,0,0) & = & 0 & (\mbox{rank $0$}).
\end{array}$
\end{center}
Therefore $\Phi$ is surjective and especially $\det\Phi(X)\neq 0$ as a
polynomial. This fact and Proposition \ref{equivariance} implies that
$\det\Phi(X)$ is the relative invariant in degree $40$.

\begin{theorem}
\label{thm-main}
{\rm (i)} The map $\Phi :\Alt_{5}^{\oplus 4}\rightarrow \Sym^{2}(\bbc^{4})$ is
surjective.

{\rm (ii)} $f(X)=\det\Phi(X)$ is the irreducible relative invariant of
the prehomogeneous vector space $(SL_{5}\times GL_{4},\Lambda_{2}\otimes
\Lambda_{1},\Alt_{5}^{\oplus 4})$ in degree $40$ corresponding to the rational
character $(\det B)^{4}$.
\end{theorem}


\bibliographystyle{amsalpha}

\newpage

\setcounter{section}{0}

\part{Archimedean local zeta functions which satisfy $\Gm$-primitive difference equations}
\label{part-localzeta}

\section*{Introduction of Part \ref{part-localzeta}}

Let $K$ be $\bbc$ or $\bbr$, $(G,\rho,V)$ a reductive prehomogeneous vector
space defined over $K$, and $V_{K}$ the set of $K$-rational points of $V$.
Let $P_{1}(x),\cdots,P_{r}(x)$ be the basic relative invariants of $(G,\rho,V)$
over $K$. For a Schwartz function $\Phi(x)$ on $V_{K}$ and
$s=(s_{1},\cdots,s_{r})\in\mathbb{C}^{r}$, the integral
\[ Z_{K}(s,\Phi)=\int_{V_{K}}\prod_{i=1}^{r}|P_{i}(x)|_{K}^{s_{i}}\Phi(x)dx \]
is called the {\em archimedean local zeta function} associated with
$(G,\rho,V)$. When we take $\Phi(x)$ as
\[ \Phi(x)=\left\{\begin{array}{ll}
\exp(-2\pi\,x{}^{t}\bar{x}) & (K=\bbc) \\
\exp(-\pi\,x{}^{t}x) & (K=\bbr), \end{array}\right. \]
$Z_{K}(s,\Phi)$ is denoted by $Z_{K}(s)$ simply.

In the case $r=1$, Igusa suggested
in \cite[\S 3, Remark]{Igusa1986} and proved in
\cite[Chapter 6]{Igusa2000} the following theorem: \\

\noindent
{\bf Theorem (Igusa).}
{\it
Let $P$ be the basic relative invariant of $(G,\rho,V)$, assuming $r=1$.
Let $d=\deg P$ and $b(s)=c\prod_{j=1}^{d}(s+\alpha_{j})$ the $b$-function
of $P(x)$.

{\rm (1)} When $K=\bbc$,
\[ Z_{\bbc}(s)=((2\pi)^{-d}c)^{s}
\prod_{j=1}^{d}\frac{\Gamma(s+\alpha_{j})}{\Gamma(\alpha_{j})}. \]

{\rm (2)} When $K=\bbr$ and when
every term of $P(x)$ is a multilinear form,
\[ Z_{\bbr}(s)=(\pi^{-d}c)^{\frac{s}{2}}
\prod_{j=1}^{d}\frac{\Gamma((s+\alpha_{j})/2)}{\Gamma(\alpha_{j}/2)}. \]
}

In this part, we extend this theorem to several
variables ($r\geq 1$). Though basically our proof presented here
is an easy modification of the proof given in \cite{Igusa2000},
a careful treatment of the Ore-Sato theorem (see Section \ref{sec-Ore-Sato})
is needed. The most important point of the proof is the fact that
$Z_{\bbc}(s)$ and $Z_{\bbr}(2s)$ satisfies a difference equation in a certain
type, called {\em $\Gm$-primitive}, or {\em hypergeometric}, which is written
by the $b$-functions.
The proof can be divided into two steps. The first step is to
characterize a desired solution of such a difference equation, written as
a product of an exponential function and gamma functions.
The next step is to prove that the characterization can be adapted to
$Z_{\bbc}(s)$ and $Z_{\bbr}(2s)$. To obtain a suitable difference equation
(especially for $Z_{\bbr}(2s)$), very delicate facts which are seen
in the proof of the Ore-Sato theorem are needed. For this reason,
we include a detailed proof of the theorem in Section \ref{sec-Ore-Sato}.
The main results will be described in Section \ref{sec-complex-LZ},
\ref{sec-real-LZ} (Theorem \ref{thm-complex-LZ} and Theorem \ref{thm-real-LZ}).

\section{The Ore-Sato theorem}
\label{sec-Ore-Sato}

Let $k$ be a field of zero characteristic and $k(s)=k(s_{1},\dotsc,s_{r})$
the rational function field of $r$ variables.
Let $\Xi$ be the free abelian group of rank $r$
($\Xi\simeq\bbz^{r}$) and $\tau_{1},\dotsc,\tau_{r}$ a basis of $\Xi$.
Then $\Xi$ acts on $k(s)$ as $k$-algebra automorphisms by $\tau_{i}f(s)
=f(s+e_{i})$ ($f(s)\in k(s)$) where $e_{1}=(1,0,\dotsc,0),\dotsc,e_{r}
=(0,\dotsc,0,1)$, the canonical basis of $k^{r}$. Let $k\Xi$ be the group
algebra of $\Xi$ over $k$ and $k(s)\#k\Xi$ the ring of linear difference
operators, which is $k(s)\otimes_{k}k\Xi$ with the semi-direct product:
$(f_{1}\otimes g)\cdot(f_{2}\otimes h)=f_{1}(gf_{2})\otimes gh$ ($g,h\in\Xi$).
We say that a $k(s)\#k\Xi$-module $V$ is {\em $\Gm$-primitive}, or
{\em hypergeometric}, iff $\dim_{k(s)}V=1$.
(We will see in Part \ref{part-PV-theory} that the
Picard-Vessiot group scheme of such a $k(s)\#k\Xi$-module is a closed subgroup
scheme of $\Gm$.) Let $k(s)^{\times}=k(s)\setminus\{0\}$.
For a fixed $k(s)$-basis $v$ of a $\Gm$-primitive $k(s)\#k\Xi$-module $V$,
we have an associated map $b_{v} : \Xi\rightarrow k(s)^{\times}$, $g\mapsto
b_{g,v}(s)$ defined by $gv=b_{g,v}(s)v$. Since $b_{1,v}(s)=1$ and $b_{gh,v}(s)
=(gb_{h,v}(s))b_{g,v}(s)$ for all $g,h\in \Xi$, $b_{v}$ is in the set
$Z^{1}(\Xi;k(s)^{\times})$ of $1$-cocycles. Let $v'$ be another $k(s)$-basis
of $V$. Then there exists an $f\in k(s)^{\times}$ such that $v'=f(s)v$.
It follows $b_{g,v'}(s)=(gf(s))f(s)^{-1}b_{g,v}(s)$ for all $g\in\Xi$; this
implies that both of $b_{v'}$ and $b_{v}$ define the same cohomology class
in $H^{1}(\Xi;k(s)^{\times})$ since the map $g\mapsto (gf(s))f(s)^{-1}$ is
in the set $B^{1}(\Xi;k(s)^{\times})$ of $1$-coboundaries.
Thus $\Gm$-primitive $k(s)\#k\Xi$-modules are classified by
$H^{1}(\Xi;k(s)^{\times})$.

An explicit description of $Z^{1}(\Xi;k(s)^{\times})$ is given by
a result called {\em the Ore-Sato theorem}, which was first obtained by
Ore \cite{Ore} for the case $r=2$ and by Sato \cite{Sato-Shintani-Muro}
for arbitrary $r$. Detailed proofs are also seen in
\cite[\S 1.1]{Loeser-Sabbah} and \cite[\S 1]{Gelfand-Graev-Retakh}.
The purpose of this section is to introduce the theorem for later use.
Since some delicate facts such as Corollary \ref{cor-Ore-Sato} are important
to us, we follow carefully the discussion given in
\cite[Appendix]{Sato-Shintani-Muro}. (Thus our statement of the theorem may be
verbose according to the interest of the reader. For a more elegant
description of $H^{1}(\Xi;\bbc(s)^{\times})$,
\cite[Proposition 1.1.4]{Loeser-Sabbah} is recommended.
See Remark \ref{rem-Ore-Sato}.)

$k(s)^{\times}$ has a natural $\bbz\Xi$-module structure as follows:
\[ \left(\sum_{i}n_{i}g_{i}\right)f(s)=\prod_{i}(g_{i}f(s))^{n_{i}}
\qquad (n_{i}\in\bbz,\quad g_{i}\in\Xi). \]
We easily see $k(s)^{\times}\simeq k^{\times}\times k(s)^{\times}/k^{\times}$
as $\bbz\Xi$-modules. Moreover, by decomposing to irreducible polynomials,
we have the following $\bbz\Xi$-module isomorphisms:
\[ k(s)^{\times}\xrightarrow{\sim} k^{\times}\oplus\bigoplus_{f}\bbz\Xi f
\xrightarrow{\sim} k^{\times}\oplus\bigoplus_{f}\bbz(\Xi/\Xi_{f}), \]
where $f$ runs over a set of irreducible polynomials which are not
translated into one another up to constant multiple by the action of $\bbz\Xi$,
and $\Xi_{f}:=\{g\in\Xi\;|\;gf=\mbox{\rm(const.)}f\}=\{g\in\Xi\;|\;gf=f\}$.
$\bbz(\Xi/\Xi_{f})$ is the group ring of $\Xi/\Xi_{f}$ over $\bbz$.

Since $\Xi$ is a finitely generated group, we have
\[ Z^{1}(\Xi;k(s)^{\times})\simeq Z^{1}(\Xi;k^{\times})\oplus
\bigoplus_{f}Z^{1}(\Xi;\bbz(\Xi/\Xi_{f})). \]
One sees that $Z^{1}(\Xi;k^{\times})$ is equal to the character group
$\Hom(\Xi,k^{\times})$. The structure of each $Z^{1}(\Xi;\bbz(\Xi/\Xi_{f}))$ is
described by the following two lemmas.

\begin{lemma}
\label{lem-Ore-Sato1}
Let $f$ be an irreducible polynomial in $k[s]=k[s_{1},\dotsc,s_{r}]$.

{\rm (i)} $\Xi/\Xi_{f}$ is a free abelian group of $\rank\Xi/\Xi_{f} > 0$.

{\rm (ii)} Take an arbitrary $\alpha\in Z^{1}(\Xi;\bbz(\Xi/\Xi_{f}))$.
Then $\alpha(g)=0$ for all $g\in\Xi_{f}$.

{\rm (iii)} If $\rank\Xi/\Xi_{f}\geq 2$, then $H^{1}(\Xi;\bbz(\Xi/\Xi_{f}))
=0$.
\end{lemma}
\begin{proof}
(i) Since $f$ is not a constant, $\Xi_{f}\neq\Xi$. Suppose $g\in\Xi$ and
$g^{n}\in\Xi_{f}$ for some positive integer $n$. There is an $m\in\bbz^{r}$
such that $gf(s)=f(s+m)$. Consider $P(t):=f(s+tnm)-f(s)$ as a polynomial
in $k(s)[t]$. Since $P(l)=g^{ln}f(s)-f(s)=0$ for all $l\in\bbz$ and since
$k(s)$ is an infinite field, we have $P\equiv 0$. Hence $0=P(1/n)=gf(s)-f(s)$,
and so $g\in\Xi_{f}$.

(ii) Take a $g'\in\Xi$ which is not in $\Xi_{f}$.
Let $\tilde{g}'$ be the image of $g'$ in $\Xi/\Xi_{f}$. Then for all $g\in
\Xi_{f}$, we have
\[ \alpha(gg')=\alpha(g')+\alpha(g)=\tilde{g}'\alpha(g)+\alpha(g'). \]
Thus $(\tilde{g}'-1)\alpha(g)=0$. Since $\bbz(\Xi/\Xi_{f})$ is an integral
domain by part (i), we have $\alpha(g)=0$.

(iii) By part (ii), we have $H^{1}(\Xi;\bbz(\Xi/\Xi_{f}))\simeq
H^{1}(\Xi/\Xi_{f};\bbz(\Xi/\Xi_{f}))$. Let $\lambda_{1},\dotsc,\lambda_{l}$
be a basis of $\Xi/\Xi_{f}$. For all $\alpha\in Z^{1}(\Xi/\Xi_{f};
\bbz(\Xi/\Xi_{f}))$, we have
\[ \alpha(\lambda_{i}\lambda_{j})
=\lambda_{i}\alpha(\lambda_{j})+\alpha(\lambda_{i})
=\lambda_{j}\alpha(\lambda_{i})+\alpha(\lambda_{j}) \]
and hence $(\lambda_{i}-1)\alpha(\lambda_{j})=(\lambda_{j}-1)
\alpha(\lambda_{i})$ for $i,j=1,\dotsc,l$. Since $\bbz(\Xi/\Xi_{f})$
is the Laurent polynomial ring, there exists an $a\in\bbz(\Xi/\Xi_{f})$ 
such that
\[ \alpha(\lambda_{i})=(\lambda_{i}-1)a\qquad (i=1,\dotsc,l). \]
Moreover, since
\[ 0=\alpha(1)=\alpha(\lambda_{i}^{-1}\lambda_{i})
=\lambda_{i}^{-1}\alpha(\lambda_{i})+\alpha(\lambda_{i}^{-1})
=(1-\lambda_{i}^{-1})a+\alpha(\lambda_{i}^{-1}), \]
we have $\alpha(\lambda_{i}^{-1})=(\lambda_{i}^{-1}-1)a$. Therefore
\[ \alpha(g)=ga-a\quad (g\in\Xi/\Xi_{f}), \]
concluding $\alpha\in B^{1}(\Xi/\Xi_{f};\bbz(\Xi/\Xi_{f}))$, follows by
induction. Indeed, if $\alpha(g')=g'a-a$ for a $g'\in\Xi/\Xi_{f}$, then
\[ \begin{array}{l}
\alpha(\lambda_{i}g')=\lambda_{i}(g'a-a)+(\lambda_{i}-1)a=\lambda_{i}g'a-a,
\vspace{1ex} \\
\alpha(\lambda_{i}^{-1}g')=\lambda_{i}^{-1}(g'a-a)+(\lambda_{i}^{-1}-1)a
=\lambda_{i}^{-1}g'a-a
\end{array} \]
for $i=1,\dotsc,l$.
\end{proof}

\begin{lemma}
\label{lem-Ore-Sato2}
Let $f$ be an irreducible polynomial in $k[s]$.
Suppose that $\rank\Xi/\Xi_{f}=1$ and let $\lambda$ be a basis of
$\Xi/\Xi_{f}$.

{\rm (i)} There exist a linear form $\mu(s)=n_{1}s_{1}+\dotsb +n_{r}s_{r}$
($n_{1},\dotsc,n_{r}\in\bbz$) and an irreducible polynomial $h(t)\in
k[t]$ of one variable $t$, such that $f(s)=h(\mu(s))$ and $\lambda f(s)
=h(\mu(s)+1)$. Moreover, there exists an $m\in\bbz^{r}$ such that
$\mu(m)=1$, i.e.\ the greatest common divisor of non-zero coefficients of
$\mu$ is $1$.

{\rm (ii)} For every $\alpha\in Z^{1}(\Xi;\bbz(\Xi/\Xi_{f}))$, there is an
$\eta\in\bbz(\Xi/\Xi_{f})$ such that
\[ \alpha(\tau_{1}^{m_{1}}\cdots\tau_{r}^{m_{r}})
=\left\{\begin{array}{cl}
\eta{\displaystyle\sum_{\nu=0}^{\mu(m)-1}}\lambda^{\nu} & (\mu(m)\geq 1),
\vspace{1ex} \\
0 & (\mu(m)=0), \vspace{1ex} \\
-\eta{\displaystyle\sum_{\nu=\mu(m)}^{-1}}\lambda^{\nu} & (\mu(m)\leq -1),
\end{array}\right. \]
for all $m=(m_{1},\dotsc,m_{r})\in\bbz^{r}$.

{\rm (iii)} $H^{1}(\Xi;\bbz(\Xi/\Xi_{f}))\simeq\bbz$.
\end{lemma}
\begin{proof}
(i) Let $\lambda_{1},\dotsc,\lambda_{r}$ be a basis of $\Xi$ such that
the image of $\lambda_{1}$ in $\Xi/\Xi_{f}$ is $\lambda$ and
$\lambda_{2},\dotsc,\lambda_{r}\in\Xi_{f}$. Then there
exists an invertible $r\times r$ matrix $(n_{ij})_{i,j}$ in $GL_{r}(\bbz)$
such that $\tau_{i}=\lambda_{1}^{n_{i1}}\cdots\lambda_{r}^{n_{ir}}$
for $i=1,\dotsc,r$. Put $\mu(s):=n_{11}s_{1}+\dotsb +n_{r1}s_{r}$ and 
take linear forms $s_{1}'(=\mu(s)),\dotsc,s_{r}'$ by
$(s_{1}',\dotsc,s_{r}')=(s_{1},\dotsc,s_{r})(n_{ij})_{i,j}$.
Then we have $\lambda_{j}s_{i}'=s_{i}'+\delta_{ij}$, where $\delta_{ij}$
denotes Kronecker's delta. There exists an irreducible polynomial
$h(t_{1},\dotsc,t_{r})\in k[t_{1},\dotsc,t_{r}]$ such that
$f(s)=h(s_{1}',\dotsc,s_{r}')$. But
\[ h(s_{1}',s_{2}'+m_{2},\dotsc,s_{r}'+m_{r})=\lambda_{2}^{m_{2}}\cdots
\lambda_{r}^{m_{r}}f(s)=f(s)=h(s_{1}',\dotsc,s_{r}') \]
for all $m_{2},\dotsc,m_{r}\in\bbz$, since $\lambda_{2},\dotsc,\lambda_{r}\in
\Xi_{f}$. Thus it follows that $h$ is actually a polynomial of one variable
$t_{1}$. The second part obviously follows by the definition of $\mu$.

(ii) By part (ii) of Lemma \ref{lem-Ore-Sato1}, we have
$Z^{1}(\Xi;\bbz(\Xi/\Xi_{f}))\simeq Z^{1}(\Xi/\Xi_{f};\bbz(\Xi/\Xi_{f}))$.
For an arbitrary $\alpha\in Z^{1}(\Xi/\Xi_{f};\bbz(\Xi/\Xi_{f}))$,
set $\eta=\alpha(\lambda)$. Then the assertion follows by induction on
$\mu(m)$, since the image of $\tau_{1}^{m_{1}}\cdots\tau_{r}^{m_{r}}$
in $\Xi/\Xi_{f}$ is $\lambda^{\mu(m)}$.

(iii) We see that $Z^{1}(\Xi/\Xi_{f};\bbz(\Xi/\Xi_{f}))\rightarrow
\bbz(\Xi/\Xi_{f})$, $\alpha\mapsto\alpha(\lambda)$ is an isomorphism.
Then $B^{1}(\Xi/\Xi_{f};\bbz(\Xi/\Xi_{f}))$ is isomorphic
to $\bbz(\Xi/\Xi_{f})(\lambda -1)$ under this isomorphism.
Hence $H^{1}(\Xi/\Xi_{f};\bbz(\Xi/\Xi_{f}))\simeq\bbz$, via the map
$\varepsilon : \bbz(\Xi/\Xi_{f})\rightarrow\bbz$,
$\sum_{i}z_{i}\lambda^{i}\mapsto\sum_{i}z_{i}$.
\end{proof}

By translating the lemmas above, we have the following theorem.

\begin{theorem}[Ore-Sato]
\label{thm-Ore-Sato}
Let $\Xi\rightarrow k(s)^{\times}$, $g\mapsto b_{g}(s)$ be a $1$-cocycle in
$Z^{1}(\Xi;k(s)^{\times})$. Write $\tau^{m}=\tau_{1}^{m_{1}}\cdots
\tau_{r}^{m_{r}}$ for $m=(m_{1},\dotsc,m_{r})\in\bbz^{r}$.
Then $b$ is written as the following form:
\begin{equation}
\label{b-expression}
b_{\tau^{m}}(s)=c(\tau^{m})\prod_{j=1}^{N}\zeta_{j}\frac{f_{j}(s+m)}{f_{j}(s)}
\prod_{i=1}^{l}\eta_{i}\cdot
\left\{\begin{array}{cl}
{\displaystyle\prod_{\nu=0}^{\mu_{i}(m)-1}}h_{i}(\mu_{i}(s)+\nu) &
(\mu_{i}(m)\geq 1) \vspace{1ex} \\
1 & (\mu_{i}(m)=0) \vspace{1ex} \\
{\displaystyle\prod_{\nu=\mu_{i}(m)}^{-1}\frac{1}{h_{i}(\mu_{i}(s)+\nu)}} &
(\mu_{i}(m)\leq -1)
\end{array}\right.
\end{equation}
for all $m\in\bbz^{r}$. Here, $c : \Xi\rightarrow k^{\times}$ is a character,
$h_{1},\dotsc, h_{l}\in k[t]$ are irreducible polynomials
of one variable, $\mu_{1},\dotsc,\mu_{l}\in \bbz s_{1}+\dotsb +\bbz s_{r}$ are
non-zero linear forms whose each greatest common divisor of non-zero
coefficients is $1$, $\eta_{i}\in\bbz(\Xi/\Xi_{\mu_{i}(s)})$ ($i=1,\dotsc,l$),
$f_{1},\dotsc,f_{N}$ are irreducible polynomials in $k[s]$ which satisfy that
$\rank \Xi/\Xi_{f_{j}} >1$ ($j=1,\dotsc,N$) and $f_{1}(s),\dotsc,f_{N}(s),
h_{1}(\mu_{1}(s)),\dotsc,h_{l}(\mu_{l}(s))$ are not translated into one another
up to constant multiple by the action of $\bbz\Xi$,
and $\zeta_{j}\in\bbz(\Xi/\Xi_{f_{j}})$ ($j=1,\dotsc,N$).
\end{theorem}

Let $\varepsilon : \bbz(\Xi/\Xi_{\mu_{i}(s)})\rightarrow\bbz$,
$\sum_{i}z_{i}\lambda^{i}\mapsto\sum_{i}z_{i}$ be the map described in
the last sentence in the proof of Lemma \ref{lem-Ore-Sato2}.
We normalize each $\eta_{i},\mu_{i},h_{i}$ in (\ref{b-expression})
so that $\varepsilon(\eta_{i})>0$ ($i=1,\dotsc,l$)
(when $\varepsilon(\eta_{i})<0$, replace $\eta_{i},\mu_{i},h_{i}(t)$ with
$-\eta_{i},-\mu_{i},h_{i}(-t-1)$). Put $f(s)=\prod_{j=1}^{N}\zeta_{j}f_{j}(s)$.

\begin{corollary}
\label{cor-Ore-Sato}
In Theorem \ref{thm-Ore-Sato}, assume $b_{\tau_{1}}(s),\dotsc,b_{\tau_{r}}(s)$
are polynomials. Take the expression (\ref{b-expression}) with
$f(s)=1$ (which necessarily holds) and $\eta_{i},\mu_{i},h_{i}$ normalized
as above. Then they satisfy the following conditions {\rm (1), (2)}.

{\rm (1)} All coefficients of each $\mu_{i}$ are non-negative integers.

{\rm (2)} For each $i=1,\dotsc,l$, let $\lambda_{i}$ be the basis of
$\Xi/\Xi_{\mu_{i}(s)}$ such that $\lambda_{i}\mu_{i}(s)=\mu_{i}(s)+1$.
Write $\eta_{i}=\sum_{j=n}^{n'}z_{ij}\lambda_{i}^{j}$ ($z_{ij}\in\bbz$,
$z_{in},z_{in'}\neq 0$). Then $h_{i}(\mu_{i}(s)+n)^{z_{in}}$ and
$h_{i}(\mu_{i}(s)+n'+\mu_{i}(m)-1)^{z_{in'}}$ are not canceled in the
product $\prod_{j=n}^{n'}\prod_{\nu=0}^{\mu_{i}(m)-1}
h_{i}(\mu_{i}(s)+j+\nu)^{z_{ij}}$ ($0\neq m\in\bbz_{\geq 0}^{r}$).
Hence $z_{in}$ and $z_{in'}$ are positive integers.
\end{corollary}
\begin{proof}
(1) Necessarily $\mu_{i}(e_{1})\varepsilon(\eta_{i}),\dotsc,
\mu_{i}(e_{r})\varepsilon(\eta_{i})$ are non-negative integers for
$i=1,\dotsc,l$ by the assumption.

(2) This is easily seen.
\end{proof}

\begin{remark}
\label{rem-Ore-Sato}
(i) The $1$-cocycle given by (\ref{b-expression}) and the following one both
define the same cohomology class in $H^{1}(\Xi;k(s)^{\times})$:
\[ \tau^{m}\mapsto c(\tau^{m})\prod_{i=1}^{l}
\left\{\begin{array}{cl}
{\displaystyle\prod_{\nu=0}^{\mu_{i}(m)-1}}
h_{i}(\mu_{i}(s)+\nu)^{\varepsilon(\eta_{i})} &
(\mu_{i}(m)\geq 1) \vspace{1ex} \\
1 & (\mu_{i}(m)=0) \vspace{1ex} \\
{\displaystyle\prod_{\nu=\mu_{i}(m)}^{-1}}
h_{i}(\mu_{i}(s)+\nu)^{-\varepsilon(\eta_{i})} & (\mu_{i}(m)\leq -1).
\end{array}\right.
\]
In \cite[\S 1]{Gelfand-Graev-Retakh}, $\eta_{i},\mu_{i},h_{i}$
are normalized to be $\varepsilon(\eta_{i})<0$. Another normalization is
given in \cite[Proposition 1.1.4]{Loeser-Sabbah} with a restriction on
$\mu_{i}$ instead of $\varepsilon(\eta_{i})$. To describe the group structure
of cohomology classes, such descriptions are more elegant.
But, in this article, the Ore-Sato theorem should be introduced
in the presented form for a reason which will arise later (especially in
Lemma \ref{lem-zero-or-one}). Here we are interested in
$Z^{1}(\Xi;k(s)^{\times})$ rather than $H^{1}(\Xi;k(s)^{\times})$.

(ii) The assumption in Corollary \ref{cor-Ore-Sato} does not imply
that each $\eta_{i}h_{i}(\mu_{i}(s))$ is a polynomial; for example, let $r=2$
and consider the $1$-cocycle defined by (\ref{b-expression}) with $c=1$,
$f(s)=1$, $l=1$, $h_{1}(t)=t$, $\mu_{1}(s)=2s_{1}+3s_{2}$, $\eta_{1}
=1-\lambda +\lambda^{2}$ (here $\lambda$ is the basis of
$\Xi/\Xi_{\mu_{1}(s)}$ represented by $\tau_{1}^{-1}\tau_{2}$).
In addition, we observe that the conditions (1), (2)
in the corollary are not sufficient for the assumption; consider
the $1$-cocycle where $\eta_{1}$ in the above example is replaced with
$\eta_{1}=1-\lambda -\lambda^{2}+\lambda^{3}+\lambda^{4}$.
(But these two examples define the same cohomology class in
$H^{1}(\Xi;k(s)^{\times})$).
\end{remark}

Let $k=\bbc$. Take a $1$-cocycle $b\in Z^{1}(\Xi;\bbc(s)^{\times})$ and keep
the notation in Theorem \ref{thm-Ore-Sato}. Write $h_{i}=t+\alpha_{i}$
($\alpha_{i}\in\bbc$) and $c(\tau^{m})=c_{1}^{m_{1}}\dotsm c_{r}^{m_{r}}$
($c_{1},\dotsc,c_{r}\in\bbc^{\times}$).
Then the $\Gm$-primitive difference equation associated with
$b$ has a solution
\[ \gamma(s)=c_{1}^{s_{1}}\dotsm c_{r}^{s_{r}}f(s)
\prod_{i=1}^{l}\eta_{i}\Gamma(\mu_{i}(s)+\alpha_{i}). \]
Assume $b_{\tau_{1}},\dotsc,b_{\tau_{r}}$ are polynomials. In this case,
$\gamma(s)$ has no zeros by Corollary \ref{cor-Ore-Sato}.
There are several methods to characterize $\gamma(s)$ among other solutions;
by asymptotic behavior or by log convexity:
M.\ Fujigami \cite{Fujigami} gave a generalization
of the Bohr-Mollerup theorem \cite[Theorem 2.1]{Artin}.
Further methods may be suggested in \cite[\S 6]{Artin}.

\section{Reductive prehomogeneous vector spaces and its $b$-functions}
\label{sec-b-functions}

Let $V$ be an $n$-dimensional $\bbc$-vector space and $G$ a connected reductive
linear algebraic group over $\bbc$. Suppose $(G,\rho,V)$ is a reductive
prehomogeneous vector space. Then, by the definition, there exists a proper
algebraic subset $S$ of $V$ such that $V\setminus S$ is a single $G$-orbit.
Let $S_{0}$ denote the union of the irreducible components of $S$ with
codimension $1$. We always assume that $S_{0}$ is not empty. Since $\rho(G)$
is connected reductive, it is self-adjoint with respect to a $\bbc$-basis of
$V$ by the theorem of Mostow \cite{Mostow}. By such a basis, we identify
the coordinate ring $\bbc[V]$ of $V$ with $\bbc[x]=\bbc[x_{1},\dotsc,x_{n}]$,
the polynomial ring of $n$ variables. Let $P_{1},\dotsc,P_{r}\in\bbc[x]$
be irreducible polynomials which define the irreducible components of $S_{0}$.
These polynomials are relative invariants and every relative invariant is
uniquely expressed in the form $cP_{1}(x)^{m_{1}}\dotsm P_{r}(x)^{m_{r}}$
($c\in\bbc^{\times}$, $(m_{1},\dotsc,m_{r})\in\bbz^{r}$).
In this sense, $P_{1},\dotsc,P_{r}$ are called
the {\em basic relative invariants} of $(G,\rho,V)$
\cite[Definition 2.10]{Kimura}. They are necessarily homogeneous polynomials
\cite[Corollary 2.7]{Kimura}. Let $\bar{P}_{1},\dotsc,\bar{P}_{r}$ be
the polynomials obtained by complex conjugation of coefficients of $P_{1},
\dotsc,P_{r}$ respectively. By the choice of the basis of $V$,
those are the basic relative invariants of the dual
prehomogeneous vector space $(G,\rho^{\ast},V^{\ast})$, when we identify
$\bbc[V^{\ast}]$ with $\bbc[x]$ by the dual basis (see \cite[Lemma 1.5]{Gyoja}
or \cite[Proposition 2.21]{Kimura}).

In the following, we consider $P_{1}(x)^{s_{1}}\dotsm P_{r}(x)^{s_{r}}$ as
a many-valued holomorphic function on $\bbc^{r}\times(V\setminus S_{0})$.
Write $\grad_{x}=(\partial/\partial x_{1},\dotsc,\partial/\partial x_{r})$.
For each $m=(m_{1},\dotsc,m_{r})\in\bbz_{\geq 0}^{r}$,
\[ \frac{\bar{P}_{1}(\grad_{x})^{m_{1}}\dotsm\bar{P}_{r}(\grad_{x})^{m_{r}}
(P_{1}(x)^{s_{1}+m_{1}}\dotsm P_{r}(x)^{s_{r}+m_{r}})}{P_{1}(x)^{s_{1}}
\dotsm P_{r}(x)^{s_{r}}}\in\bbc[s_{1},\dotsc,s_{r};x_{1},\dotsc,x_{n}] \]
is (absolutely) invariant under the action of $G$. Since every absolute
invariant in $\bbc[x]$ is a constant \cite[Proposition 2.4]{Kimura},
it is independent of $x$. Thus there exists a polynomial $b_{m}(s)\in
\bbc[s_{1},\dotsc,s_{r}]$ such that
\[ \bar{P}_{1}(\grad_{x})^{m_{1}}\dotsm\bar{P}_{r}(\grad_{x})^{m_{r}}
(P_{1}(x)^{s_{1}+m_{1}}\dotsm P_{r}(x)^{s_{r}+m_{r}})
=b_{m}(s)P_{1}(x)^{s_{1}}\dotsm P_{r}(x)^{s_{r}}. \]
Moreover, it is known that the degree of $b_{m}(s)$ (on $s$) is equal to the
degree of $\prod_{i=1}^{r}P_{i}(x)^{m_{i}}$ (on $x$); this follows from
an easy modification of the proof given in
\cite[Lemma 1.7]{Gyoja} or \cite[Proposition 2.22]{Kimura}.

\begin{definition}
\label{def-b-functions}
The polynomials $b_{m}(s)$ are called the {\em $b$-functions} of
$P_{1},\dotsc,P_{r}$.
\end{definition}

By calculating
\[ \bar{P}_{1}(\grad_{x})^{m_{1}+m_{1}'}\dotsm
\bar{P}_{r}(\grad_{x})^{m_{r}+m_{r}'}(P_{1}(x)^{s_{1}+m_{1}+m_{1}'}
\dotsm P_{r}(x)^{s_{r}+m_{r}+m_{r}'}) \]
in two ways, we have
\[ b_{m+m'}(s)=b_{m'}(s+m)b_{m}(s) \]
for all $m=(m_{1},\dotsc,m_{r}),m'=(m_{1}',\dotsc,m_{r}')\in
\bbz_{\geq 0}^{r}$. Hence the map $\bbz_{\geq 0}^{r}\rightarrow
\bbc(s)^{\times}$, $m\mapsto b_{m}(s)$ is uniquely extended to a
$1$-cocycle $b : \bbz^{r}\rightarrow\bbc(s)^{\times}$. By the Ore-Sato
theorem (Theorem \ref{thm-Ore-Sato}), each $b_{m}(s)$ ($0\neq m
=(m_{1},\dotsc,m_{r})\in\bbz_{\geq 0}^{r}$) is written as
\begin{equation}
\label{b-expression2}
b_{m}(s)=c_{1}^{m_{1}}\dotsm c_{r}^{m_{r}}\prod_{i=1}^{l}\prod_{j}
\prod_{\nu=0}^{\mu_{i}(m)-1}(\mu_{i}(s)+\alpha_{i}+j+\nu)^{z_{ij}},
\end{equation}
where the notations are taken as in \S \ref{sec-Ore-Sato} and $\eta_{i}
=\sum_{j}z_{ij}\lambda_{i}^{j}$ satisfying
the conditions in Corollary \ref{cor-Ore-Sato}. We take them
so that $(\mu_{1}(s)+\alpha_{1}),\dotsc,(\mu_{l}(s)+\alpha_{l})$ are
not translated into one another by the action of $\Xi$.
It is known that each $\alpha_{i}+j$ (with $z_{ij}>0$) is a positive rational
number (due to M.\ Kashiwara \cite{Kashiwara}).
Moreover, we have $c_{1},\dotsc,c_{r}\in\bbr_{> 0}$ since
\[ b_{m}(0)=\bar{P}_{1}(\grad_{x})^{m_{1}}\dotsm\bar{P}_{r}(\grad_{x})^{m_{r}}
(P_{1}(x)^{m_{1}}\dotsb P_{r}(x)^{m_{r}})
\in\bbr_{> 0} \]
for all $m=(m_{1},\dotsc,m_{r})\in\bbz_{\geq 0}^{r}$.
Let $d_{k}$ be the degree of $P_{k}(x)$ for $k=1,\dotsc,r$.
Then the observation just before Definition \ref{def-b-functions} implies that
\begin{equation}
\label{eq-degree-b-functions}
d_{k}=\sum_{i=1}^{l}\mu_{i}(e_{k})\varepsilon(\eta_{i})\qquad
(k=1,\dotsc,r),
\end{equation}
where $e_{1}=(1,0,\dotsc,0),\dotsc,e_{r}=(0,\dotsc,0,1)$ and
$\varepsilon(\eta_{i})=\sum_{j}z_{ij}$.

\section{Local zeta functions over $\bbc$}
\label{sec-complex-LZ}

We identify $V$ with $\bbc^{n}$ by the basis fixed in \S \ref{sec-b-functions}.
Let $dx$ denote the Haar measure on $V$ normalized to satisfy
\[ \int_{V}\exp(-2\pi x{}^{t}\bar{x})dx = 1, \]
where $x=(x_{1},\dotsc,x_{r})\in\bbc^{n}$ and ${}^{t}\bar{x}$ denotes
the transposition of the complex conjugate of $x$.
Let $|\cdot|_{\bbc}$ be the valuation of $\bbc$ defined by
$|z|_{\bbc}=z\bar{z}=|z|^{2}$ for $z\in\bbc$. The integral
\[ Z_{\bbc}(s)=\int_{V}|P_{1}(x)|_{\bbc}^{s_{1}}\dotsm|P_{r}(x)|_{\bbc}^{s_{r}}
\exp(-2\pi x{}^{t}\bar{x})dx \]
converges when $s\in\{(s_{1},\dotsc,s_{r})\in\bbc^{r}
\;|\;\re(s_{1}),\dotsc,\re(s_{r})>0\}$, and hence $Z_{\bbc}(s)$ is
a holomorphic function on this region.
Our purpose is to show that $Z_{\bbc}(s)$ is equal to
\[ \gamma_{\bbc}(s):=\prod_{k=1}^{r}((2\pi)^{-d_{k}}c_{k})^{s_{k}}
\prod_{i=1}^{l}\prod_{j}\left(
\frac{\Gamma(\mu_{i}(s)+\alpha_{i}+j)}{\Gamma(\alpha_{i}+j)}\right)^{z_{ij}} \]
with the notations in (\ref{b-expression2}).

\begin{theorem}
\label{thm-complex-LZ}
$Z_{\bbc}(s)$ has a meromorphic continuation to $\bbc^{r}$ and
$Z_{\bbc}(s)=\gamma_{\bbc}(s)$.
\end{theorem}
Actually, the first part of this theorem has been known (see \cite{Bernshtein})
and our aim is to obtain the second part. To prove this theorem,
we need the following lemma:

\begin{lemma}
\label{lem-complex-LZ}
For $s=(s_{1},\dotsc,s_{r})\in\bbc^{r}$ and
$m=(m_{1},\dotsc,m_{r})\in\bbz_{\geq 0}^{r}$, we have
\[ \begin{array}{ll}
 & \bar{P}_{1}(\grad_{x})^{m_{1}}\dotsm\bar{P}_{r}(\grad_{x})^{m_{r}}
(|P_{1}(x)|_{\bbc}^{s_{1}}\dotsm|P_{r}(x)|_{\bbc}^{s_{r}}
P_{1}(x)^{m_{1}}\dotsm P_{r}(x)^{m_{r}}) \vspace{1ex} \\
= & b_{m}(s)|P_{1}(x)|_{\bbc}^{s_{1}}\dotsm |P_{1}(x)|_{\bbc}^{s_{1}}
\end{array} \]
on $V\setminus S_{0}$.
\end{lemma}
\begin{proof}
Locally we choose the branching of the value of $\log P_{i}(x)$ and
$\log\overline{P_{i}(x)}=\log\bar{P}_{i}(\bar{x})$ so that
$|P_{i}(x)|_{\bbc}^{s_{i}}=\bar{P}_{i}(\bar{x})^{s_{i}}P_{i}(x)^{s_{i}}$
holds. Since $\bar{P}_{j}(\grad_{x})$ and
$\bar{P}_{i}(\bar{x})^{s_{i}}$ commute as differential operators
for $i,j=1,\cdots,r$, we have
\begin{eqnarray*}
 & & \bar{P}_{1}(\grad_{x})^{m_{1}}\dotsm\bar{P}_{r}(\grad_{x})^{m_{r}}
(|P_{1}(x)|_{\bbc}^{s_{1}}\dotsm|P_{r}(x)|_{\bbc}^{s_{r}}
P_{1}(x)^{m_{1}}\dotsm P_{r}(x)^{m_{r}}) \\
 & = & \bar{P}_{1}(\grad_{x})^{m_{1}}\dotsm\bar{P}_{r}(\grad_{x})^{m_{r}}
(\bar{P}_{1}(\bar{x})^{s_{1}}\dotsm\bar{P}_{r}(\bar{x})^{s_{r}}
P_{1}(x)^{s_{1}+m_{1}}\dotsm P_{r}(x)^{s_{r}+m_{r}}) \\
 & = & \bar{P}_{1}(\bar{x})^{s_{1}}\dotsm\bar{P}_{r}(\bar{x})^{s_{r}}
\bar{P}_{1}(\grad_{x})^{m_{1}}\dotsm\bar{P}_{r}(\grad_{x})^{m_{r}}
(P_{1}(x)^{s_{1}+m_{1}}\dotsm P_{r}(x)^{s_{r}+m_{r}}) \\
 & = & \bar{P}_{1}(\bar{x})^{s_{1}}\dotsm\bar{P}_{r}(\bar{x})^{s_{r}}
b_{m}(s)P_{1}(x)^{s_{1}}\dotsm P_{r}(x)^{s_{r}} \\
 & = & b_{m}(s)|P_{1}(x)|_{\bbc}^{s_{1}}\dotsm |P_{1}(x)|_{\bbc}^{s_{1}}.
\end{eqnarray*}
\end{proof}

\begin{proof}[Proof of Theorem \ref{thm-complex-LZ}]
By Lemma \ref{lem-complex-LZ}, we have
\begin{eqnarray*}
 & & b_{m}(s)Z_{\bbc}(s) \\
 & = & \int_{V}\left\{\left(\prod_{i=1}^{r}
\bar{P}_{i}(\grad_{x})^{m_{i}}\right)\left(\prod_{j=1}^{r}
|P_{j}(x)|_{\bbc}^{s_{j}}P_{j}(x)^{m_{j}}\right)\right\}
\exp(-2\pi x{}^{t}\bar{x})dx \\
 & = & \int_{V}\left(\prod_{j=1}^{r}|P_{j}(x)|_{\bbc}^{s_{j}}P_{j}(x)^{m_{j}}
\right)\left(\prod_{i=1}^{r}\bar{P}_{i}(-\grad_{x})^{m_{i}}\right)
\exp(-2\pi x{}^{t}\bar{x})dx \\
 & = & (2\pi)^{\sum_{k=1}^{r}d_{k}m_{k}}
\int_{V}\prod_{j=1}^{r}|P_{j}(x)|_{\bbc}^{s_{j}}P_{j}(x)^{m_{j}}
\prod_{i=1}^{r}\bar{P}_{i}(\bar{x})^{m_{i}}\exp(-2\pi x{}^{t}\bar{x})dx \\
 & = & (2\pi)^{\sum_{k=1}^{r}d_{k}m_{k}}
\int_{V}\prod_{j=1}^{r}|P_{j}(x)|_{\bbc}^{s_{j}+m_{j}}
\exp(-2\pi x{}^{t}\bar{x})dx \\
 & = & (2\pi)^{\sum_{k=1}^{r}d_{k}m_{k}}Z_{\bbc}(s+m)
\end{eqnarray*}
for $(s_{1},\dotsc,s_{r})\in\{(s_{1},\dotsc,s_{r})\in\bbc^{r}\;|\;\re(s_{1}),
\dotsc,\re(s_{r})>0\}$ and $m=(m_{1},\dotsc,m_{r})\in\bbz_{\geq 0}^{r}$.
Hence $Z_{\bbc}(s)$ satisfies the following equation:
\begin{equation}
\label{Z-difference-equation1}
Z_{\bbc}(s+m)=(2\pi)^{-\sum_{k=1}^{r}d_{k}m_{k}}b_{m}(s)Z_{\bbc}(s).
\end{equation}
Then $Z_{\bbc}(s)$ has a meromorphic continuation to $\bbc^{r}$
by this equation (as in \cite{Bernshtein}).
Furthermore, (\ref{Z-difference-equation1}) implies that both $Z_{\bbc}(s)$
and $\gamma_{\bbc}(s)$ satisfy the same ($\Gm$-primitive) difference equation.
Therefore $C(s):=Z_{\bbc}(s)/\gamma_{\bbc}(s)$ is a holomorphic and periodic
function with periods $e_{1},\dotsc,e_{r}$.
To show that $C(s)$ is a constant function, we investigate the asymptotic
behavior of $C(s)$ on the strip $\cS=\{(s_{1},\dotsc,s_{r})\in\mathbb{C}^{r}
\;|\; 1\leq \re(s_{i})\leq 2\quad (i=1,\cdots,r)\}$. 

Put $S^{n-1}(\bbc)=\{x\in V\;|\;x{}^{t}\bar{x}=1\}$
($\simeq SO(2n,\bbr)/SO(2n-1,\bbr)$ as real manifolds).
We identify $V\setminus\{0\}$ with $\bbr_{>0}\times
S^{n-1}(\bbc)$ via $x\mapsto (\xi,u)=(\sqrt{x{}^{t}\bar{x}},
x/\sqrt{x{}^{t}\bar{x}})$. Take the $SO(2n,\bbr)$-invariant measure $du$ on
$S^{n-1}(\bbc)$ such that $dx = 2^{n}\xi^{2n-1}d\xi du$ on $V\setminus\{0\}$.
Let
\[ \psi(s)=2^{n-1}\int_{S^{n-1}(\bbc)}
|P_{1}(u)|_{\bbc}^{s_{1}}\dotsm|P_{r}(u)|_{\bbc}^{s_{r}}du. \]
Then we have
\begin{eqnarray*}
Z_{\bbc}(s) & = &
\int_{V\setminus\{0\}}|P_{1}(x)|_{\bbc}^{s_{1}}\dotsm|P_{r}(x)|_{\bbc}^{s_{r}}
\exp(-2\pi x{}^{t}\bar{x})dx \\
 & = & 2\psi(s)\int_{0}^{\infty}
\xi^{2(\sum_{k=1}^{r}d_{k}s_{k}+n)-1}\exp(-2\pi \xi^{2})d\xi \\
 & = & (2\pi)^{-\sum_{k=1}^{r}d_{k}s_{k}-n}\psi(s)
\Gamma(d_{1}s_{1}+\dotsb +d_{r}s_{r}+n)
\end{eqnarray*}
when $\re(s_{i})>0$ ($i=1,\cdots,r$).
Since
\[ |(2\pi)^{-\sum_{k=1}^{r}d_{k}s_{k}-n}\psi(s)|\leq
(2\pi)^{-\sum_{k=1}^{r}d_{k}\re(s_{k})-n}2^{n-1}\int_{S^{n-1}(\bbc)}
|P_{1}(u)|_{\bbc}^{\re(s_{1})}\dotsm|P_{r}(u)|_{\bbc}^{\re(s_{r})}du, \]
$|(2\pi)^{-\sum_{k=1}^{r}d_{k}s_{k}-n}\psi(s)|$ is bounded in $\cS$.
Let $a_{1},\dotsc,a_{r}$ be arbitrary positive real numbers and
$t_{0}$ a real variable.
Then the well-known asymptotic behavior of the gamma function and
(\ref{eq-degree-b-functions}) imply that
\[ C(1+\sqrt{-1}a_{1}t_{0},\dotsc,1+\sqrt{-1}a_{r}t_{0})=o(\exp(|t_{0}|))
\qquad (|t_{0}|\rightarrow\infty). \]
Unless $C(s)$ is a constant function, this is impossible.
Thus the proof is completed since $C(0)=1$.
\end{proof}

\section{Local zeta functions over $\bbr$}
\label{sec-real-LZ}

In this section, we assume that the prehomogeneous vector space $(G,\rho,V)$
is {\em defined over $\bbr$} (in the sense of \cite[\S 1]{FSato1982} or
\cite[\S 2.1]{Kimura}) and replace $P_{1},\dotsc,P_{r}$ and $b$ as follows.
Since $S_{0}$ is defined over $\bbr$ (see \cite[Lemma 1.1]{FSato1982}),
we can take irreducible polynomials $P_{1},\dotsc,P_{r}\in\bbr[x]$ which
define the $\bbr$-irreducible components of $S_{0}$ (possibly $r$ becomes
smaller). Here, we are assuming the basis of $V$ is fixed so that
the $\bbr$-rational points of $\rho(G)$ is self-adjoint with respect to
the induced $\bbr$-basis of $V_{\bbr}$, the $\bbr$-rational points of $V$.
Those $P_{1},\dotsc,P_{r}$ are often called {\em the basic
relative invariants of $(G,\rho,V)$ over $\bbr$}.
They are also considered as the basic relative invariants of $(G,\rho^{\ast},
V^{\ast})$ over $\bbr$. Then there exist polynomials $b_{m}(s)$ such that
\[ P_{1}(\grad_{x})^{m_{1}}\dotsm P_{r}(\grad_{x})^{m_{r}}
(P_{1}(x)^{s_{1}+m_{1}}\dotsm P_{r}(x)^{s_{r}+m_{r}})
=b_{m}(s)P_{1}(x)^{s_{1}}\dotsm P_{r}(x)^{s_{r}} \]
for $m=(m_{1},\dotsc,m_{r})\in\bbz_{\geq 0}^{r}$.
All properties on $b_{m}(s)$ described in \S \ref{sec-b-functions} 
also hold. Keep the notations in (\ref{b-expression2}) and in
(\ref{eq-degree-b-functions}).

Let $dx$ denote the Lebesgue measure on $V_{\bbr}$ identified with
$\bbr^{r}$. As in the previous section, the integral
\[ Z_{\bbr}(s)=\int_{V_{\bbr}}|P_{1}(x)|^{s_{1}}\dotsm|P_{r}(x)|^{s_{r}}
\exp(-\pi x{}^{t}x)dx \]
converges when $(s_{1},\dotsc,s_{r})\in\{(s_{1},\dotsc,s_{r})\in\bbc^{r}
\;|\;\re(s_{1}),\dotsc,\re(s_{r})>0\}$, and hence $Z_{\bbr}(s)$ is
a holomorphic function on this region. From now on, we assume the following:
\begin{assumption}
\label{assum-multiplicity-free}
Every term of $P_{i}(x)$ is a multilinear form on $x$ for $i=1,\dotsc,r$,
i.e.\ each $P_{i}(x)$ is of the form:
\[ P_{i}(x)=\sum_{1\leq j_{1}<\dotsb <j_{d_{i}}\leq n}
a_{j_{1}\cdots j_{d_{i}}}x_{j_{1}}\dotsm x_{j_{d_{i}}}
\qquad (i=1,\dotsc,r). \]
\end{assumption}
We will see that $Z_{\bbr}(2s)$ satisfies a certain $\Gm$-primitive difference
equation in such a case.
By the assumption above, we see
\[ P_{i}(\grad_{x})\exp(-\pi x{}^{t}x)=(-2\pi)^{d_{i}}P_{i}(x)
\exp(-\pi x{}^{t}x)\qquad (i=1,\cdots,r). \]
Thus
\begin{eqnarray*}
b_{e_{i}}(s)Z_{\bbr}(s) & = & \int_{V_{\bbr}}\{P_{i}(\grad_{x})
(|P_{1}(x)|^{s_{1}}\dotsb|P_{r}(x)|^{s_{r}}P_{i}(x))\}\exp(-\pi x{}^{t}x)dx \\
 & = & (-1)^{d_{i}}\int_{V_{\bbr}}(|P_{1}(x)|^{s_{1}}\dotsb|P_{r}(x)|^{s_{r}}
P_{i}(x))P_{i}(\grad_{x})\exp(-\pi x{}^{t}x)dx \\
 & = & (2\pi)^{d_{i}}\int_{V_{\bbr}}|P_{1}(x)|^{s_{1}}\dotsb|P_{r}(x)|^{s_{r}}
P_{i}(x)P_{i}(x)\exp(-\pi x{}^{t}x)dx \\
 & = & (2\pi)^{d_{i}}\int_{V_{\bbr}}|P_{1}(x)|^{s_{1}}\dotsb
|P_{i}(x)|^{s_{i}+2}\dotsb|P_{r}(x)|^{s_{r}}\exp(-\pi x{}^{t}x)dx \\
 & = & (2\pi)^{d_{i}}Z_{\bbr}(s+2e_{i}),
\end{eqnarray*}
for $i=1,\dotsc,r$. Hence $Z_{\bbr}(s)$ satisfies the equation
\begin{equation}
\label{Z-difference-equation2}
Z_{\bbr}(s+2e_{i})=(2\pi)^{-d_{i}}b_{e_{i}}(s)Z_{\bbr}(s)\qquad
(i=1,\dotsc,r).
\end{equation}
We can obtain a meromorphic continuation of $Z_{\bbr}(s)$ to $\bbc^{r}$.
By considering $Z_{\mathbb{R}}(s+2e_{j}+2e_{k})$ in two ways, we have
\begin{equation}
\label{special-cocycle}
b_{e_{j}}(s)b_{e_{k}}(s+2e_{j})=b_{e_{j}}(s+2e_{k})b_{e_{k}}(s)
\quad (j,k=1,\cdots,r).
\end{equation}

\begin{lemma}
\label{lem-zero-or-one}
When Assumption \ref{assum-multiplicity-free} holds,
$\mu_{i}(e_{j})$ ($i=1,\cdots,l$, $j=1,\cdots,r$) are equal to
either $0$ or $1$
and hence each $\eta_{i}\cdot(\mu_{i}(s)+\alpha_{i})$ is a polynomial.
\end{lemma}
\begin{proof}
For each $i=1,\cdots,l$, take any $j,k$ such that $\mu_{i}(e_{j}),
\mu_{i}(e_{k})>0$ and write $\eta_{i}=\sum_{u=\kappa}^{\kappa'}
z_{iu}\lambda_{i}^{u}$ ($z_{i\kappa},z_{i\kappa'}\neq 0$).
By (\ref{special-cocycle}), we have
\begin{eqnarray*}
 & & \prod_{u=\kappa}^{\kappa'}\left(\prod_{v=0}^{\mu_{i}(e_{j})-1}
(\mu_{i}(s)+\alpha_{i}+u+v)
\prod_{w=0}^{\mu_{i}(e_{k})-1}
(\mu_{i}(s)+2\mu_{i}(e_{j})+\alpha_{i}+u+w)\right)^{z_{iu}} \\
 & = & \prod_{u=\kappa}^{\kappa'}\left(\prod_{v=0}^{\mu_{i}(e_{j})-1}
(\mu_{i}(s)+2\mu_{i}(e_{k})+\alpha_{i}+u+v)
\prod_{w=0}^{\mu_{i}(e_{k})-1}
(\mu_{i}(s)+\alpha_{i}+u+w)\right)^{z_{iu}}.
\end{eqnarray*}
In each side, the constant terms of the factors
\[ \begin{array}{l}
(\mu_{i}(s)+2\mu_{i}(e_{j})+\alpha_{i}+\kappa'+\mu_{i}(e_{k})-1), \vspace{1ex}
\\
(\mu_{i}(s)+2\mu_{i}(e_{k})+\alpha_{i}+\kappa'+\mu_{i}(e_{j})-1)
\end{array} \]
are maximal respectively; recall that these factors are not canceled 
(Corollary \ref{cor-Ore-Sato}). Hence the two factors coincide and
we have $\mu_{i}(e_{j})=\mu_{i}(e_{k})$. Therefore the all non-zero numbers
among $\mu_{i}(e_{1}),\dotsc,\mu_{i}(e_{r})$ coincide. This proves the
lemma since the greatest common divisor of them is $1$.
\end{proof}

\begin{theorem}
\label{thm-real-LZ}
When Assumption \ref{assum-multiplicity-free} holds, we have
\[ Z_{\mathbb{R}}(s)
=\prod_{k=1}^{r}(\pi^{-d_{i}}c_{i})^{\frac{s_{i}}{2}}\prod_{i=1}^{l}\prod_{j}
\left(\frac{\Gamma((\mu_{i}(s)+\alpha_{i}+j)/2)}{\Gamma((\alpha_{i}+j)/2)}
\right)^{z_{ij}}. \]
\end{theorem}
\begin{proof}
Let $\beta : \bbz^{r}\rightarrow\bbc(s)^{\times}$ be the map given by
\[ \beta_{m}(s)=\prod_{k=1}^{r}(\pi^{-d_{k}}c_{k})^{m_{k}}
\prod_{i=1}^{l}\prod_{j}\left\{
\begin{array}{cl}
{\displaystyle\prod_{\nu=0}^{\mu_{i}(m)-1}
(\mu_{i}(s)+\frac{\alpha_{i}+j}{2}+\nu)^{z_{ij}}} & (\mu_{i}(m)\geq 1)
\vspace{1ex} \\
1 & (\mu_{i}(m)=0) \vspace{1ex} \\
{\displaystyle\prod_{\nu=\mu_{i}(m)}^{-1}
(\mu_{i}(s)+\frac{\alpha_{i}+j}{2}+\nu)^{-z_{ij}}} &
(\mu_{i}(m)\leq -1).
\end{array}\right. \]
We see that $\beta$ is a $1$-cocycle since $Z_{\bbr}(2s)$ satisfies
the $\Gm$-primitive difference equation
\[ Z_{\bbr}(2(s+m))=\beta_{m}(s)Z_{\bbr}(2s) \]
by (\ref{Z-difference-equation2}) and Lemma \ref{lem-zero-or-one}.
In addition, the second assertion of Lemma \ref{lem-zero-or-one} implies
that all $z_{ij}$ is non-negative. Let
\[ \gamma_{\bbr}(s):=
\prod_{k=1}^{r}(\pi^{-d_{i}}c_{i})^{\frac{s_{i}}{2}}\prod_{i=1}^{l}\prod_{j}
\left(\frac{\Gamma((\mu_{i}(s)+\alpha_{i}+j)/2)}{\Gamma((\alpha_{i}+j)/2)}
\right)^{z_{ij}}. \]
Then both $\gamma_{\bbr}(2s)$ and $Z_{\bbr}(2s)$ satisfy the same difference
equation.

Put $S^{n-1}(\bbr)=\{x\in V_{\bbr}\;|\;x{}^{t}x=1\}$
($\simeq SO(n,\bbr)/SO(n-1,\bbr)$).
We identify $V_{\bbr}\setminus\{0\}$ with $\bbr_{>0}\times
S^{n-1}(\bbr)$ via $x\mapsto (\xi,u)=(\sqrt{x{}^{t}x},
x/\sqrt{x{}^{t}x})$. Take the $SO(n,\bbr)$-invariant measure $du$ on
$S^{n-1}(\bbr)$ such that $dx = \xi^{n-1}d\xi du$ on $V_{\bbr}\setminus\{0\}$.
Let
\[ \psi(s)=\frac{1}{2}\int_{S^{n-1}(\bbr)}
|P_{1}(u)|^{s_{1}}\dotsm|P_{r}(u)|^{s_{r}}du. \]
Then we have
\begin{eqnarray*}
Z_{\bbr}(s) & = &
\int_{V_{\bbr}\setminus\{0\}}|P_{1}(x)|^{s_{1}}\dotsm|P_{r}(x)|^{s_{r}}
\exp(-\pi x{}^{t}x)dx \\
 & = & 2\psi(s)\int_{0}^{\infty}
\xi^{\sum_{k=1}^{r}d_{k}s_{k}+n-1}\exp(-\pi \xi^{2})d\xi \\
 & = & \pi^{(-\sum_{k=1}^{r}d_{k}s_{k}-n)/2}\psi(s)
\Gamma\left(\frac{d_{1}s_{1}+\dotsb +d_{r}s_{r}+n}{2}\right)
\end{eqnarray*}
when $\re(s_{i})>0$ ($i=1,\cdots,r$). Therefore, similarly to the proof of
Theorem \ref{thm-complex-LZ}, we obtain that $Z_{\bbr}(2s)/\gamma_{\bbr}(2s)
=1$.
\end{proof}


\bibliographystyle{amsalpha}

\newpage

\setcounter{section}{0}

\part{Picard-Vessiot theories for artinian simple module algebras}
\label{part-PV-theory}

\section*{Introduction of Part \ref{part-PV-theory}}

The purpose of this part is to develop a unified Picard-Vessiot theory,
including Picard-Vessiot theories for differential equations and for
difference equations. The presented result was obtained by the author and
Masuoka. The contents from Section \ref{sec-basic-notions} to
Section \ref{sec-splitting} are made by retouching
\cite{Amano-Masuoka} including some omitted proofs. Section
\ref{sec-Liouville} was recently published as another paper \cite{Liouville},
in which some descriptions, such as the definition of $\sG$-extension, etc,
are improved than this article (but a sentence ``...the smoothness assumption;
namely, it is {\em spanned by} divided power sequences..." in Introduction
of \cite{Liouville} may be a little mistake).

In the usual sense, the ``Picard-Vessiot theory" means a Galois theory
for linear ordinary differential equations. See \cite{vanderPut-Singer2003}
for modern treatment. For example, consider the following differential
equation over $\bbc$:
\begin{equation}
\label{intro-example1}
y''(x)-y'(x)-y(x)=0.
\end{equation}
Let $\bbc[\partial]$ ($\partial=d/dx$) be the ring of differential operators
with constant coefficients. The differential module ($\bbc[\partial]$-module)
associated with the equation (\ref{intro-example1}) is
\[ \bbc[\partial]/\bbc[\partial](\partial^{2}-\partial-1)
\simeq
(\bbc[\partial]/\bbc[\partial](\partial-{\textstyle\frac{1+\sqrt{5}}{2}}))
\oplus
(\bbc[\partial]/\bbc[\partial](\partial-{\textstyle\frac{1-\sqrt{5}}{2}})).
\]
Thus the space of solutions is given by the $2$-dimensional
$\bbc$-vector space $\bbc\alpha +\bbc\beta$ with
$\alpha=e^{\frac{1+\sqrt{5}}{2}x}$, $\beta=e^{\frac{1-\sqrt{5}}{2}x}$.
A differential field (i.e.\ a field given a derivation)
$L$ including this space is called a {\em splitting field} for the equation. 
If $L$ is minimal with this property, it is called a {\em minimal splitting
field}. For the equation above, $L=\bbc(\alpha,\beta)$ is a minimal splitting
field. Like Galois extensions, $L/\bbc$ is then an extension of a special
type, called a {\em Picard-Vessiot extension}. For such an extension,
we can take a Galois group, called the {\em differential Galois group}
(or the {\em Picard-Vessiot group}) as an algebraic group defined by
$\Aut_{\bbc[\partial],\mbox{\scriptsize $\bbc$-alg}}(L)=:G(L/\bbc)$,
where $\Aut_{\bbc[\partial],\mbox{\scriptsize $\bbc$-alg}}$ denotes
the $\bbc[\partial]$-linear and $\bbc$-algebra automorphisms.
We have $G(L/\bbc)=\bG_{m}\times\bG_{m}$ in this case, and we can obtain the
Galois correspondence between closed subgroups of $G(L/\bbc)$ and
intermediate differential fields of $L/\bbc$. For example, the differential
field $\bbc(\alpha)$ corresponds to $\{1\}\times\bG_{m}$ (or $\bG_{m}\times
\{1\}$ according to the choice of the group action).

An analogous theory for {\em difference} equations is also known.
Bia{\l}ynicki-Birula \cite{Bialynicki-Birula} and Franke
\cite{Franke} first developed such a theory for inversive difference
fields, i.e.\ fields given an automorphism
(though Bia{\l}ynicki-Birula's paper was intended for more general theory,
not only for difference fields). A definition of Picard-Vessiot extensions
of inversive difference fields and Galois correspondences were obtained
there. But the theory had a difficulty on the existence of suitable solution
fields. For example, consider the Fibonacci recurrence
\begin{equation}
\label{intro-example2}
a(n+2)-a(n+1)-a(n)=0.
\end{equation}
Let $\bbc[\tau,\tau^{-1}]$ ($\tau : n\mapsto n+1$) be the ring of difference
operators with constant coefficients.
The difference module ($\bbc[\tau,\tau^{-1}]$-module)
associated with the equation (\ref{intro-example2}) is
\[ \begin{array}{ll}
\bbc[\tau,\tau^{-1}]/\bbc[\tau,\tau^{-1}](\tau^{2}-\tau-1) \\
\simeq (\bbc[\tau,\tau^{-1}]/\bbc[\tau,\tau^{-1}](\tau-\frac{1+\sqrt{5}}{2}))
\oplus (\bbc[\tau,\tau^{-1}]/\bbc[\tau,\tau^{-1}](\tau-\frac{1-\sqrt{5}}{2})).
\end{array} \]
Let $\cS_{\bbc}$ denote the ring of complex sequences
(see \cite[Example 3]{vanderPut-Singer1997}).
The space of solutions in $\cS_{\bbc}$ is given by the $2$-dimensional
$\bbc$-vector space $\bbc\alpha +\bbc\beta$ with
$\alpha=\{(\frac{1+\sqrt{5}}{2})^{n}\}$, $\beta=\{(\frac{1-\sqrt{5}}{2})^{n}\}
\in\cS_{\bbc}$. But one can not take any splitting field which becomes
a Picard-Vessiot extension for this equation.
If a subring in $\cS_{\bbc}$ contains $\alpha,\beta$,
then it has a zero divisor:
\[ (\alpha\beta-1)(\alpha\beta+1)=0. \]
On the other hand, if we take any inversive difference field which includes
a $2$-dimensional $\bbc$-vector space of solutions of (\ref{intro-example2}),
then it necessarily contains a new constant
(see \cite[p.\ 2]{vanderPut-Singer1997}).
However, overcoming this difficulty, the Picard-Vessiot theory for difference
equations in modern sense was developed by van der Put and Singer
\cite{vanderPut-Singer1997} with the notion of {\em Picard-Vessiot rings},
as follows. 
Consider the Laurent polynomial ring $\bbc[x,y,(xy)^{-1}]$ 
as an inversive difference ring by 
\[ \tau x=\frac{1+\sqrt{5}}{2}x,\quad \tau y =\frac{1-\sqrt{5}}{2}y. \]
On sees that $\lan (xy-1)(xy+1)\ran\subset \bbc[x,y,(xy)^{-1}]$ is a maximal
difference ideal. Put
\[ A=\bbc[x,y,(xy)^{-1}]/\lan (xy-1)(xy+1)\ran. \]
Then $A$ is a Picard-Vessiot ring for the equation (\ref{intro-example2})
in their sense (see \cite[Definition 1.4]{vanderPut-Singer1997}).
The total quotient ring of a Picard-Vessiot ring is called the {\em
total Picard-Vessiot ring} \cite[Definition 1.22]{vanderPut-Singer1997}.
For the equation (\ref{intro-example2}), we have the following total
Picard-Vessiot ring:
\begin{equation}
\label{intro-example-AS}
\begin{array}{rcl}
Q(A) & \simeq & \bbc(\alpha)\times\bbc(\alpha) \vspace{1ex} \\
x & \mapsto & (\alpha,\alpha) \vspace{1ex} \\
y & \mapsto & (\alpha^{-1},-\alpha^{-1}). \vspace{1ex}
\end{array}
\end{equation}
Then the {\em difference Galois group} for the equation can be defined by
$G(Q(A)/\bbc):=\Aut_{\bbc[\tau,\tau^{-1}],\mbox{\scriptsize $\bbc$-alg}}(Q(A))
$. In this case, we have $G(Q(A)/\bbc)\simeq
\Gm\times\bbz/2\bbz$. We obtain the Galois correspondence between 
closed subgroups of $G(Q(A)/\bbc)$ and intermediate difference subrings of
$Q(A)/\bbc$ such that every non-zero divisor is invertible
(see \cite[Theorem 1.29]{vanderPut-Singer1997}).
For example, $\bbc(\alpha)$ ($=\bbc(\alpha)(1,1)$) corresponds to
$\{1\}\times\bbz/2\bbz$ and $\bbc\times\bbc$ corresponds to $\Gm\times\{1\}$.

A unified approach to both differential and difference cases was first
attempted by Bia{\l}ynicki-Birula \cite{Bialynicki-Birula}, though it was
a theory for field extensions. Including the case that the solution algebras
can have zero divisors, Andr\'e \cite{Andre} gave such a unified approach
from the viewpoint of non-commutative differential geometry with the
theory of tannakian categories \cite{Deligne, Deligne-Milne}.
Alternatively we develop a unified Picard-Vessiot theory by a different way
based on Takeuchi's Hopf algebraic approach:
Takeuchi \cite{Takeuchi1989} beautifully clarified
the heart of the Picard-Vessiot theory in the generalized context of
$C$-ferential fields, intrinsically defining PV extensions and the minimal
splitting fields of $C$-ferential modules. By replacing linear algebraic
groups with affine group schemes (or equivalently commutative Hopf algebras),
he succeeded in removing from many of the results the assumptions of
finite generation, zero characteristic and algebraic closedness.
For a cocommutative coalgebra $C$ with a specific grouplike $1_{C}$,
a {\em $C$-ferential field} \cite{Takeuchi1989}
is a field given a unital, measuring action by $C$; the concept includes
differential fields, $\Delta$-fields \cite{Kolchin1973}, fields with
higher derivations \cite{Okugawa}, and difference fields (even non-inversive
ones are included). However, it was also a theory for field extensions and
the assumption that the tensor bialgebra $T(C^{+})$
\cite[p.\ 485]{Takeuchi1989} is a Birkhoff-Witt coalgebra (see
\cite[p.\ 504]{Takeuchi1989} or Assumption \ref{assum-BW-type}), is required
for the existence theorem of minimal splitting algebras.

In this article, we consider module algebras over a
cocommutative, pointed smooth Hopf algebra $D$. Thus $D$ is of the form
$D=D^{1}\# RG$ over a fixed field, say $R$, where $G$ is the group of
grouplikes in $D$, and the irreducible component $D^{1}$ containing $1$ is
a Birkhoff-Witt coalgebra. An inversive difference ring
which includes $R$ in its constants is precisely a $D$-module algebra
when $D^{1}=R$ and $G$ is the free group with one generator.
Differential rings are also within our scope, though only in characteristic
zero because of the smoothness assumption. Precisely a differential ring
which includes $R$ (of zero characteristic) in its constants is a $D$-module
algebra when $D^{1}=R[\partial]$ with a primitive $\partial$ and
$G$ is trivial. Algebras with higher derivations of infinite length fit in the
assumption, in arbitrary characteristic. 
An algebra (over $R$) with $R$-linear higher derivations $d_{0}=\mathrm{id},
d_{1},d_{2},\dotsc$ of infinite length is precisely a module algebra over the
Hopf algebra $R\lan d_{1},d_{2},\dotsc\ran$, which denotes the
(non-commutative) free algebra generated by $d_{1},d_{2},\dotsc$, and in
which $1,d_{1},d_{2},\dotsc$ form a divided power sequence. This Hopf algebra
becomes a Birkhoff-Witt coalgebra in arbitrary characteristic.

Throughout this article, $D$-module algebras are all supposed to be
commutative.
A $D$-module algebra $K$ is said to be {\em artinian simple} (AS) if
it is artinian as a ring and simple as a $D$-module algebra.
The last condition means that $K$ has no non-trivial $D$-stable ideal.
For example, the total Picard-Vessiot ring considered in
(\ref{intro-example-AS}) is an AS $\bbc[\tau,\tau^{-1}]$-module algebra.
Of course differential fields over $\bbc$ are AS $\bbc[\partial]$-module
algebras. In this sense we can generalize and unify the Picard-Vessiot
theories for differential and difference equations, involving the theory
of van der Put and Singer \cite{vanderPut-Singer1997}.

Let $L$ be an AS $D$-module algebra.
If $P\subset L$ is a maximal ideal, then one will see that $L_{1}:=L/P$ is a
module field over the Hopf subalgebra $D(G_{P}):=D^{1}\# RG_{P}$, where
$G_{P}$ denotes the subgroup (necessarily of finite index) of the stabilizers
of $P$. Moreover, $L$ can recover from $L_{1}$, so as
\[ L=D\otimes_{D(G_{P})}L_{1}=\bigoplus_{g\in G/G_{P}}g\otimes L_{1}, \]
where the product in $L$ recovers from the component-wise product
$(g\otimes a)(g\otimes b)=g\otimes ab$ in the last direct sum; see Section
\ref{sec-AS}. (For example, when $D=\bbc[\tau,\tau^{-1}]$ and $L=Q(A)\simeq
\bbc(\alpha)\times\bbc(\alpha)$ as above, take $P=\lan xy-1\ran\subset Q(A)$.
Then $G_{P}=\{g^{2}\;|\;g\in G\}\xrightarrow{\sim}2\bbz$ under the group
isomorphism $G\simeq\bbz$.) The $D$-invariants $L^{D}=\{a\in L\;|\; da=\vep(d)a
\quad \mbox{for all}\quad d\in D\}$ (where $\vep$ denotes the counit of $D$)
in $L$ form a subfield, such that $L^{D}\simeq L_{1}^{D(G_{P})}$.
Following \cite{Takeuchi1989},
we say that an inclusion $K\subset L$ of AS $D$-module algebras is a
{\em Picard-Vessiot (PV) extension} iff $K^{D}=L^{D}$ and there exists a
(necessarily unique) $D$-module algebra $K\subset A\subset L$ such that the
total quotient ring $Q(A)$ equals $L$, and $H:=(A\otimes_{K}A)^{D}$ generates
the left (or right) $A$-module $A\otimes_{K}A$. Then $H$ has a natural
structure of a commutative Hopf algebra over $K^{D}$ ($=L^{D}$), with which
$A/K$ is a right $H$-Galois extension; see Proposition
\ref{prop-PV-Hopf-algebra}. (In the example (\ref{intro-example-AS}),
we have $H=\bbc[z_{1},z_{2}]/\lan z_{1}^{2}z_{2}^{2}-1\ran$ with
grouplikes $z_{1}=x\otimes x^{-1}$, $z_{2}=y\otimes y^{-1}$.)
If an inclusion $K\subset L$ of AS $D$-module algebras is a PV extension,
then the induced inclusion $K/P\cap K\subset L/P$ of $D(G_{P})$-module fields
is a PV extension, where $P$ is an arbitrary maximal ideal of $L$.
The converse holds true if $G_{P}$ is normal in $G_{P\cap K}$;
see Proposition \ref{prop-interlacing2} and Theorem \ref{thm-copying}.

As our main theorems we prove: \vspace{2ex}

\noindent
{\em Galois Correspondence} (Theorem \ref{thm-Galois-correspondence1}):
Given a PV extension $L/K$ of AS $D$-module algebras, there is a 1-1
correspondence between the intermediate AS $D$-module algebras $K\subset M
\subset L$ and the Hopf ideals $I$ in the associated Hopf algebra $H$;
$L/M$ is then a PV extension with the associated Hopf algebra $H/I$
(Proposition \ref{prop-Galois-correspondence2}). This has the obvious
interpretation in terms of the affine group scheme $\bG(L/K)=\Spec H$
corresponding to $H$, and $\bG(L/K)$ is isomorphic
to the automorphism group scheme $\bAutDlinKalg(A)$ (see Section
\ref{sec-translation}).
\vspace{2ex}

\noindent
{\em Characterization} (Theorem \ref{thm-characterization}):
An inclusion $K\subset L$ of AS $D$-module algebras with $K^{D}=L^{D}$ is
a finitely generated PV extension iff $L/K$ is a minimal splitting
algebra for some $K\# D$-module $V$ of finite $K$-free rank, say $n$;
this means that $\dim_{L^{D}}\Hom_{K\#D}(V,L)=n$ and $L$ is
``minimal" with this property
(see Proposition \ref{prop-splitting-algebras2}). \vspace{2ex}

\noindent
{\em Tensor Equivalence} (Theorem \ref{equivalence-to-representations}):
If this is the case, the symmetric tensor category $\cM_{\mathrm{fin}}^{H}$
of finite-dimensional right comodules over the associated Hopf algebra $H$
(or equivalently that category $\mathrm{Rep}_{\bG(L/K)}$ of finite-dimensional
linear representations of $\bG(L/K)$) is equivalent to the abelian, rigid
tensor full subcategory $\{\{V\}\}$ ``generated" by $V$, in the tensor
category $({}_{K\# D}\cM,\otimes_{K},K)$ of $K\# D$-modules;
cf.\ \cite[Theorem 2.33]{vanderPut-Singer2003}. \vspace{2ex}

\noindent
{\em Unique Existence} (Theorem \ref{thm-existence}):
Suppose that $K^{D}$ is an algebraically closed field. For every $K\# D$-module
$V$ of finite $K$-free rank, there is a unique (up to isomorphism) minimal
splitting algebra $L/K$ which is a (finitely generated) PV extension.
\vspace{2ex}

One cannot overestimate the influence of the article \cite{Takeuchi1989} by
Takeuchi on this article of ours. Especially the main theorems above except
the third are very parallel to results in \cite{Takeuchi1989}, including their
proofs. A $C$-ferential field is equivalent to a module field over the tensor
bialgebra $T(C^{+})$.
We remark that even if $K,L$ are fields, the first two theorems above do not
imply the corresponding results in \cite{Takeuchi1989} since the notion of
$C$-ferential fields is more general than $D$-module fields in the sense of
ours. The fourth only generalizes \cite[Theorems 4.5, 4.6]{Takeuchi1989}
in which $T(C^{+})$ is supposed to be of Birkhoff-Witt type.

The last section (Section \ref{sec-Liouville}) treats the solvability
theory for Liouville extensions. The notion of Liouville extensions of
differential fields first appeared in Kolchin's historical work on the
Picard-Vessiot theory \cite{Kolchin1948}. An extension of
differential fields (of zero characteristic) is called Liouville iff it
contains no new constants and it is obtained by iterating integrations,
exponentiations, and algebraic extensions. It was shown that
a Picard-Vessiot extension of differential fields is Liouville iff
the connected component of its differential Galois group is solvable.
By the Lie-Kolchin triangularization theorem and others
\cite[Ch.\ I]{Kolchin1948}, we can characterize several types of
Liouville extensions in matrix theoretical way. For example, a Liouville
extension is obtained only by iterating integrations iff its differential
Galois group is unipotent.
For the case of an arbitrary characteristic, Okugawa \cite{Okugawa} studied the
Picard-Vessiot theory for fields with higher derivations of infinite length,
and obtained similar results on Liouville extensions.

Liouville extensions of difference fields were first
studied by Franke \cite{Franke}.
In the context of \cite{vanderPut-Singer1997},
Hendriks and Singer \cite{Hendriks-Singer}
studied on Liouville solutions of difference equations with rational function
coefficients. They defined the notion of ``Liouvillian sequences"
and showed that a linear difference equation can be solved in terms of
such sequences iff the difference Galois group is solvable.
(Moreover, they gave an algorithm to find such Liouville solutions,
using Petkov\v{s}ek's algorithm \cite{Petkovsek}.)

In the last section, we define the notion of Liouville extensions of
AS $D$-module algebras and prove a solvability theorem in the unified context.

When we study Liouville extensions in terms of affine group schemes,
we will meet the following difficulty: the Lie-Kolchin triangularization
theorem can not be extended generally to affine group schemes
(see \cite[Ch.\ 10]{Waterhouse}).
Certainly there are gaps between the triangulability
and the connected solvability, even if the base field is algebraically closed.
So we need some intermediate notions and have to study how they are related
each other. In Section \ref{secLGS}, we define ``Liouville group schemes"
so that it is suitable for Liouville extensions defined later,
and study how strong the definition is.
An algebraic affine group scheme $\bG$ over a field
$k$ is called {\em ($k$-)Liouville} (cf.\ \cite[p.\ 374]{Kolchin1973})
iff there exists a normal chain of closed subgroup schemes
$\bG=\bG_{0}\vtr\bG_{1}\vtr\dotsb\vtr\bG_{r}=\{1\}$
such that each $\bG_{i-1}/\bG_{i}$ ($i=1,\dotsc,r$) is at least one of the
following types: finite etale, a closed subgroup scheme of $\Ga$,
or a closed subgroup scheme of $\Gm$. When $k$ is algebraically closed,
$\bG$ is Liouville iff the connected component $\bG^{\circ}$ is solvable
(Proposition \ref{solvability2}). But in general it does not holds; we show
this fact by examples. For connected affine group schemes, we will see the
condition to be Liouville is properly stronger than the solvability but
weaker than the triangulability.

Let $L\supset K$ be an inclusion of AS $D$-module algebras. For finitely many
elements $x_{1},\dotsc,x_{n}\in L$, let $K\lan x_{1},\dotsc,x_{n}\ran$ denote
the smallest AS $D$-module subalgebra in $L$ including both $K$ and $x_{1},
\dotsc,x_{n}$. $L/K$ is called a {\em $\Ga$-primitive extension} (resp., a
{\em $\Gm$-primitive extension}) iff there is an $x\in L$ such that
$d(x)\in K$ for all $d\in D^{+}=\Ker\varepsilon$
(resp., $x$ is a non-zero divisor of $L$ (which is necessarily invertible) and
$d(x)x^{-1}\in K$ for all $d\in D$) and $L=K\lan x\ran$. We say that
$L/K$ is a {\em finite etale extension} iff $L$ is a separable
$K$-algebra in the sense of \cite{DeMeyer-Ingraham}, i.e.\ $L$ is a projective
$L\otimes_{K}L$-module.
Then we define {\em Liouville extension} as such a
finitely generated extension $L/K$ that $L^{D}=K^{D}$ and there exists
a sequence of AS $D$-module algebras
$K=L_{0}\subset L_{1}\subset\dotsb\subset L_{r}=L$
such that each $L_{i}/L_{i-1}$ ($i=1,\dotsc,r$) is at least one of the
following types: $\Ga$-primitive extension, $\Gm$-primitive extension,
or finite etale extension. As the last one of the main theorems, we will show
the following:
\vspace{2ex}

\noindent
{\em Solvability} (Theorem \ref{mainTh}):
Let $L/K$ be a finitely generated PV extension.
Then the following are equivalent:
\begin{enumerate}
\renewcommand{\labelenumi}{\rm(\alph{enumi})}
\item $L/K$ is a Liouville extension.
\item There exists a Liouville extension $F/K$ such that $L\subset F$.
\item $\bG(L/K)$ is Liouville.
\end{enumerate}
When $k$ is algebraically closed, these are equivalent to:
\begin{enumerate}
\renewcommand{\labelenumi}{\rm(d)}
\item $\bG(L/K)^{\circ}$ is solvable.
\end{enumerate}

Moreover we will characterize ten types of Liouville extensions
just being compatible with \cite[\S 24--27]{Kolchin1948};
see Definition \ref{defLE}, Corollary \ref{precisely} and its following
paragraph.

\subsection*{Conventions}

Throughout this part, we always work over an arbitrarily
fixed field $R$. All vector spaces, algebras and coalgebras are defined 
at least over $R$. All algebras are associative and have the identity element.
All modules over an algebra are unital, left modules unless otherwise stated.
All separable algebras are taken in the sense of \cite{DeMeyer-Ingraham};
see also \cite[Ch.\ 6]{Waterhouse}.

The notation $\Hom_{\cR}$ (resp.\ $\End_{\cR}$) with a ring $\cR$ always
denotes the set of all $\cR$-linear maps (resp.\ $\cR$-linear endomorphisms),
but the unadorned $\Hom$ may indicate group homomorphisms or homomorphisms
of group schemes. Algebra (resp.\ coalgebra) maps are always dentoed by
$\Alg$ (resp.\ $\Coalg$). The notation $\Aut$ indicates
automorphisms in some sense; for example, $\AutDlinKalg$ means $D$-linear and
$K$-algebra automorphisms. $\bAut$ in the bold style indicates an associated
group functor as in \cite[(7.6)]{Waterhouse}.
Coalgebra structures are denoted by $(\Delta,\varepsilon)$.
If we need to specify a coalgebra (or a coring) $C$,
the notation $(\Delta_{C},\varepsilon_{C})$ is also used.
For a coalgebra $C$, $C^{+}$ denotes $\Ker\vep$.
The antipode of a Hopf algebra is denoted by $S$.
We use the {\em sigma notation} (see \cite[\S 1.2, pp.\,10--11]{Sweedler1969}
or \cite[\S 1.4, pp.\ 6--7]{Montgomery}):
\[ \Delta(c)=\sum_{(c)}c_{(1)}\otimes c_{(2)} \qquad \mbox{etc.} \]
When $(M,\lambda)$ is a right (resp.\ left) $C$-comodule, $\lambda(m)$
($m\in M$) is denoted by the sigma notation
\[ \lambda(m)=\sum_{(m)}m_{(0)}\otimes m_{(1)}\in M\otimes_{R}C\quad
\mbox{(resp.}\ \ \lambda(m)=\sum_{(m)}m_{(-1)}\otimes m_{(0)}
\in C\otimes_{R}M). \]

By ``a symmetric tensor category $(\gA,\otimes,I)$" we mean that
$(\gA,\otimes)$ is a symmetric tensor (or monoidal) category
\cite[\S 10.4, p.\ 199]{Montgomery} with a fixed unit object $I$.
We can define algebras, coalgebras, etc., in $(\gA,\otimes,I)$ naturally
by commutative diagrams. For an algebra $A$ in $(\gA,\otimes,I)$,
left $A$-modules (resp.\ right $A$-modules, resp.\ $(A,A)$-bimodules)
in $(\gA,\otimes,I)$ can also be defined and the category of them
is denoted by ${}_{A}\gA$ (resp.\ $\gA_{A}$, resp.\ ${}_{A}\gA_{A}$).
For a ring $\cR$, ${}_{\cR}\cM$ (resp.\ $\cM_{\cR}$) denotes the category of
left (resp.\ right) $\cR$-modules. For a coalgebra $C$, $\cM^{C}$
(resp.\ ${}^{C}\cM$) denotes the category of right (resp.\ left) $C$-comodules.
Moreover, further notations, such as ${}_{A}\cM^{H}$, ${}_{A}^{H}\cM$, etc.,
which indicate categories of relative Hopf modules are used as in
\cite[\S 8.5]{Montgomery}.

\section{Basic notions and results on $D$-module algebras}
\label{sec-basic-notions}

Let $D$ be a cocommutative bialgebra.
An algebra $A$ is called a {\em $D$-module algebra}
(see \cite[\S 7.2, p.\,153]{Sweedler1969} or \cite[\S 4.1]{Montgomery}) iff
$A$ has a $D$-module structure $D\otimes_{R}A\rightarrow A$ and the action of
$D$ measures $A$ to $A$. The last condition means that
\[ \rho_{A} : A \rightarrow \Hom_{R}(D,A),\quad
a \mapsto [d\mapsto da] \]
is an algebra map, where $\Hom_{R}(D,A)$ is considered as an
algebra with the convolution product (see \cite[pp.\,69--70]{Sweedler1969}
or \cite[\S 1.4]{Montgomery});
or in other words, using the sigma notation,
\[ d(ab)=\sum_{(d)}(d_{(1)}a)(d_{(2)}b),\qquad d(1)=\varepsilon(d)1 \]
hold for all $d\in D$ and $a,b\in A$. Throughout this part, we assume
$D$-module algebras are commutative unless otherwise stated.
Note that the algebra $\Hom_{R}(D,A)$ is commutative in our situation, and
it has a $D$-module algebra structure given by
\begin{equation}
\label{eq-convolution}
(d\varphi)(c)=\varphi(cd)\qquad (c,d\in D,\ \varphi\in\Hom_{R}(D,A)).
\end{equation}
One sees $\rho_{A}$ is an injective $D$-module algebra map.
For a $D$-module algebra $A$, the {\em smash product} $A\# D$ means the algebra
which is $A\otimes_{R}D$ with the semi-direct product:
\[ (a\# c)(b\# d)=\sum_{(c)}a(c_{(1)}b)\# c_{(2)}d \]
(see \cite[pp.\ 155--156]{Sweedler1969} or \cite[\S 4.1]{Montgomery}).
For $A\# D$-modules $V,W\in{}_{A\# D}\cM$, we have an $A\# D$-module
structure on $V\otimes_{A}W$ given by
\[ (a\#d)(v\otimes w)=a\sum_{(d)}d_{(1)}v\otimes d_{(2)}w\qquad
(a\in A,\ d\in D,\ v\in V,\ w\in W). \]
Thus we have an abelian symmetric tensor category
\cite[Definition 1.15]{Deligne-Milne}
$({}_{A\# D}\cM,\otimes_{A},A)$ with the canonical symmetry
$V\otimes_{A}W\rightarrow W\otimes_{A}V$, $v\otimes w\mapsto w\otimes v$.
For a $D$-module $V$,
\[ V^{D}:=\{v\in V\;|\;dv=\varepsilon(d)v\quad \mbox{for all}\quad
d\in D\} \]
is called the {\em constants} (or the {\em $D$-invariants}) of $V$.
Especially $A^{D}$ becomes an algebra. If $V$ is an $A\# D$-module,
we see $\Hom_{A\# D}(A,V)\xrightarrow{\sim}V^{D}$, $\varphi\mapsto
\varphi(1)$ is an $A^{D}$-module isomorphism and in particular
$\End_{A\# D}(A)\simeq A^{D}$ as algebras. The functor $(-)^{D} : {}_{A\# D}\cM
\rightarrow {}_{A^{D}}\cM$ is an exact functor since $A$ is a projective
$A\# D$-module (indeed, $A\# D\simeq A\oplus(A\otimes_{R}D^{+})$ as
$A\# D$-modules via $a\# d\mapsto (a\vep(d),a\otimes(d-\vep(d)))$). 
Let $B$ be a $D$-module algebra including $A$ as a $D$-module subalgebra,
$V$ an $A\# D$-module, and $W$ a $B\# D$-module.
If $D$ is a Hopf algebra with the antipode $S$, then $\Hom_{A}(V,W)$ has a
$B\# D$-module structure given by the $D$-conjugation:
\begin{equation}
\label{eq-internalHom}
((b\# d)\varphi)(v)=b\sum_{(d)}d_{(1)}(\varphi(S(d_{(2)})v)\qquad
(v\in V) 
\end{equation}
for $b\in B$, $d\in D$, and $\varphi\in\Hom_{A}(V,W)$;
see \cite[Proposition 1.8]{Takeuchi1989}. We see
$\Hom_{A}(V,W)^{D}=\Hom_{A\# D}(V,W)$. Especially $\Hom_{A}$
is an internal Hom of $({}_{A\# D}\cM,\otimes_{A},A)$ in such a case
(see \cite[p.\ 109]{Deligne-Milne} for the definition of internal Hom).

The following proposition, like Schur's lemma, is very important:
\begin{proposition}
\label{prop-simple-objects}
Let $\gA$ be an abelian subcategory of ${}_{\bbz}\cM$ (the category of additive
groups) whose inclusion functor is additive.
An object $X$ in $\gA$ is simple iff
\begin{enumerate}
\renewcommand{\labelenumi}{\rm(\alph{enumi})}
\item the endomorphism ring $E:=\End_{\gA}(X)$ is a division ring, and
\item for every object $Y$ in $\gA$, the evaluation map
\[ \mathrm{ev} : \Hom_{\gA}(X,Y)\otimes_{E}X\rightarrow Y \]
is injective.
\end{enumerate}
\end{proposition}
\begin{proof}
(``If" part.) Let $Y$ be a subobject of $X$.
Since $\Hom_{\gA}(X,Y)$ is a right ideal of $E=\End_{\gA}(X)$,
$\Hom_{\gA}(X,Y)$ equals $0$ or $E$ by (a).
If $\Hom_{\gA}(X,Y)=E$, we have $Y=X$ since $X\simeq E\otimes_{E}X\rightarrow
Y$ is injective.
If $\Hom_{\gA}(X,Y)=0$, then $E\rightarrow \Hom_{\gA}(X,X/Y)$
is injective.
Since all $E$-modules are flat, we have that $X\simeq E\otimes_{E}X
\rightarrow\Hom_{\gA}(X,X/Y)\otimes_{E}X\rightarrow X/Y$ is injective
and hence $Y=0$.

(``Only if" part.) (a) For $0\neq f\in E$, $\im(f)$ is a non-zero
subobject of $X$ and hence $\im(f)=X$, $\Ker(f)=0$.
Thus $f$ is invertible.

(b) Since $X$ is simple, each $0\neq f\in\Hom_{\gA}(X,Y)$ is monic and
hence $\im(f)$ is isomorphic to $X$. It suffices to prove that the sum
$\sum_{i=1}^{r}\im(f_{i})\subset Y$ is direct if $f_{1},\dotsc,f_{r}$ are
$E$-linearly independent in $\Hom_{\gA}(X,Y)$. To prove this,
we shall use induction on $r$. When $r=1$, the assertion is clear.
When $r>1$, suppose that the assertion is true for $\{f_{1},\dotsc,f_{r-1}\}$.
Seeking a contradiction, assume $\im(f_{r})\cap \sum_{i=1}^{r-1}
\im(f_{i})\neq 0$. Since $\im(f_{r})$ is simple,
we have $\im(f_{r})\subset \bigoplus_{i=1}^{r-1}\im(f_{i})$.
Then there exist $\varphi_{i}\in E$ ($i=1,\dotsc,r-1$) such that
the diagram
\[ \begin{array}{c} \begin{CD}
\im(f_{r}) @>{\mbox{\scriptsize inclusion}}>>
{\displaystyle\bigoplus_{i=1}^{r-1}}\im(f_{i})
@>{\mbox{\scriptsize projection}}>> \im(f_{i}) \\
 @A{f_{r}}AA @. @AA{f_{i}}A
\end{CD} \\
\begin{CD}
 X @>{\hspace{7em} \varphi_{i} \hspace{7em}}>> X
\end{CD} \end{array} \vspace{2ex} \]
commutes for $i=1,\dotsc,r-1$ since $f_{i} : X\rightarrow \im(f_{i})$
is invertible. Then $f_{r}=f_{1}\circ\varphi_{1}+\dotsb +f_{r-1}\circ
\varphi_{r-1}$. This contradicts that $f_{1},\dotsc,f_{r}$ are
$E$-linearly independent.
\end{proof}

\begin{remark}
I heard the above proposition from Professor A.\ Masuoka.
Though it seems well-known, an explicit citation was not found as far as
I searched. It is said that Professor Masuoka knew this by a comment from
Professor T.\ Brezi\'nski
on \cite[Theorem 1.1 and the Theorem on p.\ 232]{Masuoka-Yanai}.
\end{remark}

\begin{definition}
A $D$-module algebra $A$ is called {\em simple} iff it is simple in
${}_{A\# D}\cM$, i.e.\ $A$ has no non-trivial $D$-stable ideal.
\end{definition}

The next corollary follows immediately from Propositoin
\ref{prop-simple-objects}.
\begin{corollary}
\label{cor-simple-module-algebras}
A $D$-module algebra $A$ is simple iff
\begin{enumerate}
\renewcommand{\labelenumi}{\rm(\alph{enumi})}
\item $A^{D}$ is a field, and
\item for every $A\# D$-module $Y$, the map
\[ \begin{array}{lll}
 & Y^{D}\otimes_{A^{D}}A\rightarrow Y, & y\otimes a\mapsto ay \vspace{1ex} \\
\mbox{(or} & A\otimes_{A^{D}}Y^{D}\rightarrow Y, & a\otimes y\mapsto ay
\mbox{)}
\end{array} \]
is injective.
\end{enumerate}
\end{corollary}
\begin{proof}
Recall that $\End_{A\# D}(A)\simeq A^{D}$ and $\Hom_{A\# D}(A,Y)\simeq Y^{D}$.
The evaluation map is identified with the map in (b) above.
\end{proof}

Let $A$ be a $D$-module algebra and $\rho_{A} : A\rightarrow\Hom_{R}(D,A)$
the associated algebra map. Then $\Hom_{R}(D,A)$ has two $A$-module
structures $A\otimes_{R}\Hom_{R}(D,A)\rightarrow\Hom_{R}(D,A)$, given by
(I) $a\otimes\varphi\mapsto\rho_{A}(a)\ast\varphi=\varphi\ast\rho_{A}(a)$,
and (II) $a\otimes\varphi\mapsto a\varphi=[d\mapsto a\varphi(d)]$.
The structure given by (II)
can be considered through the following algebra isomorphism:
\[ \sigma : A\xrightarrow{\sim}\Hom_{R}(D,A)^{D},\quad a\mapsto a\vep_{D}, \]
whose inverse is $\varphi\mapsto\varphi(1)$. Here we are
taking the $D$-module structure on $\Hom_{R}(D,A)$ in the sense of
(\ref{eq-convolution}). As in \cite[Corollary 1.4]{Takeuchi1989}, the next
lemma follows from Corollary \ref{cor-simple-module-algebras}.
\begin{lemma}
\label{lem-beta-injection}
If $A$ is simple, then the following map:
\[ \beta : A\otimes_{A^{D}}A\rightarrow\Hom_{R}(D,A),\quad
a\otimes b\mapsto\rho_{A}(b)\ast(a\vep_{D})=a\rho_{A}(b) \]
is a two-sided $A$-linear (left through $\sigma$, right through $\rho_{A}$)
injection.
\end{lemma}
\begin{proof}
Consider $Y=\Hom_{R}(D,A)$ as an $A\# D$-module by the $A$-module structure
given by (I) and by the $D$-module structure in the sense of
(\ref{eq-convolution}):
\[ (a\# d)\varphi=\rho_{A}(a)\ast(d\varphi)\qquad
(a\in A,\ d\in D,\ \varphi\in Y). \]
This is well-defined:
\[ \begin{array}{ll}
 & (a'\# d')((a\# d)\varphi))=\rho_{A}(a')\ast{\displaystyle
\sum_{(d')}}(d_{(1)}'\rho_{A}(a))\ast(d_{(2)}'d\varphi) \vspace{1ex} \\
= & \rho_{A}(a')\ast{\displaystyle\sum_{(d')}}\rho_{A}(d_{(1)}'a)
\ast(d_{(2)}'d\varphi)={\displaystyle\sum_{(d')}}\rho_{A}(a'(d_{(1)}'a))
\ast(d_{(2)}'d\varphi) \vspace{1ex} \\
= & ((a'\# d')(a\# d))\varphi.
\end{array} \]
Then $\beta$ is injective since $A\simeq Y^{D}$ through $\sigma$.
\end{proof}

This lemma has an application when one needs to think of the $A^{D}$-linear
dependence of elements in $A$. Takeuchi generalized the Wronskian
(and Casoratian) criterion as follows:

\begin{proposition}(\cite[Proposition 1.5]{Takeuchi1989})
\label{Wronskian}
Let $K$ be a $D$-module field. Then $a_{1},\dotsc,a_{n}\in K$ are
$K^{D}$-linearly independent iff there exist $h_{1},\dotsc,h_{n}\in D$
such that $\det(h_{i}(a_{j}))_{i,j}\neq 0$.
\end{proposition}
\begin{proof}
We include the proof for convenience.

(``If" part.) If $\sum_{j=1}^{n}c_{j}a_{j}=0$ ($c_{1},\dotsc,c_{n}\in K^{D}$),
then $\sum_{j=1}^{n}c_{j}h_{i}(a_{j})=0$ for $i=1,\dotsc,n$. Since the
matrix $(h_{i}(a_{j}))_{i,j}$ is invertible, we have $c_{1}=\dotsb =c_{n}=0$.

(``Only if" part.) Put $W=K^{D}a_{1}+\dotsb +K^{D}a_{n}$, an $n$-dimensional
$K^{D}$-vector subspace of $K$. Consider the $K$-linear injection
\[ \beta : K\otimes_{K^{D}}W\rightarrow\Hom_{R}(D,K)\simeq
\Hom_{K}(K\otimes_{R}D,K), \quad b\otimes w\mapsto [d\mapsto b(dw)] \]
to which $\beta$ in Lemma \ref{lem-beta-injection} restricts.
Let $\{d_{\alpha}\}_{\alpha\in\Lambda}$ be an $R$-basis of $D$ and
$d_{\alpha}^{\vee}$ be the dual of $d_{\alpha}$ in $\Hom_{R}(D,K)$.
Then $\Hom_{K}(K\otimes_{R}D,K)=\prod_{\alpha\in\Lambda}Kd_{\alpha}^{\vee}$
as a $K$-vector space.
Notice that $\rho_{A}(a)=\sum_{\alpha\in\Lambda}(d_{\alpha}a)
d_{\alpha}^{\vee}$ ($a\in K$). Since $\rho_{A}(a_{1}),\dotsc,\rho_{A}(a_{n})$
are $K$-linearly independent and since $K$ is a field, we obtain a $K$-basis
$v_{1},\dotsc,v_{n}$ of $K\otimes_{K^{D}}W$ such that
\begin{eqnarray*}
\beta(v_{1}) & = & h_{1}^{\vee}
{}+\sum_{\substack{\alpha\in\Lambda \\ d_{\alpha}\neq h_{1}}}c_{1,\alpha}
d_{\alpha}^{\vee} \\
\beta(v_{2}) & = & \qquad h_{2}^{\vee}
{}+\sum_{\substack{\alpha\in\Lambda \\ d_{\alpha}\neq h_{1},h_{2}}}c_{2,\alpha}
d_{\alpha}^{\vee} \\
\vdots & &  \\
\beta(v_{n}) & = & \qquad\qquad h_{n}^{\vee}
{}+\sum_{\substack{\alpha\in\Lambda \\ d_{\alpha}\neq h_{1},\dotsc,h_{n}}}
c_{n,\alpha}d_{\alpha}^{\vee}
\end{eqnarray*}
for some $n$ elements $h_{1},\dotsc,h_{n}\in\{d_{\alpha}\}_{\alpha\in
\Lambda}$ by sweeping-out. Consider the transposed $K$-linear map of $\beta$:
\[ \gamma : K\otimes_{R}D\rightarrow\Hom_{K}(K\otimes_{K^{D}}W,K),
\quad a\otimes d\mapsto [b\otimes w\mapsto ab(dw)]. \]
Let $v_{1}^{\vee},\dotsc,v_{n}^{\vee}$ be the dual basis of
$v_{1},\dotsc,v_{n}$.
Then we have ${}^{t}(\gamma(h_{1}),\dotsc,\gamma(h_{n}))
=T{}^{t}(v_{1}^{\vee},\dotsc,v_{n}^{\vee})$ with a strictly lower triangular
matrix $T\in GL_{n}(K)$.
Thus $\gamma(h_{1}),\dotsc,\gamma(h_{n})$ form a $K$-basis for
$\Hom_{K}(K\otimes_{K^{D}}W,K)\simeq\Hom_{K^{D}}(W,K)$. The $K$-isomorphism
\[ \begin{array}{rcccl}
K^{n} & \xrightarrow{\sim} & \Hom_{K^{D}}(W,K) & \xrightarrow{\sim} &
K^{n} \\
(c_{i})_{i} & \mapsto & {\displaystyle\sum_{i=1}^{n}}c_{i}
\gamma(h_{i}) & \mapsto &
\left({\displaystyle\sum_{i=1}^{n}}c_{i}h_{i}(a_{j})\right)_{j}
\end{array} \]
is precisely the right multiplication of matrix $(h_{i}(a_{j}))_{i,j}$.
It follows that the matrix has an inverse.
\end{proof}

\begin{remark}
In the above proof, we see that $\Ker\gamma$ is a left ideal of $K\#D$.
If $D=R[\partial]$ with one primitive $\partial$, then $\Ker\gamma$
is generated by a monic differential operator of order $n$. Thus we can
take $h_{1}=1,h_{2}=\partial,\dotsc,h_{n}=\partial^{n-1}$ in such a case.
Namely we have the Wronskian criterion in the usual sense.
Similarly we also have the ordinary Casoratian criterion for difference fields.
\end{remark}

\section{Tensor equivalences associated to Hopf subalgebras}
\label{sec-relative-Hopf-modules}

In what follows we assume that $D$ is a cocommutative Hopf algebra.
Let $C$ be a Hopf subalgebra of $D$.
A coalgebra in the tensor category $({}_{D}\cM,\otimes_{R},R)$ is called a
{\em $D$-module coalgebra}. $D$ is a $D$-module coalgebra, and
$\bar{D}:=D/DC^{+}$ is its quotient.
The $R$-abelian category ${}^{\bar{D}}_{D}\cM$ of left $(\bar{D},D)$-Hopf
modules is defined as follows (see \cite[pp.\ 454--455]{Takeuchi1979} or
\cite[\S 8.5]{Montgomery}): \\

\noindent
\underline{Objects}. An object of ${}^{\bar{D}}_{D}\cM$
is a left $D$-module $M$ which is also a left $\bar{D}$-comodule with a
structure $\lambda_{M}$, say, such that
\[ \lambda_{M}(dm)=\Delta(d)\lambda_{M}(m)=\sum_{(d)}\sum_{(m)}d_{(1)}m_{(-1)}
\otimes d_{(2)}m_{(0)}\in \bar{D}\otimes_{R}M \]
for all $d\in D$ and $m\in M$. \\

\noindent
\underline{Morphisms}. Morphisms of ${}^{\bar{D}}_{D}\cM$ are $D$-module
and $\bar{D}$-comodule maps. \\

Given objects $M,N$ in ${}^{\bar{D}}_{D}\cM$, let
$M\Box_{\bar{D}}N$ denote the cotensor product; this is by definition the
equalizer of the two $\bar{D}$-colinear maps 
\[ M\otimes N\substack{\longrightarrow \\ \longrightarrow}
\bar{D}\otimes M\otimes N \]
given by the structure maps of $M,N$, or in other words,
\[ M\Box_{\bar{D}}N
=\left\{\sum_{i}x_{i}\otimes y_{i}\right.\;\left|\;
\sum_{i}\sum_{(x_{i})}(x_{i})_{(-1)}\otimes (x_{i})_{(0)}\otimes y_{i}
=\sum_{i}\sum_{(y_{i})}(y_{i})_{(-1)}\otimes x_{i}\otimes(y_{i})_{(0)}\right\}.
\]
This is a $D$-submodule of $M\otimes N$, and is further an object in
${}^{\bar{D}}_{D}\cM$. We see that ${}^{\bar{D}}_{D}\cM =({}^{\bar{D}}_{D}\cM,
\Box_{\bar{D}},\bar{D})$ is a symmetric tensor category.
Indeed, the associativity $(M\Box_{\bar{D}}N)\Box_{\bar{D}}L\xrightarrow{\sim}
M\Box_{\bar{D}}(N\Box_{\bar{D}}L)$ and the symmetry
$M\Box_{\bar{D}}N\xrightarrow{\sim}N\Box_{\bar{D}}M$
are induced by those of $({}_{D}\cM,\otimes_{R})$ naturally.
We have isomorphisms
\[ \begin{array}{l}
M\Box_{\bar{D}}\bar{D}\xrightarrow{\sim}M,\quad
{\displaystyle\sum_{i}}m_{i}\otimes a_{i}\mapsto
{\displaystyle\sum_{i}}m_{i}\varepsilon(a_{i}), \vspace{1.5ex} \\
\bar{D}\Box_{\bar{D}}M\xrightarrow{\sim}M,\quad
{\displaystyle\sum_{i}}a_{i}\otimes m_{i}\mapsto
{\displaystyle\sum_{i}}\varepsilon(a_{i})m_{i},
\end{array} \]
whose inverses are obtained by $\lambda_{M}$. Thus $\bar{D}$ is a unit object.

For an object $V$ in ${}_{C}\cM$, define
\[ \Phi(V) = D\otimes_{C}V. \]
This is naturally an object in ${}^{\bar{D}}_{D}\cM$. We thus have an
$R$-linear functor
\[ \Phi : {}_{C}\cM\rightarrow {}^{\bar{D}}_{D}\cM. \]

\begin{proposition}
\label{prop-Hopf-descent}
$\Phi$ is an equivalence of symmetric tensor categories.
\end{proposition}
\begin{proof}
By \cite[Theorem 2 and 4]{Takeuchi1979}, $\Phi$ is a category equivalence; its
quasi-inverse $N\mapsto\Psi(N)$ is given by
\[ \Psi(N)=\{n\in N\;|\;\lambda_{N}(n)=\bar{1}\otimes n\ \ \mbox{in}\ \ \bar{D}
\otimes N\}, \]
where $\lambda_{N} : N\rightarrow \bar{D}\otimes N$ is the structure map
on $N$. It is easy to see that
\[ \begin{array}{l}
\Psi(M)\otimes\Psi(N)\rightarrow\Psi(M\Box_{\bar{D}}N),\quad m\otimes n
\mapsto m\otimes n, \vspace{1ex} \\
R\rightarrow\Psi(\bar{D}),\quad 1\mapsto\bar{1}
\end{array} \]
are isomorphisms in ${}_{C}\cM$. We see that the isomorphisms, as tensor
structures, make $\Psi$ an equivalence of symmetric tensor categories.
\end{proof}

Let $D^{1}$ denote the irreducible component in $D$ containing 1; this is
the largest irreducible Hopf subalgebra. If the characteristic $\ch R$ of $R$
is zero, then $D^{1}=U(\g)$, the universal envelope of the Lie algebra $\g
=P(D)$ of all primitives in $D$; see \cite[Ch.\ XIII]{Sweedler1969} or
\cite[\S 5.6]{Montgomery}. Let $G=G(D)$ denote the group of all grouplikes
in $D$. 

In what follows we suppose:
\begin{assumption}
\label{assum-pointed}
$D$ is pointed, so that
\[ D = D^{1}\# RG, \]
the smash product with respect to the conjugate action by $G$ on $D^{1}$;
see \cite[Theorem 8.1.5]{Sweedler1969} or \cite[Corollary 5.6.4]{Montgomery}.
\end{assumption}

In the following, we take as $C$ a Hopf subalgebra of the form
\[ C=D(G_{1}):=D^{1}\# RG_{1}, \]
where $G_{1}\subset G$ is a subgroup {\em of finite index}.
The equivalence $\Phi$ will be denoted by
\begin{equation}
\label{eq-copying-functor}
\Phi_{G_{1}} : {}_{D(G_{1})}\cM\xrightarrow{\approx}{}^{\bar{D}}_{D}\cM,
\end{equation}
if one needs to specify $G_{1}$.

The vector space $R(G/G_{1})$ freely spanned by the set $G/G_{1}$ of left 
cosets is a quotient left $D$-module coalgebra of $D$ along the map
$D=D^{1}\# RG\rightarrow R(G/G_{1})$ which is given by the counit
$\varepsilon : D^{1}\rightarrow R$ and the natural projection $G\rightarrow
G/G_{1}$. Since the map induces an isomorphism $\bar{D}\xrightarrow{\sim}
R(G/G_{1})$, left $\bar{D}$-comodules are identified with $(G/G_{1})$-graded
vector spaces: for $N\in {}^{\bar{D}}\cM$,
\[ N=\bigoplus_{s\in G/G_{1}}N_{s}\qquad
(N_{s}=\{n\in N\;|\;\lambda_{N}(n)=s\otimes n\}). \]
An object in ${}^{\bar{D}}_{D}\cM$ is a left $D$-module
$N=\bigoplus_{s\in G/G_{1}}N_{s}$ which satisfy
that $gN_{s}\subset N_{gs}$ ($g\in G$, $s\in G/G_{1}$).
If $M=\bigoplus_{s\in G/G_{1}}M_{s}$ is another object in
${}^{\bar{D}}_{D}\cM$, then
\[ M\Box_{\bar{D}}N=\bigoplus_{s\in G/G_{1}}M_{s}\otimes N_{s}. \]
We have $D=\bigoplus_{g\in G/G_{1}}gD(G_{1})$.

\begin{notation}
\label{notation-left-cosets}
Here and in what follows, $g\in G/G_{1}$ means that $g$ lies
in a fixed system of those representatives in $G$ for the left cosets
$G/G_{1}$ which include the neutral element $1$ in $G$.
\end{notation}

The neutral component $N_{1}=\Psi(N)$ in $N$ is a $D(G_{1})$-submodule.
We have the identification
\[ \Phi(N_{1})=\bigoplus_{g\in G/G_{1}}g\otimes N_{1}. \]
Here $D$ acts on the right-hand side so that if $d\in D^{1}$,
\[ d(g\otimes n)=g\otimes (g^{-1}dg)n\quad (n\in N_{1}), \]
and if $h\in G$,
\[ h(g\otimes n)=g'\otimes tn\quad (n\in N_{1}), \]
where $g'$ is a representative and $t\in G_{1}$ such that $hg=g't$.
Hence, by Proposition \ref{prop-Hopf-descent}, we have an isomorphism
$\Phi(N_{1})=\bigoplus_{g\in G/G_{1}}g\otimes N_{1}\xrightarrow{\sim}N$
in ${}^{\bar{D}}_{D}\cM$, given by $g\otimes n\mapsto gn$.

An algebra $A$ in $({}^{\bar{D}}_{D}\cM,\Box_{\bar{D}},\bar{D})$
is precisely such a $D$-module algebra that is the direct product
$\prod_{s\in G/G_{1}}A_{s}$
of $D^{1}$-module algebras $A_{s}$ ($s\in G/G_{1}$), satisfying $gA_{s}
\subset A_{gs}$ ($g\in G$). It is identified with $\Phi(A_{1})
=\bigoplus_{g\in G/G_{1}}g\otimes A_{1}$, which is endowed with the
component-wise product. We observe that $e_{g}=g\otimes 1\in\Phi(A_{1})$ are
orthogonal central idempotents.

Let $A=\Phi(A_{1})$ be as above. An $A_{1}$-module $V$ in ${}_{D(G_{1})}\cM$
is precisely a module over the algebra $A_{1}\# D(G_{1})$ of smash product:
${}_{A_{1}}({}_{D(G_{1})}\cM)={}_{A_{1}\# D(G_{1})}\cM$.
$\Phi(V)$ is naturally an $A$-module in ${}^{\bar{D}}_{D}\cM$; this is in
particular an $A\# D$-module.

\begin{proposition}
\label{prop-copying-functor}
Let $A_{1}$ be a $D(G_{1})$-module algebra and
$A=\Phi(A_{1})$. The functor
\[ \Phi : {}_{A_{1}\# D(G_{1})}\cM\rightarrow {}_{A\# D}\cM \]
is an equivalence of $R$-abelian categories.
\end{proposition}
\begin{proof}
By Proposition \ref{prop-Hopf-descent}, it suffices to prove that the category
${}_{A}({}^{\bar{D}}_{D}\cM)$ of $A$-modules in
${}^{\bar{D}}_{D}\cM$ is isomorphic to ${}_{A}({}_{D}\cM)={}_{A\# D}\cM$.
Given $N$ in ${}_{A\# D}\cM$,
define $N_{g}=e_{g}N$ ($e_{g}=g\otimes 1\in A=\Phi(A_{1})$, $g\in G/G_{1}$).
Then $N=\bigoplus_{g\in G/G_{1}}N_{g}$ so that $N$ is in
${}_{A}({}^{\bar{D}}_{D}\cM)$. This gives the desired isomorphism.
\end{proof}

This proposition can be extended as follows:
\begin{proposition}
\label{prop-copying-functor2}
Let $A=\Phi(A_{1})$ be as above. The functor
\[ \Phi : ({}_{A_{1}}({}_{D(G_{1})}\cM)_{A_{1}},\otimes_{A_{1}},A_{1})
\rightarrow ({}_{A}({}_{D}^{\bar{D}}\cM)_{A},\otimes_{A},A) \]
is a tensor equivalence.
\end{proposition}
\begin{proof}
For $V,W\in {}_{A_{1}}({}_{D(G_{1})}\cM)_{A_{1}}$, we easily see
\[ \Phi(V)\otimes_{A}\Phi(W)\simeq\sum_{g\in G/G_{1}}g\otimes(V\otimes_{A_{1}}
W)=\Phi(V\otimes_{A_{1}}W) \]
in $({}_{A}({}_{D}^{\bar{D}}\cM)_{A},\otimes_{A},A)$.
\end{proof}

We see that the functor $\Phi$ preserves constants and simple module
algebras:
\begin{lemma}
\label{lem-preservation}
{\rm (i)} Let $V$ be a $D(G_{1})$-module. Then an isomorphism $V^{D(G_{1})}
\xrightarrow{\sim}\Phi(V)^{D}$ is given by
$v\mapsto\sum_{g\in G/G_{1}}g\otimes v$.

{\rm (ii)} Let $A_{1}$ be a commutative $D(G_{1})$-module algebra. Then
$A_{1}$ is a simple $D(G_{1})$-module algebra iff $\Phi(A_{1})$ is a simple
$D$-module algebra.
\end{lemma}
\begin{proof}
(i) If $\sum_{g}g\otimes v_{g}\in\Phi(V)^{D}$, one sees first $v_{1}\in V^{C}$,
and then $v_{g}=v_{1}$ for all $g\in G/G_{1}$.

(ii) This directly follows from Proposition \ref{prop-copying-functor}.
\end{proof}

\section{Artinian simple $D$-module algebras}
\label{sec-AS}

Let $D=D^{1}\# RG$ be a cocommutative pointed Hopf algebra as in
the previous section. In what follows we further assume:
\begin{assumption}
\label{assum-BW-type}
The irreducible Hopf algebra $D^{1}$ is of Birkhoff-Witt type.
\end{assumption}

This means that every primitive element of $D$ lies in a divided power
sequence of infinite length; an infinite sequence $\{1=d_{0},d_{1},\dotsc,
d_{n},\dotsc\}$ in $D^{1}$ is called a {\em divided power
sequence} if $\Delta(d_{n})=\sum_{i=0}^{n}d_{i}\otimes d_{n-i}$
(see \cite[p.\ 268]{Sweedler1969}). This assumption is necessarily satisfied if
$\ch R = 0$ (for each primitive $\partial\in P(D)$, $\{1,\partial,
\partial^{2}/2,\dotsc,\partial^{n}/n!,\dotsc\}$ is a divided power sequence of
infinite length).
If $\ch R = p > 0$, this is equivalent to the Verschiebung map
$D^{1}\rightarrow R^{1/p}\otimes D^{1}$ being surjective;
see \cite{Heyneman-Sweedler}. The assumption is also equivalent to saying
that $D$ is smooth as a cocommutative coalgebra.

Moreover this implies that, for a commutative algebra $A$,
the $A$-algebra $\Hom_{R}(D^{1},A)$ with the convolution product is the
projective limit of power series $A$-algebras
(see \cite[p.\ 278]{Sweedler1969}). Thus, if $A$ is a domain (resp.\ reduced),
then $\Hom_{R}(D^{1},A)$ is also a domain (resp.\ reduced). Furthermore,
$\Hom_{R}(D,A)$ is isomorphic to the direct product of
$A$-algebras isomorphic to $\Hom_{R}(D^{1},A)$ indexed by $G$:
\[ \begin{array}{rcccl}
\Hom_{R}(D,A) & \xrightarrow{\sim} & \Hom_{R}(RG,\Hom_{R}(D^{1},A))
 & \xrightarrow{\sim} & {\displaystyle\prod_{g\in G}}\Hom_{R}(D^{1},A)
\vspace{1ex} \\
\varphi & \mapsto & [g\mapsto [d\mapsto\varphi(dg)]] & \mapsto &
([d\mapsto\varphi(dg)])_{g}.
\end{array} \]
Hence, if $A$ is reduced, then $\Hom_{R}(D,A)$ is also reduced.
(These facts implies that $D^{1}$ and $D$ are convolutionally reduced
in the sense of \cite[Definition 5.2]{Tyc-Wisniewski}.)

As in \cite[p.\ 505]{Takeuchi1989}, we have the following:
\begin{lemma}
\label{lem-stable-ideals}
Let $A$ be a $D$-module algebra and $\rho_{A} : A\rightarrow\Hom_{R}(D,A)$
the algebra map associated with the structure on $A$.

{\rm (i)} If $J\subset A$ is an ideal, then $\rho_{A}^{-1}(\Hom_{R}(D,J))$ is
a $D$-stable ideal, which is maximal among $D$-stable ideals included in $J$.
Therefore $J$ is a $D$-stable ideal iff $J=\rho_{A}^{-1}(\Hom_{R}(D,J))$.

{\rm (ii)} If $I\subset A$ is a $D$-stable ideal, then also the radical
$\sqrt{I}$ is a $D$-stable ideal.

{\rm (iii)} If $P\subset A$ is a prime ideal, then
$(\rho_{A}^{1})^{-1}(\Hom_{R}(D^{1},P))$ is a prime $D^{1}$-stable ideal.
Here $\rho_{A}^{1} : A\rightarrow\Hom_{R}(D^{1},A)$ is the algebra map
associated with the $D^{1}$-module algebra structure on $A$.
\end{lemma}
\begin{proof}
(i) This is easily seen.

(ii) Since the algebra $\Hom_{R}(D,A/\sqrt{I})\simeq
\Hom_{R}(D,A)/\Hom_{R}(D,\sqrt{I})$ is reduced, we have
$\Hom_{R}(D,\sqrt{I})$ is a radical ideal of $\Hom_{R}(D,A)$.
Hence its pull-back $\rho_{A}^{-1}(\Hom_{R}(D,\sqrt{I}))$ is also a radical
ideal. By part (i), it includes $I$.
On the other hand, $\rho_{A}^{-1}(\Hom_{R}(D,\sqrt{I}))$ is included 
in $\sqrt{I}$. Therefore $\rho_{A}^{-1}(\Hom_{R}(D,\sqrt{I}))=\sqrt{I}$.

(iii) Since $\Hom_{R}(D^{1},A/P)\simeq\Hom_{R}(D^{1},A)/\Hom_{R}(D^{1},P)$
is a domain, $\Hom_{R}(D^{1},P)$ is a prime ideal. Thus its pull-back
$(\rho_{A}^{1})^{-1}(\Hom_{R}(D^{1},P))$ is also prime.
\end{proof}

Let $K$ be a $D$-module algebra and $\Omega(K)$ the set of all minimal prime
ideals in $K$. Then $G$ acts on $\Omega(K)$. Let $G_{\Omega(K)}$ denote
the normal subgroup consisting of those elements in $G$ which stabilize every
$P\in\Omega(K)$.

\begin{proposition}
\label{prop-NS}
Suppose that $K$ is noetherian as a ring and simple as a $D$-module algebra.
Then $\Omega(K)$ is a finite set.

{\rm (i)} The action of $G$ on $\Omega(K)$ is transitive, so that the subgroups
$G_{P}$ of stabilizers of $P\in\Omega(K)$ are conjugate to each other.

{\rm (ii)} Every $P\in\Omega(K)$ is $D^{1}$-stable, so that $K/P$ is a
$D(G_{P})$-module domain. This is simple as a $D(G_{\Omega(K)})$-module
algebra.

{\rm (iii)} Let $P\in\Omega(K)$, and set $K_{1}=K/P$. Then we have a natural
isomorphism of $D$-module algebras,
\[ K\simeq\Phi_{G_{P}}(K_{1}). \]
\end{proposition}
\begin{proof}
(ii) Let $\rho : K\rightarrow\Hom_{R}(D^{1},K)$ be the algebra map
associated with the $D^{1}$-module algebra structure on $K$. Put
$P'=\rho^{-1}(\Hom_{R}(D^{1},P))$. By Lemma \ref{lem-stable-ideals} (iii),
$P'$ is a $D^{1}$-stable prime ideal included in $P$.
Then we have $P=P'$ by the minimality of $P$. Hence $P$
is $D^{1}$-stable. (This also follows from
\cite[Theorem 5.9 (2)]{Tyc-Wisniewski}.)

Let $P\subset J\subsetneq K$ be a $D(G_{\Omega(K)})$-stable ideal.
Then, $\bigcap_{g\in G/G_{\Omega(K)}}gJ$ is
$D$-stable, and hence is zero. Since $P$ is prime, there exists $g$ such that
$gJ\subset P$, and so $P\subset J\subset g^{-1}P$. By the minimality of
$g^{-1}P$, we have $P=J$ ($=g^{-1}P$).
Thus $K/P$ is a simple $D(G_{\Omega(K)})$-module algebra.

(i) Let $P\in\Omega(K)$. We see
\begin{equation}
\label{eq-intersection-of-prime-ideals}
\bigcap_{g\in G}gP = \bigcap_{Q\in\Omega(K)}Q = 0,
\end{equation}
since the intersections are both $D$-stable. The first equality implies
$\{gP\;|\;g\in G\}=\Omega(K)$; this proves (i).

(iii) By (i), $g\mapsto gP$ gives a bijection $G/G_{P}\xrightarrow{\sim}
\Omega(K)$. If $Q$ and $Q'$ in $\Omega(K)$ are distinct, then ($Q\subsetneq$)
$Q+Q'=K$, by (ii). This together with (\ref{eq-intersection-of-prime-ideals})
proves that the natural map gives an isomorphism,
\[ K\xrightarrow{\sim}\prod_{Q\in\Omega(K)}K/Q =\prod_{g\in G/G_{P}}K/gP. \]
Obviously, $\Phi_{G_{P}}(K_{1})$ is isomorphic to the last direct product.
\end{proof}

For a commutative ring $K$ in general, we say that $K$ is {\em total} iff 
every non-zero divisor in $K$ is invertible.

\begin{corollary}
\label{cor-total}
Let $K$ be a noetherian simple $D$-module algebra as above.
Then the following are equivalent.
\begin{enumerate}
\renewcommand{\labelenumi}{\rm(\alph{enumi})}
\item $K$ is total;
\item $K$ is artinian as a ring;
\item The Krull dimension $\Kdim(K)=0$, or in other words $\Omega(K)$ equals
the set of all maximal ideals in $K$.
\end{enumerate}

If these conditions are satisfied, every $K\# D$-module is free as a
$K$-module.
\end{corollary}
\begin{proof}
Each condition is equivalent to that for any/some $P\in\Omega(K)$, $K/P$ is
a field. The last assertion holds true by part (iii) of the last proposition
and by Proposition \ref{prop-copying-functor}.
\end{proof}

\begin{definition}
\label{def-AS}
A $D$-module algebra $K$ is said to be {\em AS} iff it is artinian and simple.
By the corollary above, this is equivalent to that $K$ is total, noetherian
and simple.
\end{definition}

For later use we prove some results. The following lemma is a particular
case of \cite[Theorem 3.4]{Tyc-Wisniewski}.

\begin{lemma}
\label{lem-localization}
Let $A$ be a $D$-module algebra, and let $T\subset A$ be a $G$-stable
multiplicative subset which contains no zero-divisors.
The $D$-module algebra structure on $A$ can be uniquely extended to the
localization $T^{-1}A$ of $A$ by $T$.
($D^{1}$ may not be of Birkhoff-Witt type.)
\end{lemma}
\begin{proof}
Let $\rho : A\rightarrow\Hom_{R}(D,A)\subset\Hom_{R}(D,T^{-1}A)$ be
the algebra map associated with the $D$-module algebra structure on $A$.
For each $t\in T$, we see $\rho(t)(g)=g(t)$ are invertible in $T^{-1}A$
for all $g\in G$.
Hence $\rho(t)$ ($t\in T$) is invertible in $\Hom_{R}(D,T^{-1}A)$ by
\cite[Corollary 9.2.4]{Sweedler1969}. This implies that $\rho$ 
is uniquely extended to an algebra map
$\tilde{\rho} : T^{-1}A\rightarrow \Hom(D,T^{-1}A)$ so that
$\tilde{\rho}(1/t)\ast\rho(t)=\vep$ ($t\in T$); cf.\ the proof of
\cite[Proposition 1.9]{Takeuchi1989}.
We have thus obtained the measuring action
\[ d(a/t)=\tilde{\rho}(a/t)(d)\qquad (d\in D,\quad a\in A,\quad t\in T) \]
by $D$ on $T^{-1}A$. It remains to prove that this makes $T^{-1}A$ a
$D$-module. We have only to see that
\[ cd(1/t)=c(d(1/t))\qquad (c,d\in D,\quad t\in T). \]
This holds, since the two maps $D\otimes D\rightarrow T^{-1}A$, given by
$c\otimes d\mapsto cd(1/t)$ and $c\otimes d\mapsto c(d(1/t))$ coincide,
being the convolution-inverse of $c\otimes d\mapsto cdt$.

For convenience, we describe how to extend the action of $D$ explicitly.
Let $t\in T$. The action of $D$ on $1/t$ is given by:
\begin{eqnarray*}
g(1/t) & = & 1/g(t)\qquad (g\in G), \\
d(1/t) & = & \frac{\varepsilon(d)}{t}-\frac{dt-\varepsilon(d)t}{t^{2}}
{}+\frac{1}{t^{3}}\sum_{(d)}(d_{(1)}t-\varepsilon(d_{(1)})t)
(d_{(2)}t-\varepsilon(d_{(2)})t) \\
 & & {}-\frac{1}{t^{4}}\sum_{(d)}(d_{(1)}t-\varepsilon(d_{(1)})t)
(d_{(2)}t-\varepsilon(d_{(2)})t)(d_{(3)}t-\varepsilon(d_{(3)})t)+\dotsb\qquad
(d\in D^{1}).
\end{eqnarray*}
We observe that the right hand side of $d(1/t)$ ($d\in D^{1}$)
in the equation is a finite sum by the coradical filtration; see the proof of
\cite[Lemma 9.2.3]{Sweedler1969} or \cite[Lemma 5.2.10]{Montgomery}.
\end{proof}

\begin{lemma}
\label{AS-subalgebra}
Let $L$ be an AS $D$-module algebra, and let $K\subset L$ be a $D$-module
subalgebra. If $K$ is total, then $K$ is AS.
\end{lemma}
\begin{proof}
Given an element $x\neq 0$ in $L=\prod_{P\in\Omega(L)}L/P$, define the support
of $x$ by
\begin{equation}
\label{eq-support}
\supp(x)=\{P\in\Omega(L)\;|\;x\not\in P\}.
\end{equation}
One sees that $x$ is a non-zero divisor iff $\supp(x)=\Omega(L)$.

Choose an element $x\neq 0$ in $K$ with minimal support. Then for $g\in G$,
the supports $\supp(x)$ and $\supp(gx)$ are either equal or disjoint,
according to $x(gx)$ being non-zero or zero. By Proposition \ref{prop-NS} (i),
we have those elements $x,g_{1}x,\dotsc,g_{r}x$ in $K$ with disjoint supports,
whose sum is a non-zero divisor. Let $y$ be the inverse of the sum; this is
indeed in $K$, since $K$ is total. We see that $e:=xy$ is a (primitive)
idempotent in $K$ with $\supp(e)=\supp(x)$. By the minimality of the support,
each non-zero element in $eK$ has $\supp(x)$ as its support, and hence has an
inverse in $eK$, just as $x$ above. We have $K=\prod_{i=1}^{r}g_{i}eK$, the
direct product of the fields $g_{i}eK$; this proves the lemma.
\end{proof}

\begin{corollary}
\label{total-quotient-ring}
Let $A$ be a $D$-module subalgebra in an AS $D$-module algebra $L$.

{\rm (i)} Every non-zero divisor $x$ in $A$ has full support:
$\supp(x)=\Omega(L)$ (see (\ref{eq-support})).

{\rm (ii)} Let $K=Q(A)$ denote the total quotient ring of $A$; this is realized
in $L$ by {\rm (i)}. Then $K$ is an AS $D$-module subalgebra of $L$.
\end{corollary}
\begin{proof}
Let $T$ be the set of all non-zero divisors in $A$. Then, $K=T^{-1}A$.

(i) Choose an $x\in T$ such that $\supp(x)$ is minimal in
$\{\supp(t)\;|\;t\in T\}$. If $\supp(x)\neq\Omega(L)$, then there is a
$g\in G$ such that $\supp(gx)\cap\supp(x)=\emptyset$, which implies
$x(gx)=0$, a contradiction.

(ii) Let $\rho_{L}: L\rightarrow\Hom(D,L)$ be the algebra map associated to
the $D$-module algebra structure on $L$. It restricts to $\rho : A\rightarrow
\Hom(D,A)$ associated to $A$. If $t\in T$, $\rho_{L}(1/t)$ is the inverse
of $\rho(t)$ in $\Hom(D,L)$, and hence is contained in $\Hom(D,T^{-1}A)$
by the proof of Lemma \ref{lem-localization}.
This implies that $K$ ($=T^{-1}A$) is a $D$-module subalgebra of $L$.
$K$ is AS by Lemma \ref{AS-subalgebra}.
\end{proof}

\section{Sweedler's correspondence theorem}
\label{sec-Sweedler}

Let $K\subset A$ be an inclusion of $D$-module algebras.
Then $A\otimes_{K}A$ has a coalgebra structure in the symmetric tensor
category $({}_{A}({}_{D}\cM)_{A},\otimes_{A},A)$ given by
\[ \begin{array}{l}
\Delta : A\otimes_{K}A \rightarrow (A\otimes_{K}A)\otimes_{A}(A\otimes_{K}A),
\quad a\otimes b \mapsto (a\otimes 1)\otimes (1\otimes b), \vspace{1ex} \\
\vep : A\otimes_{K}A\rightarrow A,\quad a\otimes b\mapsto ab.
\end{array} \]
In particular $A\otimes_{K}A$ is an {\em $A$-coring} (a coalgebra
in $({}_{A}\cM_{A},\otimes_{A},A)$, the category of $(A,A)$-bimodules);
see \cite{Sweedler1975}.
The next theorem is an analog on AS $D$-module algebras of Sweedler's
correspondence theorem \cite[Theorem 2.1]{Sweedler1975},
which play a key role to obtain the Galois correspondence later.

\begin{theorem}
\label{thm-Sweedler}
Let $K\subset L$ be an inclusion of AS $D$-module algebras.
Let $\cC_{L/K}$ be the set of all $D$-stable coideals of $L\otimes_{K}L$ and
$\cA_{L/K}$ the set of all intermediate AS $D$-module algebras of $L/K$.

{\rm (i)} For $M\in\cA_{L/K}$, we have $J_{M}:=\Ker(L\otimes_{K}L
\twoheadrightarrow L\otimes_{M}L)\in\cC_{L/K}$.

{\rm (ii)} For $J\in\cC_{L/K}$, let $\pi : L\otimes_{K}L\twoheadrightarrow
L\otimes_{K}L/J$ be the canonical surjection.
Then $M_{J}:=\{a\in L\;|\;a\pi(1\otimes 1)
=\pi(1\otimes 1)a\}\in\cA_{L/K}$.

{\rm (iii)} The correspondence 
\[ \begin{array}{rcl}
\cC_{L/K} & \leftrightarrow & \cA_{L/K} \vspace{1ex} \\
J & \rightarrow & M_{J} \vspace{1ex} \\
J_{M} & \leftarrow & M
\end{array} \]
is bijective.
\end{theorem}
\begin{proof}
(i) Both $L\otimes_{K}L$ and $L\otimes_{M}L$ are coalgebras in
$({}_{L}({}_{D}\cM)_{L},\otimes_{L},L)$ and obviously
$L\otimes_{K}L\twoheadrightarrow L\otimes_{M}L$ is a $D$-linear $L$-coring map.
Hence its kernel $J_{M}$ is a $D$-stable coideal of $L\otimes_{K}L$.

(ii) We easily see that $M_{J}$ is a subalgebra of $L$ which contains $K$.
For any $d\in D$ and $a\in M_{J}$, we have $d(a)\pi(1\otimes 1)
=d(a\pi(1\otimes 1))=d(\pi(1\otimes 1)a)=\pi(1\otimes 1)d(a)$.
Thus $M_{J}$ is a $D$-module subalgebra of $L$.
Let $t$ be a non-zero divisor in $M_{J}$.
By Corollary \ref{total-quotient-ring} (i), $t$ is invertible in $L$. We see
$t^{-1}\pi(1\otimes 1)=t^{-1}\pi(1\otimes 1)tt^{-1}=t^{-1}t\pi(1\otimes 1)
t^{-1}=\pi(1\otimes 1)t^{-1}$ and hence $t^{-1}\in M_{J}$. This implies that
$M_{J}$ is total. Therefore $M_{J}$ is an intermediate AS $D$-module algebra
of $L/K$ by Lemma \ref{AS-subalgebra}.

(iii) Take an $M\in\cA_{L/K}$.
For all $a\in M$, we have $a\otimes 1 - 1\otimes a\in J_{M}$. Hence $M\subset
M_{J_{M}}$. By the definition of $M_{J_{M}}$, one sees
$M_{J_{M}}\otimes_{M}M_{J_{M}}\simeq M_{J_{M}}\otimes_{M}M$.
Since $M_{J_{M}}$ is an $M\# D$-module, it is a free $M$-module
(see Corollary \ref{cor-total}). Hence $M_{J_{M}}=M$.

Conversely, take $J\in\cC_{L/K}$.
Let $\xi : L\otimes_{M_{J}}L\rightarrow L\otimes_{K}L/J$, $a\otimes b\mapsto
a\pi(1\otimes 1)b$, which is a surjective $D$-linear $L$-coring map.
Then we have $J_{M_{J}}\subset J$ by chasing the following commutative diagram:
\[ \begin{CD}
0 @>>> J_{M_{J}} @>>> L\otimes_{K}L @>>> L\otimes_{M_{J}}L @>>> 0
\qquad \mbox{(exact)} \\
@. @. @| @VV{\xi}V @. \\
0 @>>> J @>>> L\otimes_{K}L @>>> L\otimes_{K}L/J @>>> 0 \qquad \mbox{(exact).}
\end{CD} \]
If we prove that $\xi$ is injective, then $J=J_{M_{J}}$ follows.

For a fixed $P\in\Omega(M_{J})$, put $M'=M_{J}/P$ ($=\Psi_{G_{P}}(M_{J})$),
$L'=L/PL$ ($=\Psi_{G_{P}}(L)$), and $C=\Psi_{G_{P}}(L\otimes_{K}L/J)
=(L\otimes_{K}L/J)e_{1}=e_{1}(L\otimes_{K}L/J)e_{1}$
(where $e_{1}\in M_{J}$ is the primitive idempotent such that $M'=M_{J}e_{1}$).
Then $C$ is a coalgebra in $({}_{L'}({}_{D(G_{P})}\cM)_{L'},\otimes_{L'},L')$
by Corollary \ref{prop-copying-functor2}. It suffices to prove
that $\xi'=\Psi_{G_{P}}(\xi) : L'\otimes_{M'}L'\rightarrow C$ is injective.

Regarding $C$ merely as an $L'$-coring, let $\gA$ be the category of right
$C$-comodules in $({}_{L'}\cM_{L'},\otimes_{L'},L')$. Then $\gA$ is an
abelian category since $C$ is a left free $L'$-module. $L'$ has a natural
$C$-comodule structure given by
\[ \lambda : L'\rightarrow L'\otimes_{L'}C\simeq C,\quad
a\mapsto \pi(1\otimes 1)a\quad (=e_{1}\pi(1\otimes 1)e_{1}a). \]
We see $\End_{\gA}(L')\xrightarrow{\sim} M'$, $f\mapsto f(1)$ is an algebra
isomorphism. ($f\in\End_{L'}(L')$ is a $C$-comodule map iff
$f(1)\pi(1\otimes 1)=\pi(1\otimes 1)f(1)$.)
On the other hand, $\Hom_{\gA}(L',C)\xrightarrow{\sim} L'$,
$f\mapsto (\vep\circ f)(1)$ is an $M'$-module isomorphism whose inverse
is given by $a\mapsto [b\mapsto a\pi(1\otimes 1)b]$. Indeed,
for $f\in \Hom_{\gA}(L',C)$,
\begin{eqnarray*}
 & & (\vep\circ f)(1)\pi(1\otimes 1)b = (\vep\circ f)(1)\lambda(b)
= (\vep\circ f)(1){\displaystyle\sum_{(b)}}b_{(0)}\otimes b_{(1)} \\
 & = & {\displaystyle\sum_{(b)}}(\vep\circ f)(b_{0})\otimes b_{(1)} 
= ((\vep\otimes\id)\circ(f\otimes\id))(\lambda(b))
= (\vep\otimes\id)(\Delta_{C}(f(b))) = f(b).
\end{eqnarray*}
We will show that $L'$ is simple in $\gA$, concluding that $\xi'$
is injective by Proposition \ref{prop-simple-objects}.
Every simple subobject of $L'$ is of the form $eL'$, where $e$ is an
idempotent of $L'$. Since $\lambda$ is $D(G_{P})$-linear, we see that
$g(eL')$ is also a simple object for each $g\in G$. Each $g(eL')$ coincides
or trivially intersects with $eL'$ since $g(eL')\cap eL'$ is also a right
$C$-comodule. It follows from Proposition \ref{prop-NS} (i) that $L'$
is semisimple in $\gA$. But the endomorphism ring $\End_{\gA}(L')\simeq M'$
is a field. This implies that $L'$ is a simple object in $\gA$.
\end{proof}

\begin{example}
\label{Fibonacci}
Let $R=\bbq$ and $D=\bbq[\tau,\tau^{-1}]$, the ring of linear difference 
operators. Take $\alpha=\{(\frac{1+\sqrt{5}}{2})^{n}\}\in\cS_{\bbc}$ as in
Introduction. Then $\bbq=\bbq(1,1)\subset L=\bbq(\sqrt{5},\alpha)\times
\bbq(\sqrt{5},\alpha)$ (where $\tau(1,0)=(0,1)$, $\tau(0,1)=(1,0)$,
$L^{D}=\bbq(\sqrt{5})=\bbq(\sqrt{5})(1,1)$) is an inclusion of AS
$D$-module algebras.
Write $e_{1}=(1,0), e_{2}=(0,1)\in L$.
Then $\cA_{L/\bbq}$ and $\cC_{L/\bbq}$ correspond as the following:
\begin{center}
\begin{tabular}{|c|c|}
\hline
$\cA_{L/\bbq}$ & $\cC_{L/\bbq}$ \\ \hline
$L$ & $\mathrm{Span}_{L,L}\{\sqrt{5}\otimes 1-1\otimes\sqrt{5},
\ e_{1}\otimes e_{2},\ e_{2}\otimes e_{1},
\ \alpha\otimes 1-1\otimes\alpha\}$  \\ \hline
$\bbq(\sqrt{5},\alpha)$ & $\mathrm{Span}_{L,L}\{\sqrt{5}\otimes 1
{}-1\otimes\sqrt{5},\ \alpha\otimes 1-1\otimes\alpha\}$  \\ \hline
$\bbq(\sqrt{5})\times\bbq(\sqrt{5})$ & $\mathrm{Span}_{L,L}\{
e_{1}\otimes e_{2},\ e_{2}\otimes e_{1},\ \sqrt{5}\otimes 1-1\otimes\sqrt{5}\}$ 
\\ \hline
$\bbq(\alpha)$ & $\mathrm{Span}_{L,L}\{\alpha\otimes 1-1\otimes\alpha\}$
\\ \hline
$\bbq(\sqrt{5})$ & $\mathrm{Span}_{L,L}\{\sqrt{5}\otimes 1-1\otimes
\sqrt{5}\}$ \\ \hline
$\bbq\times\bbq$ & $\mathrm{Span}_{L,L}\{e_{1}\otimes e_{2},
\ e_{2}\otimes e_{1}\}$ \\ \hline
$\bbq$ & $0$ \\ \hline
\end{tabular}
\end{center}
We will see that $L/\bbq$ is not a Picard-Vessiot extension but
$L/\bbq(\sqrt{5})$ is.
\end{example}

\section{Galois correspondence for Picard-Vessiot extensions}
\label{sec-PV-extensions}

Let $K\subset A$ be an inclusion of $D$-module algebras. Then $A\otimes_{K}A$
has an algebra structure naturally and become a $D$-module algebra since
$D$ is cocommutative. Thus $(A\otimes_{K}A)^{D}$ is a $K^{D}$-subalgebra
of $A\otimes_{K}A$.

\begin{definition}
\label{def-PV-extension}
Let $K\subset L$ be an inclusion of AS $D$-module algebras.
We say that $L/K$ is a {\em Picard-Vessiot}, or {\em PV}, {\em extension}
if the following conditions are satisfied:
\begin{enumerate}
\renewcommand{\labelenumi}{\rm(\alph{enumi})}
\item $K^{D}=L^{D}$; this will be denoted by $k$.
\item There exists a $D$-module subalgebra $A\subset L$ including $K$, such
that the total quotient ring $Q(A)$ of $A$ equals $L$, and the $k$-subalgebra
$H:=(A\otimes_{K}A)^{D}$ generates the left (or equivalently right)
$A$-module $A\otimes_{K}A$: $A\cdot H=A\otimes_{K}A$
(or $H\cdot A=A\otimes_{K}A$).
\end{enumerate}
\end{definition}

\begin{proposition}
\label{prop-PV-Hopf-algebra}
Let $L/K$ be a PV extension of AS $D$-module algebras and take $A,H$ as
in the condition (b) above.

{\rm (i)} The product map $\mu : A\otimes_{k}H\rightarrow A\otimes_{K}A$,
$\mu(a\otimes h)=a\cdot h$ is a $D$-linear isomorphism.

{\rm (ii)} The $A$-coring structure maps $\Delta,\vep$ of $A\otimes_{K}A$
induce $k$-algebra maps $\Delta_{H} : H\rightarrow H\otimes_{k}H$,
$\vep_{H} : H\rightarrow k$. Then $(H,\Delta_{H},\vep_{H})$ becomes a
commutative Hopf algebra over $k$. The antipode is induced from the twist
map $\tw : A\otimes_{K}A\rightarrow A\otimes_{K}A$, $a\otimes b\mapsto
b\otimes a$.

{\rm (iii)} The $k$-algebra map $\theta : A\rightarrow A\otimes_{k}H$,
$\theta(a)=\mu^{-1}(1\otimes a)$ makes $A$ a right $H$-comodule.
$A/K$ is necessarily a right $H$-Galois extension \cite[Sect.\ 8.1]{Montgomery}
in the sense that
\[ {}_{A}\theta : A\otimes_{K}A\rightarrow A\otimes_{k}H,\quad
{}_{A}\theta(a\otimes b)=a\theta(b) \]
is an isomorphism. A Hopf algebra structure on $H$ with this property
is unique.

{\rm (iv)} Such an algebra $A$ that satisfies the condition (b) above is
unique.
\end{proposition}
\begin{proof}
(i) By Corollary \ref{cor-simple-module-algebras}, the natural map
$L\otimes_{k}(L\otimes_{K}A)^{D}\rightarrow L\otimes_{K}A$ is injective.
Since the map $\mu$ is its restriction, it is injective. On the other hand,
$\mu$ is surjective by the condition (b). Note that this can be uniquely
extended to an isomorphism $L\otimes_{k}H\xrightarrow{\sim}L\otimes_{K}A$.

(ii) Since $A^{D}=k$ by the condition (a), $\vep$ maps $H$ into $k$.
The twofolds of $\mu$:
\[ \varphi : A\otimes_{k}H\otimes_{k}H\xrightarrow{\mu\otimes\id}
A\otimes_{K}A\otimes_{k}H\xrightarrow{\id\otimes\mu}
A\otimes_{K}A\otimes_{K}A \]
is a $D$-linear isomorphism. This induces an algebra isomorphism
$\varphi_{1}:=\id\otimes\mu|_{H\otimes_{k}H} : H\otimes_{k}H\xrightarrow{\sim}
(A\otimes_{K}A\otimes_{K}A)^{D}$. Similarly the threefolds
of $\mu$ induces an algebra isomorphism
$\varphi_{2} : H\otimes_{k}H\otimes_{k}H
\xrightarrow{\sim}(A\otimes_{K}A\otimes_{K}A\otimes_{K}A)^{D}$. Since
$\Delta : A\otimes_{K}A\rightarrow (A\otimes_{K}A)\otimes_{A}(A\otimes_{K}A)
\simeq A\otimes_{K}A\otimes_{K}A$ maps $H$ into
$(A\otimes_{K}A\otimes_{K}A)^{D}$, a $k$-algebra map
$\Delta_{H} : H\rightarrow H\otimes_{k}H$ is induced by $\Delta_{H}
=\varphi_{1}^{-1}\circ\Delta|_{H}$. We see
\[ \begin{array}{l}
\varphi_{2}\circ (\Delta_{H}\otimes\id_{H})\circ\Delta_{H}
=(\Delta\otimes\id)\circ\varphi_{1}\circ\Delta_{H}
=(\Delta\otimes\id_{A})\circ\Delta|_{H},\quad \\
\varphi_{2}\circ (\id_{H}\otimes\Delta_{H})\circ\Delta_{H}
=(\id\otimes\Delta)\circ\varphi_{1}\circ\Delta_{H}
=(\id_{A}\otimes\Delta)\circ\Delta|_{H}.
\end{array} \]
Then we have $(\Delta_{H}\otimes\id)\circ\Delta_{H}
=(\id\otimes\Delta_{H})\circ\Delta_{H}$ by the coassociativity of
$\Delta$. The counitary property is easily seen.
Therefore $(H,\Delta_{H},\vep_{H})$ is a commutative bialgebra.
Since $D$ is cocommutative, $\tw$ maps $H$ into $H$. Put $S=\tw|_{H}$.
For $w=\sum_{i}(a_{i}\otimes a_{i}')\otimes(b_{i}'\otimes b_{i})\in
H\otimes_{k}H$, we see
\[ \varphi_{1}(w)=\sum_{i}a_{i}\otimes a_{i}'b_{i}'\otimes b_{i},\quad
m((\id\otimes S)(w))=\sum_{i}a_{i}b_{i}\otimes a_{i}'b_{i}' \]
where $m$ denotes the multiplication of $H$. Thus, for
$h=\sum_{i}a_{i}\otimes b_{i}\in H$,
\[ m((\id\otimes S)(\Delta_{H}(h)))
=m((\id\otimes S)(\varphi_{1}^{-1}(\sum_{i}a_{i}\otimes 1 \otimes
b_{i})))=\sum_{i}a_{i}b_{i}\otimes 1. \]
This implies $\id\ast S=m\circ (\id\otimes S)\circ\Delta_{H}
=u\varepsilon$ where $u : k\rightarrow H$ is the unit map of $H$.
We have $S\ast\id
=u\varepsilon$ similarly. Therefore $S$ is the antipode of $H$.

(iii) We see
\[ (\varphi\circ(\theta\otimes\id)\circ\theta)(a)
=1\otimes 1\otimes a
=(\Delta\circ\mu\circ\theta)(a)
=(\varphi\circ(\id\otimes\Delta_{H})\circ\theta)(a) \]
for all $a\in A$. Thus we have $(\theta\otimes\id)\circ\theta
=(\id\otimes\Delta_{H})\circ\theta$. On the other hand,
\[ (\id\otimes\varepsilon_{H})\circ\theta =\varepsilon\circ\mu
\circ\theta =\id. \]
Therefore $(A,\theta)$ is an $H$-comodule. 
The map ${}_{A}\theta$, being $\mu^{-1}$, is an isomorphism.
Since this interprets $\theta$ into the natural right
$A\otimes_{K}A$-comodule structure $A\rightarrow A\otimes_{A}(A\otimes_{K}A)
\simeq A\otimes_{K}A$, $a\mapsto 1\otimes a$ on $A$, we see the described
uniqueness of the structure on $H$.

(iv) This follows in the same way as \cite[Lemma 2.5]{Takeuchi1989}. We
include the proof for convenience.
If $A,B$ satisfy the condition (b), then also $AB$ satisfy it.
Thus we may assume $A\subset B$. Put $H_{A}=(A\otimes_{K}A)^{D}$,
$H_{B}=(B\otimes_{K}B)^{D}$, the corresponding Hopf algebras. Then $H_{A}$ is
a Hopf subalgebra of $H_{B}$.
Hence $H_{B}/H_{A}$ is a faithfully flat extension
(see \cite[Theorem 3.1]{Takeuchi1972} or \cite[Ch.\ 14]{Waterhouse}).
The extension $(L\otimes_{K}B)/(L\otimes_{K}A)$ is identified
with $(L\otimes_{k}H_{B})/(L\otimes_{k}H_{A})$ through the $\mu$-isomorphism.
It follows that $B/A$ is a faithfully flat
extension since $L$ is a free $K$-module. Hence $aA=aB\cap A$
for all $a\in A$ since the canonical map $A/aA\rightarrow B\otimes_{A}(A/aA)
\simeq B/aB$ is injective. For any $b\in B$, there exists a non-zero divisor
$a\in A$ such that $ab\in A$.
Since $ab\in aB\cap A=aA$ and since $a$ is a non-zero divisor, $b\in A$
follows. Therefore we have $A=B$.
\end{proof}

\begin{definition}
\label{def-principal-algebra}
$A$ (resp., $H$) is called {\em the principal algebra}
(resp., {\em the Hopf algebra}) {\em for} $L/K$. To indicate these we say
that $(L/K,A,H)$ is a PV extension. The associated affine group
scheme $\bG(L/K):=\Spec H$ is called {\em the PV group scheme for $L/K$}.
\end{definition}

\begin{theorem}
\label{thm-Galois-correspondence1}
Let $L/K$ be a PV extension of AS $D$-module algebras with
the Hopf algebra $H$. Let $\cA_{L/K}$ be the set of intermediate
AS $D$-module algebras of $L/K$ and $\cH\cI_{H}$ the set of all Hopf ideals of
$H$. Then $\cA_{L/K}$ and $\cH\cI_{H}$ correspond bijectively as follows:
\[ \begin{array}{l}
\cA_{L/K}\rightarrow\cH\cI_{H},\quad
M \mapsto H\cap\Ker(L\otimes_{K}L\rightarrow L\otimes_{M}L), \vspace{1.5ex} \\
\cH\cI_{H}\rightarrow\cA_{L/K},\quad I\mapsto
\{x\in L\;|\;x\otimes 1 - 1\otimes x\in I\cdot (L\otimes_{K}L)\}.
\end{array} \]
\end{theorem}

This theorem is obtained as the composite of 1-1 correspondences given by
Theorem \ref{thm-Sweedler} and the next proposition.
For a commutative algebra $A$ (resp.\ a $D$-module algebra $B$),
let $\cI(A)$ (resp.\ $\cI_{D}(B)$) denote the set of all ideals of $A$
(resp.\ $D$-stable ideals of $B$).

\begin{proposition}
\label{prop-ideal-correspondence}
Let $(L/K,A,H)$ be a PV extension.

{\rm (i)} $\cI(H)$ and $\cI_{D}(L\otimes_{K}L)$ correspond bijectively
as follows:
\[ \begin{array}{l}
\cI(H)\rightarrow\cI_{D}(L\otimes_{K}L),\quad
I \mapsto I\cdot (L\otimes_{K}L), \vspace{1.5ex} \\
\cI_{D}(L\otimes_{K}L)\rightarrow\cI(H),\quad J\mapsto J\cap H.
\end{array} \]

{\rm (ii)} Under the correspondence, $J$ is a $D$-stable coideal iff $I$
is a Hopf ideal.
\end{proposition}
\begin{proof} This follows in the same way as
\cite[Proposition 2.6]{Takeuchi1989}. We include the proof for convenience.

(i) Since $L$ is the total quotient ring of $A$, we have
$\cI_{D}(L\otimes_{K}L)\subset\cI_{D}(A\otimes_{K}A)$.
Furthermore, $\cI_{D}(L\otimes_{K}A)\cap\cI_{D}(A\otimes_{K}L)
=\cI_{D}(L\otimes_{K}L)$ in $\cI_{D}(A\otimes_{K}A)$.
Considering the $\mu$-isomorphism, we claim the map
\[ \cI(H)\rightarrow\cI_{D}(A\otimes_{k}H)\xrightarrow{\sim}
\cI_{D}(A\otimes_{K}A),\quad
I\mapsto A\otimes_{k}I\mapsto I\cdot(A\otimes_{K}A) \]
is injective with the image $\cI_{D}(L\otimes_{k}H)\simeq
\cI_{D}(L\otimes_{K}A)$. The injectivity is clear.
Since $A\otimes_{k}I = (L\otimes_{k}I)\cap
(A\otimes_{k}H)$, the image is contained in $\cI_{D}(L\otimes_{k}H)$.
Then it suffices to prove that every $D$-stable ideal of $L\otimes_{k}H$
is written as $L\otimes_{k}I$ by some $I\in\cI(H)$.
Let $\ga\subset L\otimes_{k}H$ be a $D$-stable ideal and take the canonical
map $\varphi : H\xrightarrow{\sim}(L\otimes_{k}H)^{D}\twoheadrightarrow
((L\otimes_{k}H)/\ga)^{D}$. Put $I=\Ker\varphi = \ga\cap H$, an ideal
of $H$. Since $L\otimes_{k}((L\otimes_{k}H)/\ga)^{D}\rightarrow
(L\otimes_{k}H)/\ga$ is injective (Corollary \ref{cor-simple-module-algebras}),
we have $\ga = L\otimes_{k}I$ by chasing the following diagram:
\[ \begin{CD}
0 @>>> L\otimes_{k}I @>>> L\otimes_{k}H @>{\mathrm{id}\otimes\varphi}>>
L\otimes_{k}((L\otimes_{k}H)/\ga)^{D} @>>> 0 \\
 @. @. @| @VVV @. \\
0 @>>> \ga @>>> L\otimes_{k}H @>>> (L\otimes_{k}H)/\ga @>>> 0.
\end{CD} \]
Then the claim is proved.
By symmetry, we see the image of $\cI(H)\rightarrow
\cI_{D}(A\otimes_{K}A)$ is also equal to $\cI_{D}(A\otimes_{K}L)$.
It follows $\cI_{D}(L\otimes_{K}A)=\cI_{D}(A\otimes_{K}L)
=\cI_{D}(L\otimes_{K}L)$ in $\cI_{D}(A\otimes_{K}A)$, proving (i).

(ii) By the similar discussion to (i), we have that $\cI(H\otimes_{k}H)$
and $\cI_{D}(L\otimes_{K}L\otimes_{K}L)$ correspond bijectively. If
$I\leftrightarrow J$ in (i), then $I\otimes_{k}H\leftrightarrow
J\otimes_{K}L$ and $H\otimes_{k}I\leftrightarrow L\otimes_{K}J$.
Therefore, $\Delta_{H}(I)\subset I\otimes_{k}H+H\otimes_{k}I$ iff
$\Delta J\subset J\otimes_{K}L+L\otimes_{K}J$. On the other hand,
$\Ker((L\otimes_{K}L)\otimes_{L}(L\otimes_{K}L)\rightarrow (L\otimes_{K}L/J)
\otimes_{L}(L\otimes_{K}L/J))=J\otimes_{L}(L\otimes_{K}L)
{}+(L\otimes_{K}L)\otimes_{L}J=J\otimes_{K}L+L\otimes_{K}J$
holds since $J$, $L\otimes_{K}L$, and $J\otimes_{K}L/J$ are free $L$-modules.
It follows that $I$ is a biideal of $H$ iff $J$ is a $D$-stable
coideal of $L\otimes_{K}L$. It is known that every biideal of a commutative
Hopf algebra over a field is a Hopf ideal (see \cite[Theorem 1 (iv)]{Nichols}).
\end{proof}

Actually, $\cI_{D}(A\otimes_{K}A)=\cI_{D}(L\otimes_{K}L)$ holds in the above
(replace $L$ with $A$ in the proof of (i) and use Proposition
\ref{prop-principal-algebra} (i)).

\begin{example}
In Example \ref{Fibonacci}, if we put $K=\bbq(\sqrt{5})$, then $L/K$ is a
PV extension. The principal algebra and the Hopf algebra are given by
$A=K[\alpha,\alpha^{-1}]\times K[\alpha,\alpha^{-1}]$ and
$H=K[g_{1},g_{2}]$ with grouplikes
$g_{1}=\alpha\otimes\alpha^{-1}$, $g_{2}=(e_{1}-e_{2})\otimes(e_{1}-e_{2})$.
In this case, $\cA_{L/K}$, $\cC_{L/K}$, $\cH\cI_{H}$ correspond as follows:
\begin{center}
\begin{tabular}{|c|c|c|}
\hline
$\cA_{L/K}$ & $\cC_{L/K}$ & $\cH\cI_{H}$ \\ \hline
$L$ & $\mathrm{Span}_{L,L}\{e_{1}\otimes e_{2},\ e_{2}\otimes e_{1},
\ \alpha\otimes 1-1\otimes\alpha\}$
 & $H^{+}=\lan g_{1}-1,g_{2}-1\ran$ \\ \hline
$K(\alpha)$ & $\mathrm{Span}_{L,L}\{\alpha\otimes 1-1\otimes\alpha\}$ &
$\lan g_{1}-1\ran$ \\ \hline
$K\times K$ & $\mathrm{Span}_{L,L}\{e_{1}\otimes e_{2},\ e_{2}\otimes e_{1}\}$
 & $\lan g_{2}-1\ran$ \\ \hline
$K$ & 0 & 0 \\ \hline
\end{tabular}
\end{center}
\end{example}

\vspace{2ex}

\begin{proposition}
\label{prop-Galois-correspondence2}
Let $(L/K,A,H)$ be a PV extension. Suppose $\cA_{L/K}\ni M\leftrightarrow
I\in\cH\cI_{H}$ in Theorem \ref{thm-Galois-correspondence1}.

{\rm (i)} $(L/M,AM,H/I)$ is a PV extension.

{\rm (ii)} $A^{\coinv H/I}=\{a\in A\;|\;\theta(a)-a\otimes 1\in
A\otimes_{k}I\}=A\cap M$, and the $\mu$-isomorphism $A\otimes_{k}H
\xrightarrow{\sim}A\otimes_{K}A$ induces an isomorphism
$A\otimes_{k}H^{\coinv H/I}\xrightarrow{\sim}A\otimes_{K}(A\cap M)$.
\end{proposition}
\begin{proof}
(i) We have an isomorphism $L\otimes_{k}H/I\xrightarrow{\bar{\mu}}
L\otimes_{M}AM$ by considering the following diagram:
\[ \begin{CD}
0 @>>> L\otimes_{k}I @>>> L\otimes_{k}H @>>> L\otimes_{k}H/I @>>> 0 \\
@. @VV{\mu}V @VV{\mu}V @VV{\bar{\mu}}V @. \\
0 @>>> I\cdot(L\otimes_{K}A) @>>> L\otimes_{K}A @>>> L\otimes_{M}AM @>>> 0.
\end{CD} \]
Restrict the diagram as
\[ \begin{CD}
0 @>>> AM\otimes_{k}I @>>> AM\otimes_{k}H @>>> AM\otimes_{k}H/I @>>> 0 \\
@. @VV{\mu}V @VV{\mu}V @VV{\bar{\mu}}V @. \\
0 @>>> I\cdot(AM\otimes_{K}A) @>>> AM\otimes_{K}A @>>> AM\otimes_{M}AM @>>>
0. \end{CD} \]
Since $AM\otimes_{K}A=(AM\otimes_{K}K)\cdot H$, we have
$AM\otimes_{M}AM=AM\cdot(AM\otimes_{M}AM)^{D}$ through the surjection
$AM\otimes_{K}A\rightarrow AM\otimes_{M}AM$. On the other hand,
$\bar{\mu}$ induces an isomorphism $H/I\xrightarrow{\sim}
(AM\otimes_{M}AM)^{D}$ of Hopf algebras.

(ii) This follows by considering the next diagram:
\begin{center}
\begin{picture}(300,110)
\put(0,100){$0$}
\put(10,103){\vector(1,0){30}}
\put(45,100){$A\otimes_{K}(A\cap M)$}
\put(123,103){\vector(1,0){30}}
\put(158,100){$A\otimes_{K}A$}
\put(200,106){\vector(1,0){30}}
\put(200,100){\vector(1,0){30}}
\put(235,100){$A\otimes_{K}AM\otimes_{M}AM$}
\put(270,62){\vector(0,1){30}}
\put(173,62){\line(0,1){30}}
\put(177,62){\line(0,1){30}}
\multiput(80,62)(0,4){7}{\line(0,1){2}}
\put(80,92){\vector(0,1){0}}
\put(0,50){$0$}
\put(10,53){\vector(1,0){30}}
\put(45,50){$A\otimes_{K}A^{\coinv H/I}$}
\put(123,53){\vector(1,0){30}}
\put(158,50){$A\otimes_{K}A$}
\put(200,56){\vector(1,0){30}}
\put(200,50){\vector(1,0){30}}
\put(235,50){$A\otimes_{K}A\otimes_{k}H/I$}
\put(270,12){\vector(0,1){30}}
\put(175,12){\vector(0,1){30}}
\multiput(80,12)(0,4){7}{\line(0,1){2}}
\put(80,45){\vector(0,1){0}}
\put(0,0){$0$}
\put(10,3){\vector(1,0){30}}
\put(50,0){$A\otimes_{k}H^{\coinv H/I}$}
\put(123,3){\vector(1,0){30}}
\put(158,0){$A\otimes_{k}H$}
\put(200,6){\vector(1,0){30}}
\put(200,0){\vector(1,0){30}}
\put(235,0){$A\otimes_{k}H\otimes_{k}H/I$.}
\end{picture}
\end{center}
\end{proof}

Let $H$ be a commutative Hopf algebra over $k$. It is known that normal
Hopf ideals $I$ of $H$ and Hopf subalgebras $H_{1}\subset H$ correspond
bijectively by $H_{1}=H^{\coinv H/I}={}^{\coinv H/I}H$ and $I=HH_{1}^{+}$
(see \cite{Takeuchi1972}). Let $(V,\lambda)$ be a right $H$-comodule
in general. If we put $V_{1}=\lambda^{-1}(V\otimes_{k}H_{1})$
($=V^{\coinv H/I}$), then we have
$\lambda(V_{1})\subset V_{1}\otimes_{k}H_{1}$. Indeed, write $\lambda(v)
=\sum_{i}v_{i}\otimes h_{i}\in V\otimes_{k}H_{1}$ for $v\in V_{1}$, where
$h_{i}$ are $k$-linearly independent. Then
\[ \sum_{i}\sum_{(v_{i})}(v_{i})_{(0)}\otimes(v_{i})_{(1)}
\otimes h_{i}
=\sum_{i}\sum_{(h_{i})}v_{i}\otimes(h_{i})_{(1)}
\otimes (h_{i})_{(2)}\in V\otimes_{k}H_{1}\otimes_{k}H_{1}. \]
This implies $v_{i}\in V_{1}$.
As in \cite[Theorem 2.9]{Takeuchi1989}, we have the following proposition.

\begin{proposition}
\label{prop-Galois-correspondence3}
Let $(L/K,A,H)$ be a PV extension and $H_{1}\subset H$ a Hopf subalgebra.
Put $I=HH_{1}^{+}$ and $A_{1}=\theta^{-1}(A\otimes_{k}H_{1})=A^{\coinv H/I}$.
Let $L_{1}$ be the total quotient ring of $A_{1}$ so that $L_{1}$
is an intermediate AS $D$-module algebra of $L/K$.

{\rm (i)} $(L_{1}/K,A_{1},H_{1})$ is a PV extension.

{\rm (ii)} $I$ is the Hopf ideal of $H$ which corresponds to $L_{1}$.

{\rm (iii)} $H_{1}\mapsto L_{1}$ gives a 1-1 correspondence between the
Hopf subalgebras of $H$ and the intermediate AS $D$-module algebras which
are PV extensions over $K$.
\end{proposition}
\begin{proof}
(i) Since $\theta(A_{1})\subset A_{1}\otimes_{k}H_{1}$, we have
$\mu(A_{1}\otimes_{k}H_{1})\supset A_{1}\otimes_{K}A_{1}$.
Consider $A\otimes_{K}A$ as a right $H$-comodule by the structure map
$\id\otimes\theta$. Then the inclusion $H\hookrightarrow A\otimes_{K}A$
is an $H$-comodule map; recall that, for $h=\sum_{i}a_{i}\otimes_{K}b_{i}\in H
=(A\otimes_{K}A)^{D}$,
\[ \Delta(h)=\sum_{i}a_{i}\otimes_{K}\mu^{-1}(1\otimes_{K}b_{i})
=\sum_{i}a_{i}\otimes_{K}\theta(b_{i}). \]
Thus $H_{1}\subset A\otimes_{K}\theta^{-1}(A\otimes_{k}H_{1})
=A\otimes_{K}A_{1}$. Also we have $H_{1}=S(H_{1})\subset\tw(A\otimes_{K}A_{1})
=A_{1}\otimes_{K}A$. Hence $H_{1}\subset A_{1}\otimes_{K}A_{1}$ and so
$\mu(A_{1}\otimes_{k}H_{1})\subset A_{1}\otimes_{K}A_{1}$.
This implies that $\mu : A\otimes_{k}H\xrightarrow{\sim}A\otimes_{K}A$
induces a $D$-linear isomorphism $A_{1}\otimes_{k}H_{1}\xrightarrow{\sim}
A_{1}\otimes_{K}A_{1}$. Therefore $H_{1}=(A_{1}\otimes_{K}A_{1})^{D}$ and
$A_{1}\otimes_{K}A_{1}=A_{1}\cdot H_{1}$.

(ii) Let $M$ be the intermediate AS $D$-module algebra of $L/K$ which
corresponds to $I$. Then $M\supset A\cap M=A^{\coinv H/I}=A_{1}$ by
Proposition \ref{prop-Galois-correspondence2} (ii). Since $L_{1}$ is the
smallest AS $D$-module subalgebra of $L$ which includes $A_{1}$, we have
$L_{1}\subset M$. Let $I'$ be the Hopf ideal which corresponds to $L_{1}$.
Then $I=HH_{1}^{+}\subset H\cap\Ker(L\otimes_{K}L\twoheadrightarrow
L\otimes_{L_{1}}L)=I'$ since the counit $\vep : H_{1}\rightarrow k$
is a restriction of $\mathrm{mult} : A_{1}\otimes_{K}A_{1}\rightarrow A_{1}$.
Thus we have $L_{1}\supset M$.

(iii) Let $L_{1}$ be an intermediate AS $D$-module algebra of $L/K$ such that
$(L_{1}/K,A_{1},H_{1})$ is a PV extension. Since $A\otimes_{K}A=A\cdot H$
and $A_{1}\otimes_{K}A_{1}=A_{1}\cdot H_{1}$, we have
$A_{1}A\otimes_{K}A_{1}A=A_{1}A\cdot H_{1}H$. This implies that $A_{1}A$
is the principal algebra for $L/K$ and hence $A_{1}A=A$ by Proposition
\ref{prop-PV-Hopf-algebra} (iv). Thus $A_{1}\subset A$ and $H_{1}\subset H$,
a Hopf subalgebra. Since the $\mu$-isomorphism $A_{1}\otimes_{k}H_{1}
\xrightarrow{\sim}A_{1}\otimes_{K}A_{1}$ induces a left $A$-module
isomorphism $A\otimes_{k}H_{1}\xrightarrow{\sim}A\otimes_{K}A_{1}$,
we have $A_{1}=\theta^{-1}(A\otimes_{k}H_{1})$. This proves (iii).
\end{proof}

Finally in this section, we prove two important properties on principal
algebras which are used later.
\begin{proposition}
\label{prop-principal-algebra}
Let $(L/K,A,H)$ be a PV extension.

{\rm (i)} $A$ is simple as a $D$-module algebra.

{\rm (ii)} $A$ contains all primitive idempotents in $L$.
\end{proposition}
\begin{proof}
(i) The following proof is essentially the same as that of
\cite[Theorem 2.11]{Takeuchi1989}.

Let $0\neq I\subset A$ be a $D$-stable ideal. Since $L\otimes_{K}I
\in\cI_{D}(L\otimes_{K}A)$, there exists an ideal
$\ga\in\cI(H)$ such that $L\otimes_{K}I=\ga\cdot (L\otimes_{K}A)$
by the proof of Proposition \ref{prop-ideal-correspondence} (i).
But $IL=L$ since $L$ is simple and hence $L\otimes_{K}IL=L\otimes_{K}L$.
This implies that the $D$-stable ideal of $L\otimes_{K}L$ which corresponds
to $\ga$ is $L\otimes_{K}L$. Thus $\ga = H$. Therefore $L\otimes_{K}I
=H\cdot (L\otimes_{K}A)=L\otimes_{K}A$, concluding $I=A$.

(ii) Since $L$ is a localization of $A$, we have $\Omega(L)\subset
\Omega(A)$ via $P\mapsto P\cap A$.
We see $A\subset\prod_{P\in\Omega(L)}A/P\cap A$. It remains to prove
that if $P\neq Q$ in $\Omega(L)$, then the sum $J:=P\cap A + Q\cap A$ equals
$A$. If $J\subsetneq A$ on the contrary, one sees as in the proof of
Proposition \ref{prop-NS} (ii) that $J=P\cap A=Q\cap A$ by Part (i),
and so $P=Q$.
\end{proof}

\begin{remark}
Before the version four of this article, the following assertion
was written in this place (as a lemma):

{\it Let $K$ be a finite product of fields and $A$ a commutative $K$-algebra.
Let $T$ be a multiplicative subset of $A$ which has no zero divisors and put
$L=T^{-1}A$. Then every separable $K$-subalgebra of $L$ is included in $A$.
} 

But this is false (the sentence ``the proof can be reduced to the case that
$A$ is a domain" in the proof was wrong). The following is a counter example.
Let $A=\bbc[X,Y]/\lan XY\ran$ and $L$ the total quotient ring of $A$.
Then $A$ has no nontrivial idempotent. But $L$ has a nontrivial idempotent
$X/(X+Y)$. 
Excuse me for making such a mistake and I hope the reader of old versions
was free from the wrong argument.
\end{remark}

\section{Translation into affine group schemes}
\label{sec-translation}

For an inclusion of $D$-module algebras $K\subset A$, let $\AutDlinKalg(A)$
denote the group of $D$-linear $K$-algebra automorphisms of $A$.
Let $\bAutDlinKalg(A)$ denote the associated group functor over $k=K^{D}$;
it associates to each commutative $k$-algebra $T$ the automorphism group
$\AutDlinKTalg(A\otimes_{k}T)$, where $T$ is considered as a $D$-module
algebra by the trivial action $dt=\vep(d)t$ ($d\in D$, $t\in T$). 
As in \cite[Appendix]{Takeuchi1989}, we have the following:
\begin{theorem}
Let $(L/K,A,H)$ be a PV extension and $\bG(L/K)=\Spec H$ the PV group
scheme. Then the linear representation $\phi : \bG(L/K)\rightarrow
\bGL(A)$ arising from the $H$-comodule structure $\theta : A\rightarrow
A\otimes_{k}H$ gives an isomorphism $\bG(L/K)\xrightarrow{\sim}
\bAutDlinKalg(A)$ of affine $k$-group schemes. In particular,
$\bG(L/K)(k)\simeq\AutDlinKalg(A)=\AutDlinKalg(L)$.
\end{theorem}
\begin{proof}
Let $T$ be a commutative $k$-algebra. For $\alpha\in\bG(L/K)(T)=\Alg_{k}(H,T)$,
$\phi_{T}(\alpha)$ is given by
\[ \phi_{T}(\alpha) : A\otimes_{k}T\xrightarrow{\sim}A\otimes_{k}T,\quad
a\otimes t\mapsto \sum_{(a)}a_{(0)}\otimes\alpha(a_{(1)})t. \]
We easily see $\phi_{T}(\alpha)\in\AutDlinKTalg(A\otimes_{k}T)$.
We will construct the inverse $\psi : \bAutDlinKalg(A)\rightarrow
\bG(L/K)$. For an element $\beta\in\AutDlinKTalg(A\otimes_{k}T)$,
let ${}_{A}\beta$ denote the left $A$-linear extension of $\beta|_{A} :
A\rightarrow A\otimes_{k}T$. Namely, ${}_{A}\beta : A\otimes_{K}A\rightarrow
A\otimes_{k}T$, $a\otimes b\mapsto a\beta(b\otimes 1)$.
Consider the $D$-linear $A$-algebra map
\[ A\otimes_{k}H\xrightarrow{\mu}A\otimes_{K}A\xrightarrow{{}_{A}\beta}
A\otimes_{k}T. \]
We see this maps the constants $H$ into $T$. Then we have a $k$-algebra map
$\psi_{T}(\beta):=({}_{A}\beta\circ\mu)^{D}={}_{A}\beta|_{H}\in\Alg_{k}(H,T)
=\bG(L/K)(T)$ so that $(\id_{A}\otimes\psi_{T}(\beta))\circ\theta
=\beta|_{A}$. This gives a homomorphism $\psi : \bAutDlinKalg(A)\rightarrow
\bG(L/K)$. Indeed, for $\beta,\gamma\in\AutDlinKTalg(A\otimes_{k}T)$,
\[ (\gamma\circ\beta)(a\otimes 1)=\sum_{(a)}\gamma(a_{(0)}\otimes 1)
\psi_{T}(\beta)(a_{(1)})=\sum_{(a)}a_{(0)}\otimes \psi_{T}(\gamma)(a_{(1)})
\psi_{T}(\beta)(a_{(2)})\quad (a\in A). \]
One easily sees $\phi_{T}\circ\psi_{T}=\id$. For $\alpha\in\bG(L/K)(T)$, we see
${}_{A}\phi_{T}(\alpha)=(\vep\otimes\alpha)\circ(\id\otimes\theta)$ where
$\vep : A\otimes_{K}A\rightarrow A$, the counit, and hence
${}_{A}\phi_{T}(\alpha)|_{H}=(\vep_{H}\otimes\alpha)\circ\Delta_{H}=\alpha$.
This implies $\psi_{T}\circ\phi_{T}=\id$.
\end{proof}

Theorem \ref{thm-Galois-correspondence1} and Proposition
\ref{prop-Galois-correspondence2}, \ref{prop-Galois-correspondence3} can be
translated as follows.
\begin{theorem}
Let $L/K$ be a PV extension of AS $D$-module algebras.

{\rm (i)} If $M$ is an intermediate AS $D$-module algebra of $L/K$, then
$L/M$ is also a PV extension and the PV group scheme $\bG(L/M)$ is identified
with a closed subgroup scheme of $\bG(L/K)$. Then intermediate AS $D$-module
algebras of $L/K$ and closed subgroup schemes of $\bG(L/K)$ correspond
bijectively by $M\mapsto\bG(L/M)$.

{\rm (ii)} Under the correspondence above, $M/K$ is a PV extension iff
$\bG(L/K)\vtr\bG(L/M)$. If this is the case, we have an
isomorphism $\bG(M/K)\simeq \bG(L/K)/\bG(L/M)$ of affine group schemes.
\end{theorem}

\section{Copying and interlacing}
\label{sec-copying}

In this section, we investigate how PV extensions change through the
functor $\Phi$ and $\Psi$ described in Section \ref{sec-relative-Hopf-modules}.

First we easily see the following:
\begin{lemma}
\label{lem-interlacing1}
Let $G_{1}\subset G$ be a subgroup of finite index. Write $\Phi =\Phi_{G_{1}}$.
Let $K_{1}\subset L_{1}$ be an inclusion of AS $D(G_{1})$-module algebras.
Then $(L_{1}/K_{1},A_{1},H)$ is a PV extension iff
$(\Phi(L_{1})/\Phi(K_{1}),\Phi(A_{1}),H)$ is a PV extension of AS $D$-module
algebras.
\end{lemma}
\begin{proof}
This follows by Proposition \ref{prop-copying-functor2} and
Lemma \ref{lem-preservation}.
\end{proof}

\begin{remark}
\label{reduction-to-field}
Let $K\subset L$ be an inclusion of AS $D$-module algebras. Choose $\gp\in
\Omega(K)$, and let $P_{1},\dotsc,P_{r}$ be all those elements in $\Omega(L)$
that lie over $\gp$. Define $K_{1}=K/\gp$, $L_{1}=L/\gp L
=\prod_{i=1}^{r}L/P_{i}$.
Then we have an inclusion $K_{1}\subset L_{1}$ of AS $D(G_{\gp})$-module
algebras such that the induced inclusion $\Phi_{G_{\gp}}(K_{1})\subset
\Phi_{G_{\gp}}(L_{1})$ is identified with $K\subset L$. We can thus reduce
to the case where $K$ is a field, especially to discuss PV extensions
by Lemma \ref{lem-interlacing1}.
\end{remark}

\begin{example}
\label{normal-subgroup}
Let $G_{1}\subset G$ be a {\em normal} subgroup of finite index. Let $K$ be a
$D$-module field. Regarding this as a $D(G_{1})$-module algebra, define
$L=\Phi_{G_{1}}(K)$. We then have the inclusion
\[ K\hookrightarrow L=\bigoplus_{g\in G/G_{1}}g\otimes K,\quad x\mapsto 
\sum_{g}g\otimes g^{-1}x \]
of AS $D$-module algebras. If $K^{D(G_{1})}=K^{D}$, then $K^{D}=L^{D}$
($=:k$) by Lemma \ref{lem-preservation} (i). Moreover, $(L/K,L,H)$ is a PV
extension with $H=k(G/G_{1})^{\ast}$, the dual of the group algebra
$k(G/G_{1})$. In fact, we see that the elements
\[ e_{g} := \sum_{h\in G/G_{1}}(h\otimes 1)\otimes_{K}(hg\otimes 1)
\quad (g\in G/G_{1}) \]
in $L\otimes_{K}L$ are $D$-invariant, and behave as the dual basis in $H$ of
the group elements $g$ ($\in G/G_{1}$) in $k(G/G_{1})$. Namely,
$\Delta(e_{g})=\sum_{h}e_{gh^{-1}}\otimes e_{h}$, $\varepsilon(e_{g})
=\delta_{1,g}$, $S(e_{g})=e_{g^{-1}}$.
The $H$-comodule structure $\theta : L\rightarrow L\otimes_{k}H$ is given by
\[ \theta(h\otimes x)=\sum_{g}(hg^{-1}\otimes gx)\otimes_{k}e_{g}, \]
as is seen from following computation in $L\otimes_{K}L$:
\[ 1\otimes_{K}(h\otimes x)=\sum_{f}(f\otimes f^{-1}hx)\otimes_{K}(h\otimes 1)
=\sum_{g}(hg^{-1}\otimes gx)\otimes_{K}(h\otimes 1)
=\sum_{g}(hg^{-1}\otimes gx)\cdot e_{g}. \]
\end{example}

\begin{proposition}
\label{prop-interlacing2}
Let $(L/K,A,H)$ be a PV extension of AS $D$-module algebras. Choose arbitrarily
$P\in\Omega(L)$, and write $\Phi = \Phi_{G_{P}}$. Let $\gp = P\cap K$ ($\in
\Omega(K)$). Define
\[ K_{1}=K/\gp,\quad A_{1}=A/P\cap A,\quad L_{1}=L/P. \]
Then,

{\rm (i)} $A\simeq\Phi(A_{1})$.

{\rm (ii)} $\Phi(K_{1})$ is identified with the $K$-subalgebra $\hat{K}$ of
$L$ which is spanned over $K$ by the primitive idempotents in $L$.

{\rm (iii)} $(L_{1}/K_{1},A_{1},\bar{H}=H/I)$ is a PV extension of
$D(G_{P})$-module fields, where $I=H\cap\Ker(L\otimes_{K}L\rightarrow
L\otimes_{\hat{K}}L)$; cf.\ \cite[Corollary 1.16]{vanderPut-Singer1997}.

{\rm (iv)} The subalgebra of $H$
\[ B = \{h\in H\;|\;\Delta(h)\equiv h\otimes 1\mod H\otimes_{k}I\}
\quad (=H^{\coinv\bar{H}}) \]
is a separable $k$-algebra. We have a right $\bar{H}$-colinear $B$-algebra
isomorphism $H\simeq B\otimes_{k}\bar{H}$.

{\rm (v)} If $G_{P}$ is normal in $G_{\gp}$, then $B\subset H$ is a Hopf
subalgebra which is isomorphic to $k(G_{\gp}/G_{P})^{\ast}$, and we have
an extension
\[ k(G_{\gp}/G_{P})^{\ast}\rightarrowtail H\twoheadrightarrow \bar{H} \]
of Hopf algebras; cf.\ \cite[Corollary 1.17]{vanderPut-Singer1997}.
\end{proposition}
\begin{proof}
(i) This follows from Proposition \ref{prop-principal-algebra}.

(ii) This is easy to see.

(iii) By Proposition \ref{prop-Galois-correspondence2} (i),
we have a PV extension
\[ (L/\hat{K},A,\bar{H})=(\Phi(L_{1})/\Phi(K_{1}),\Phi(A_{1}),\bar{H}). \]
part (iii) now follows by Lemma \ref{lem-interlacing1}.

For the remaining (iv), (v), we may suppose that $K$ is a field and
so $\gp = 0$, $G_{\gp}=G$ by Remark \ref{reduction-to-field}.

(iv) By Proposition \ref{prop-Galois-correspondence2} (ii), we see
$A\otimes_{k}B\simeq A\otimes_{K}\hat{K}$.
Since this is a separable $A$-algebra, $B$ is a separable $k$-algebra by
\cite[Ch.\ II, Proposition 1.8]{DeMeyer-Ingraham}.

Recall that $A$ has the natural, right $\bar{H}$-comodule $k$-algebra
structure $A\xrightarrow{1\otimes -}A\otimes_{\hat{K}}A\simeq
A\otimes_{k}\bar{H}$; in fact, $A$ is also a left $\bar{H}$-comodule
$k$-algebra. We see that the map
\begin{equation}
\label{D-lin-barH-colin-alg-splitting}
\sigma : \Phi(A_{1}\otimes_{K}A_{1})=A\otimes_{\hat{K}}A\rightarrow
A\otimes_{K}A
\end{equation}
given by $g\otimes (a\otimes_{K}b)\mapsto (g\otimes a)\otimes_{K}(g\otimes b)$
($g\in G/G_{P}$) is a $D$-linear, two-sided $\bar{H}$-colinear $k$-algebra
splitting of $A\otimes_{K}A\rightarrow A\otimes_{\hat{K}}A$. The induced
$\sigma^{D} : \bar{H}\rightarrow H$ is a two-sided $\bar{H}$-colinear
$k$-algebra splitting of $H\rightarrow\bar{H}$. It follows by
\cite[Theorem 7.2.2]{Montgomery} that
\begin{equation}
\label{Doi-Takeuchi}
B\otimes_{k}\bar{H}\rightarrow H,\quad b\otimes x\mapsto b\sigma^{D}(x)
\end{equation}
gives a right $\bar{H}$-colinear $B$-algebra isomorphism.

(v) If $G_{P}$ is normal in $G$, then $(\hat{K}/K,\hat{K},k(G/G_{P})^{\ast})$
is a PV extension by Example \ref{normal-subgroup}. Hence the assertion
follows from Proposition \ref{prop-Galois-correspondence3}.
\end{proof}

\begin{theorem}
\label{thm-copying}
Let $K\subset L$ be an inclusion of AS $D$-module algebras. Choose arbitrarily
$P\in\Omega(L)$, and let $\gp = P\cap K$ ($\in\Omega(K)$). Then $L/K$ is a 
PV extension if
\begin{enumerate}
\renewcommand{\labelenumi}{\rm(\alph{enumi})}
\item $G_{P}$ is normal in $G_{\gp}$, and
\item the inclusion $K_{1}:=K/\gp\subset L_{1}:=L/P$ of $D(G_{P})$-module
fields is a PV extension.
\end{enumerate}
The converse holds true if the field $K^{D}$ ($=L^{D}$) of $D$-invariants is
separably closed.
\end{theorem}
\begin{proof}
This follows by slightly modifying the last proof, as follows. We may suppose
that $K$ is a field.

Suppose that $(L_{1}/K_{1},A_{1},\bar{H})$ is a PV extension. Define
$A=\Phi(A_{1})$ with $\Phi =\Phi_{G_{P}}$. Recall from Proposition
\ref{prop-interlacing2} that if $L/K$ is PV, the principal algebra
must be $A$. As was seen in the last proof, $A\otimes_{K}A$ is a right
$\bar{H}$-comodule $k$-algebra and the map $\sigma$ given in
(\ref{D-lin-barH-colin-alg-splitting}) induces an $\bar{H}$-colinear
$k$-algebra map $\sigma^{D} : \bar{H}\rightarrow (A\otimes_{K}A)^{D}$.
Again by \cite[Theorem 7.2.2]{Montgomery}, we have a $D$-linear and $\bar{H}$-colinear
isomorphism
\[ A\otimes_{K}\Phi(K)\otimes_{k}\bar{H}\simeq A\otimes_{K}A \]
of algebras over $A\otimes_{K}\Phi(K)$; see (\ref{Doi-Takeuchi}).
It follows that $L/K$ is a PV extension iff the natural injection
\begin{equation}
\label{natural-injection}
A\otimes_{k}(A\otimes_{K}\Phi(K))^{D}\rightarrow A\otimes_{K}\Phi(K)
\end{equation}
is surjective. If $G_{P}$ is normal in $G$, then this is surjective since
by Example \ref{normal-subgroup}, $A\otimes_{k}(\Phi(K)\otimes_{K}\Phi(K))^{D}
\rightarrow A\otimes_{K}\Phi(K)$ is already surjective.

To prove the converse, we may suppose (b) by Proposition
\ref{prop-interlacing2} (iii), and that the map given in
(\ref{natural-injection}) is an isomorphism by the argument above.
It follows that
\[ \dim_{k}(A\otimes_{K}\Phi(K))^{D}=[G:G_{P}]. \]
If $k$ is a separably closed field, then $B=(A\otimes_{K}\Phi(K))^{D}$ is
isomorphic to $k\times\dotsb\times k$, the product of $[G:G_{P}]$ copies of
$k$. The isomorphism (\ref{natural-injection}) induces a $D(G_{P})$-linear
isomorphism
\[ L_{1}\otimes_{k}B\xrightarrow{\sim} L_{1}\otimes_{K}\Phi(K)
\xrightarrow{\sim}\Phi(L_{1})=L. \]
Thus all primitive idempotents in $L$ are fixed by the action of $G_{P}$
and hence $G_{P}=G_{\Omega(L)}$. This implies that $G_{P}$ is normal in $G$.
\end{proof}

As will be seen from the following, the second half
does not necessarily hold true unless $k$ is separably closed.

\begin{example}
Let $G=D_{3}=\{1,\sigma,\sigma^{2},\tau,\sigma\tau,\sigma^{2}\tau\}$
be the dihedral group of order $6$
($\sigma^{3}=1$, $\tau^{2}=1$, $\sigma\tau=\tau\sigma^{2}$),
$G_{1}=\{1,\tau\}\subset G$, and $D=\bbq G$. Put $k=\bbq$, $K=\bbq(\sqrt{-1})$,
and identify $G_{1}$ with $\Gal(K/k)$.
Then $G$ acts on $K$ so that $\sigma a = a$ for all $a\in K$.
Take the system of representatives $\{1,\sigma,\sigma^{2}\}$
for $G/G_{1}$ and let $L=\Phi_{G_{1}}(K)=1\otimes K+\sigma\otimes K
{}+\sigma^{2}\otimes K$. Then $L/K$ is a PV extension of AS $D$-module
algebras with the Hopf algebra
\[ H=\bbq[z_{1},z_{2},z_{3}]/\lan z_{1}+z_{2}-1,z_{2}^{2}-z_{2},
z_{3}^{2}+z_{2},z_{2}z_{3}-z_{3}\ran \]
where $\vep(z_{i})=\delta_{1i}$ ($i=1,2,3$) and
\begin{eqnarray*}
\Delta(z_{1}) & = & z_{1}\otimes z_{1}+\frac{1}{2}z_{2}\otimes z_{2}
{}+\frac{1}{2}z_{3}\otimes z_{3}, \\
\Delta(z_{2}) & = & z_{1}\otimes z_{2}+z_{2}\otimes(z_{1}+\frac{1}{2}z_{2})
{}-\frac{1}{2}z_{3}\otimes z_{3}, \\
\Delta(z_{3}) & = & z_{1}\otimes z_{3}-\frac{1}{2}z_{2}\otimes z_{3}
{}+z_{3}\otimes(z_{1}-\frac{1}{2}z_{2}).
\end{eqnarray*}
Indeed, a Hopf algebra isomorphism $H\xrightarrow{\sim}(L\otimes_{K}L)^{D}$
is given by
\begin{eqnarray*}
z_{1} & \mapsto  & (1\otimes 1)\otimes_{K}(1\otimes 1)+(\sigma\otimes 1)
\otimes_{K}(\sigma\otimes 1)+(\sigma^{2}\otimes 1)\otimes_{K}(\sigma^{2}
\otimes 1),   \\
z_{2} & \mapsto & (1\otimes 1)\otimes_{K}(\sigma\otimes 1+\sigma^{2}\otimes 1)
{}+(\sigma\otimes 1)\otimes_{K}(1\otimes 1+\sigma^{2}\otimes 1) \\
 & & {}+(\sigma^{2}\otimes 1)\otimes_{K}(1\otimes 1+\sigma\otimes 1), \\
z_{3} & \mapsto & -\sqrt{-1}(1\otimes 1)\otimes_{K}(\sigma\otimes 1-\sigma^{2}
\otimes 1)+\sqrt{-1}(\sigma\otimes 1)\otimes_{K}(1\otimes 1-\sigma^{2}
\otimes 1) \\
 & & {}-\sqrt{-1}(\sigma^{2}\otimes 1)\otimes_{K}(1\otimes 1-\sigma\otimes 1).
\end{eqnarray*}
The PV group scheme $\bG(L/K)=\Spec H$ is a twisted form of $\bbz/3\bbz$
(see \cite[(6.4)]{Waterhouse}):
\begin{eqnarray*}
\bbq(\sqrt{-1})\otimes_{\bbq}H & \xrightarrow{\sim} & \bbq(\sqrt{-1})\times
\bbq(\sqrt{-1})\times\bbq(\sqrt{-1})\simeq \bbq(\sqrt{-1})(\bbz/3\bbz)^{\ast}
\\
z_{1} & \mapsto & (1,0,0) \\
\frac{1}{2}(z_{2}+\sqrt{-1}z_{3}) & \mapsto & (0,1,0) \\
\frac{1}{2}(z_{2}-\sqrt{-1}z_{3}) & \mapsto & (0,0,1).
\end{eqnarray*}
\end{example}

\section{Splitting algebras}
\label{sec-splitting}

Let $K$ be an AS $D$-module algebra and $V$ a $K\# D$-module.
The rank $\rk_{K}(V)$ of the free $K$-module $V$
will be called the $K$-rank; see Corollary \ref{cor-total}.

\begin{definition}
\label{splittig-algebra}
Let $K\subset L$ be an inclusion of AS $D$-module algebras and $V$ a
$K\# D$-module. We say that $V$ {\em splits in} $L/K$, or $L/K$ is a
{\em splitting algebra for $V$} iff there is an $L\# D$-linear injection
$L\otimes_{K}V\hookrightarrow L^{\Lambda}$ into some power $L^{\Lambda}$ of
$L$. $K\langle V\rangle$ denotes the smallest AS $D$-module subalgebra of $L$
including $K$ and $f(V)$ for all $f\in\Hom_{K\# D}(V,L)$. If $L=K\langle V
\rangle$ and $V$ splits in $L/K$, we say $L/K$ is a {\em minimal splitting
algebra for $V$}.
\end{definition}

Similarly to \cite[Proposition 3.1]{Takeuchi1989}, we have the following:

\begin{proposition}
\label{prop-splitting-algebras1}
Let $K\subset L$ be an inclusion of AS $D$-module algebras and $V$
a $K\#D$-module.

{\rm (i)} If $V$ splits in $L/K$, every $K\#D$-submodule of $V$
splits in $L/K$.

{\rm (ii)} If $V$ splits in $L/K$, it splits in $K\langle V\rangle /K$.

{\rm (iii)} $V$ splits in $L/K$ iff the canonical $L$-module map
\begin{equation}
\label{eq-dense}
L\otimes_{L^{D}}\Hom_{K\# D}(V,L)\rightarrow \Hom_{K}(V,L),
\end{equation}
has a dense image; in other words, the map
\begin{equation}
\label{eq-inj}
L\otimes_{K}V\rightarrow\Hom_{L^{D}}(\Hom_{K\# D}(V,L),L),\quad
a\otimes v\mapsto [f\mapsto af(v)]
\end{equation}
is injective.
\end{proposition}
\begin{proof}
(i) Since all $K\# D$-modules are free $K$-modules, this follows
immediately by the definition.

(ii) If $V$ splits in $L/K$, then the image of $L\otimes_{K}V\rightarrow
L^{\Lambda}$ is in $K\langle V\rangle^{\Lambda}$ by the definition. Thus
$V$ splits in $K\langle V\rangle /K$.

(iii) Recall that $\Hom_{K}(V,L)$ has an $L\#D$-module structure by
$D$-conjugation (\ref{eq-internalHom}) and the map (\ref{eq-dense}) is
necessarily injective by Corollary \ref{cor-simple-module-algebras}.

(``If" part.) Since $L^{D}$ is a field, $\Hom_{K\# D}(V,L)$ is a free
$L^{D}$-module. By taking a dual basis, we can identify
$\Hom_{L^{D}}(\Hom_{K\# D}(V,L),L)$ with some power $L^{\Lambda}$ of $L$.
Then the injective $L$-module map
$L\otimes_{K}V\rightarrow\Hom_{L^{D}}(\Hom_{K\# D}(V,L),L)
\xrightarrow{\sim}L^{\Lambda}$ can be considered as an injective
$L\#D$-module map.

(``Only if" part.) By the definition, there is an $L\#D$-linear injection
$\varphi : L\otimes_{K}V\rightarrow L^{\Lambda}$ for some power $L^{\Lambda}$.
Let $\{f_{i}\}_{i\in\Lambda}\subset\Hom_{K\# D}(V,L)$
be the family of $K\#D$-module maps induced by
$\varphi(1\otimes v)=(f_{i}(v))_{i\in\Lambda}$ ($v\in V$).
Take an arbitrary element $w=\sum_{j}a_{j}\otimes v_{j}\in L\otimes_{K}V$.
If the image of $w$ by the map (\ref{eq-inj}) is $0$, then we have
$\sum_{j}a_{j}f_{i}(v_{j})=0$ for all $i\in\Lambda$ and so
$\varphi(w)=0$ ($\Leftrightarrow$ $w=0$). Thus the map (\ref{eq-inj})
is injective.
\end{proof}

\begin{proposition}
\label{prop-splitting-algebras2}
Let $K\subset L$ be an inclusion of AS $D$-module algebras and $V$ a
$K\# D$-module with finite $K$-rank $\rk_{K}(V)=r<\infty$.
Then 
\begin{equation}
\label{eq-rank}
\dim_{L^{D}}\Hom_{K\# D}(V,L)\leq r
\end{equation}
and the following are equivalent:
\begin{enumerate}
\renewcommand{\labelenumi}{\rm(\alph{enumi})}
\item $V$ splits in $L/K$;
\item $L\otimes_{L^{D}}\Hom_{K\# D}(V,L)\xrightarrow{\sim}\Hom_{K}(V,L)$;
\item $\dim_{L^{D}}\Hom_{K\# D}(V,L)=r$;
\item There is an isomorphism $L\otimes_{K}V\xrightarrow{\sim}L^{r}$
as $L\# D$-modules;
\item $L\otimes_{L^{D}}(L\otimes_{K}V)^{D}\xrightarrow{\sim}L\otimes_{K}V$;
\item $\dim_{L^{D}}(L\otimes_{K}V)^{D}=r$;
\item There is an injective $L\#D$-module map $L\otimes_{K}V
\rightarrow L^{n}$ for some integer $n$.
\end{enumerate}
\end{proposition}
\begin{proof}
The inequality (\ref{eq-rank}) follows since $L\otimes_{L^{D}}\Hom_{K\# D}(V,L)
\rightarrow\Hom_{K}(V,L)$ is an injective $L\#D$-module map whose
cokernel is a free $L$-module with finite rank.
Then ((a) $\Leftrightarrow$ (b) $\Leftrightarrow$ (c) $\Leftrightarrow$ (d))
follows from (the proof of) Proposition \ref{prop-splitting-algebras1} (iii).
((d) $\Leftrightarrow$ (e) $\Leftrightarrow$ (f)) is easily seen.
((d) $\Rightarrow$ (g) $\Rightarrow$ (a)) is trivial.
\end{proof}

\begin{lemma}
\label{lem-splitting-algebras3}
Let $K\subset L$ be an inclusion of AS $D$-module algebras,
$V$ a $K\#D$-module, and $W$ a $K\#D$-submodule of $V$ with finite
$K$-rank. If $V$ splits in $L/K$, then the restriction
$\mathrm{res.} : \Hom_{K\#D}(V,L)\rightarrow\Hom_{K\#D}(W,L)$ is surjective.
\end{lemma}
\begin{proof}
Consider the following $L\#D$-module map:
\[ \varphi : L\otimes_{L^{D}}\Hom_{K\#D}(V,L)\rightarrow
\Hom_{K}(W,L),\quad a\otimes f\mapsto af|_{W}. \]
Notice that $\im\varphi$ is a direct summand of $\Hom_{K}(W,L)$ as
an $L$-module. The transposed $L$-linear map of $\varphi$ is given by
\[ L\otimes_{K}W\hookrightarrow L\otimes_{K}V\rightarrow
\Hom_{L^{D}}(\Hom_{K\#D}(V,L),L), \]
which is injective by Proposition \ref{prop-splitting-algebras1} (iii).
Thus $\varphi$ is surjective. Since the functor
$(-)^{D}$ is exact, we have that $\varphi^{D}=\mathrm{res.}$ is surjective.
\end{proof}

We see that the functor $\Phi$ preserves splitting algebras:
\begin{lemma}
\label{descent-splitting}
Let $G_{1}\subset G$, $K_{1}\subset L_{1}$ be as in Lemma
\ref{lem-interlacing1}. Write $\Phi = \Phi_{G_{1}}$. Then, $L_{1}/K_{1}$
is a (minimal) splitting algebra for a $K_{1}\# D(G_{1})$-module $V_{1}$,
iff $\Phi(L_{1})/\Phi(K_{1})$ is a (minimal) splitting algebra
for the $\Phi(K_{1})\# D$-module $\Phi(V_{1})$.
\end{lemma}
\begin{proof}
This easily follows from Proposition \ref{prop-copying-functor} if
one notices that
\[ \Phi(K_{1}\lan V_{1}\ran)=\Phi(K_{1})\lan\Phi(V_{1})\ran \]
to see the equivalence on minimality.
\end{proof}

Let $K\subset L$ be an inclusion of AS $D$-module algebras.
For finitely many elements $u_{1},\dotsc,u_{m}$ in $L$, let $K\lan u_{1},
\dotsc, u_{m}\ran$ denote the smallest AS $D$-module subalgebra in $L$
including $K$ and $u_{1},\dotsc,u_{m}$.

\begin{definition}
\label{def-finitely-generated-extension}
$L/K$ is said to be {\em finitely generated} iff $L$ is of the form $K\lan
u_{1},\dotsc,u_{m}\ran$. This is equivalent to that $L_{1}/K_{1}$ is 
finitely generated, where $K_{1}=K/P\cap K$, $L_{1}=L/P$ for an arbitrarily
chosen $P\in\Omega(L)$.
\end{definition}

\begin{theorem}
\label{thm-characterization}
Let $K\subset L$ be as above. Suppose $K^{D}=L^{D}$. Then the following are
equivalent:
\begin{enumerate}
\renewcommand{\labelenumi}{\rm(\alph{enumi})}
\item $L/K$ is a finitely generated PV extension;
\item $L/K$ is a minimal splitting algebra for a cyclic $K\# D$-module
of finite $K$-rank;
\item $L/K$ is a minimal splitting algebra for a $K\# D$-module of finite
$K$-rank;
\item $L=K\lan x_{ij}\ran$, where $X=(x_{ij})_{i,j}$ is a $GL_{n}$-primitive
in Kolchin's sense \cite{Kolchin1973}: $X\in GL_{n}(L)$, and for every
$d\in D$, $(dX)X^{-1}\in M_{n}(K)$ with $dX=(dx_{ij})_{i,j}$.
\end{enumerate}
\end{theorem}
\begin{proof}
We write $k=K^{D}$ ($=L^{D}$).

(a) $\Rightarrow$ (b). By Lemmas \ref{lem-interlacing1} and
\ref{descent-splitting}, we may assume that $K$ is a field. Suppose that
$(L/K,A,H)$ is a finitely generated PV extension. By Proposition
\ref{prop-interlacing2}, we have a finitely generated PV extension
$(L_{1}/K,A_{1},\bar{H})$ of module fields over $C:=D(G_{P})$ with $P\in
\Omega(L)$, such that $L=\Phi(L_{1})$, $A=\Phi(A_{1})$.

There exist those finitely many elements $u_{1},\dotsc,u_{m}$ in $A$ which
span an $H$-subcomodule over $k$, and satisfy $L=K\lan u_{1},\dotsc,u_{m}\ran$.
(In fact, take some $x_{1},\dotsc,x_{l}\in A$ which satisfy
$L=K\lan x_{1},\dotsc,x_{l}\ran$. Then there exists a finite dimensional
$H$-subcomodule $U\subset A$ such that $x_{1},\dotsc,x_{l}\in U$
(see \cite[(3.3)]{Waterhouse}). Choose $u_{1},\dotsc,u_{m}$ which span $U$
over $k$.) Set an element $\bu=(u_{1},\dotsc,u_{m})$
in $A^{m}$, and let $V=(K\# D)\bu$, the cyclic $K\# D$-submodule generated 
by $\bu$. Since $L\otimes_{K}A\simeq L\otimes_{k}H$, we see that $L/K$ is
a minimal splitting algebra for $A^{m}$, and hence for $V$.

It remains to prove that the $K$-dimension $\dim_{K}(V)$ is finite.
It suffices to prove that the natural image $V(P)$, say, of $V$ under the
projection $A^{m}\rightarrow A_{1}^{m}$ has a finite $K$-dimension, since
$V$ is naturally embedded into $\prod_{P\in\Omega(L)}V(P)$. Let $g_{1},\dotsc,
g_{s}$ be a system of representatives of the right cosets $G_{P}\backslash G$.
Then we have
\[ V = \sum_{i=1}^{s}(K\# C)g_{i}\bu. \]
Fix an $i\in\{1,\dotsc,s\}$, and let $\bw = (w_{1},\dotsc,w_{m})\in A_{1}^{m}$
denote the natural image of $g_{i}\bu$. It suffices to prove that
$W:=(K\# C)\bw$ has a finite $K$-dimension.
By re-numbering we have a $k$-basis, $w_{1},\dotsc,
w_{r}$ ($r\leq m$), of the $k$-subspace in $A_{1}$ spanned by $w_{1},\dotsc,
w_{m}$. There is a rank $r$ matrix $T$ with entries in $k$, such that $\bw
=\bw'T$ with $\bw'=(w_{1},\dotsc,w_{r})$. It suffices to prove that $W':=
(K\# C)\bw'$ has a finite $K$-dimension, since $W'\simeq W$ under the right
multiplication by $T$.

Notice that for any $g\in G$, $gu_{1},\dotsc,gu_{m}$ span an $H$-subcomodule
in $A$ since the comodule structure map $\theta : A\rightarrow A\otimes_{k}H$
is $D$-linear. It then follows that $w_{1},\dotsc,w_{r}$ form a $k$-basis of an
$\bar{H}$-subcomodule in $A_{1}$. By applying Proposition \ref{Wronskian}
for $w_{1},\dotsc,w_{r}\in L_{1}$, there exist $r$ elements
$h_{1},\dotsc,h_{r}\in C$ such that $(h_{i}(w_{j}))_{i,j}$ is an invertible
matrix. We claim that $(c\bw')(h_{i}(w_{j}))^{-1}_{i,j}\in K^{r}$
for all $c\in C$. If it follows, then $W'=Kh_{1}(\bw')+\dotsb +Kh_{r}(\bw')$
and hence we have $\dim_{K}W'<\infty$, as desired.

Let $\theta_{1} : A_{1}\rightarrow A_{1}\otimes_{k}\bar{H}$ be the
comodule structure map associated to the PV extension $L_{1}/K$.
Write
\[ \theta_{1}(w_{j})=\sum_{s=1}^{r}w_{s}\otimes_{k}z_{sj}\quad (z_{sj}\in H,
\ j=1,\dotsc,r). \]
By applying $\mu$-isomorphism $A_{1}\otimes_{k}H\xrightarrow{\sim}
A_{1}\otimes_{K}A_{1}$ in each side, we have
\begin{equation}
\label{Wronskian-Hopf-algebra}
1\otimes_{K}w_{j}=\sum_{s=1}^{r}(w_{s}\otimes_{K}1)z_{sj}\quad
\mbox{in $A_{1}\otimes_{K}A_{1}$}.
\end{equation}
Hence
\[ 1\otimes_{K}h_{i}(w_{j})=\sum_{s=1}^{r}(h_{i}(w_{s})\otimes_{K}1)z_{sj}
\quad (i,j=1,\dotsc,r), \]
i.e.\ $1\otimes_{K}(h_{i}(w_{j}))_{i,j}=((h_{i}(w_{j}))_{i,j}\otimes_{K}1)Z$
with $Z=(z_{ij})_{i,j}$. Since $(h_{i}(w_{j}))_{i,j}$ is invertible,
we have 
\[ Z=((h_{i}(w_{j}))_{i,j}^{-1}\otimes_{K}1)(1\otimes_{K}(h_{i}(w_{j}))_{i,j})
\in GL_{r}(L_{1}\otimes_{K}L_{1}). \]
On the other hand, recalling (\ref{Wronskian-Hopf-algebra}) we have
\[ 1\otimes_{K}(c\bw')=((c\bw')\otimes_{K}1)Z \]
for all $c\in C$. Thus, by multiplying $1\otimes_{K}(h_{i}(w_{j}))^{-1}_{i,j}
=Z^{-1}((h_{i}(w_{j}))_{i,j}^{-1}\otimes_{K}1)$ from the right,
\[ 1\otimes_{K}(c\bw')(h_{i}(w_{j}))_{i,j}^{-1}=((c\bw')\otimes_{K}1)ZZ^{-1}
((h_{i}(w_{j}))_{i,j}^{-1}\otimes_{K}1)=(c\bw')(h_{i}(w_{j}))^{-1}_{i,j}
\otimes_{K}1 \]
for all $c\in C$. This implies the claim above.

(b) $\Rightarrow$ (c). This is trivial.

(c) $\Rightarrow$ (d). Suppose that $L/K$ is a minimal splitting algebra for
$V$ with finite $K$-free basis $v_{1},\dotsc,v_{n}$.
By Proposition \ref{prop-splitting-algebras2} (c), we have a $k$-basis
$f_{1},\dotsc,f_{n}$ in $\Hom_{K\# D}(V,L)$. Define
\begin{equation}
\label{making-GL-primitive}
X=(x_{ij})_{i,j}=(f_{j}(v_{i}))_{i,j},\quad \bv = {}^{t}(v_{1},\dotsc,v_{n}).
\end{equation}
Then we have $X\in GL_{n}(L)$ since there is an $L$-module isomorphism
\[ L^{n}\simeq L\otimes_{K}V\xrightarrow{\sim}
\Hom_{k}(\Hom_{K\# D}(V,L),L)\simeq L^{n} \]
which is precisely the multiplication of $X$ (see the proof of
Proposition \ref{prop-splitting-algebras1}, \ref{prop-splitting-algebras2}).
If we write $dv_{i}=\sum_{s=1}^{n}c_{is}(d)v_{s}$ ($c_{is}(d)\in K$)
for $d\in D$, then we have $dx_{ij}=f_{j}(dv_{i})
=\sum_{s=1}^{n}c_{is}(d)f_{j}(v_{s})=\sum_{s=1}^{n}c_{is}(d)x_{sj}$.
This implies that $X$ is $GL_{n}$-primitive such that
\begin{equation}
\label{in-fact-GL-primitive}
(dX)X^{-1}\bv = d\bv\quad (d\in D),
\end{equation}
i.e.\ $(dX)X^{-1}=(c_{ij}(d))_{i,j}\in M_{n}(K)$.
By the definition, we have $L=K\lan V\ran=K\lan x_{ij}\ran$.

(d) $\Rightarrow$ (a). This is shown by modifying
\cite[Example 2.5c]{Takeuchi1989} as follows.

Put $Y=(y_{ij})_{i,j}=X^{-1}$ and $A=K[x_{ij},y_{ij}]$. First we shall show
that $A$ is a $D$-module subalgebra of $L$.
Define $\phi\in\Hom_{R}(D,M_{n}(K))$ by $\phi_{d}=d(X)X^{-1}=d(X)Y$ ($d\in D$).
Since $\phi_{g}^{-1}=Xg(Y)=g(g^{-1}(X)X^{-1})\in GL_{n}(K)$ for all $g\in G$,
$\phi$ is convolution-invertible in $\Hom_{R}(D,M_{n}(K))$ by
\cite[Corollary 9.2.4]{Sweedler1969}. We see that the $\psi\in
\Hom_{R}(D,M_{n}(L))$ given by $\psi_{d}=Xd(Y)$, is the inverse of $\phi$,
and so $\psi\in\Hom_{R}(D,M_{n}(K))$.
This implies that $A$ is a $D$-module subalgebra of $L$.
Since $Q(A)$, the total quotient ring of $A$, is an AS $D$-module subalgebra
of $L$ containing $K$ and $x_{11},\dotsc,x_{nn}$ (recall Corollary
\ref{total-quotient-ring}), we have $Q(A)=K\lan x_{ij}\ran=L$.

Put
\[ Z=(Y\otimes_{K}1)(1\otimes_{K}X),\ Z^{-1}=(1\otimes_{K}Y)(X\otimes_{K}1)
\in GL_{n}(A\otimes_{K}A). \]
For all $d\in D$,
\begin{eqnarray*}
d(Z) & = &
{\displaystyle\sum_{(d)}}(d_{(1)}(Y)\otimes_{K}1)(1\otimes_{K}d_{(2)}(X)) \\
 & = & {\displaystyle\sum_{(d)}}(Y\psi_{d_{(1)}}\otimes_{K}1)
(1\otimes_{K}\phi_{d_{(2)}}X) \\
 & = & {\displaystyle\sum_{(d)}}(Y\otimes_{K}1)(\psi_{d_{(1)}}\otimes_{K}1)
(1\otimes_{K}\phi_{d_{(2)}})(1\otimes_{K}X) \\
 & = & {\displaystyle\sum_{(d)}}(Y\otimes_{K}1)(\psi_{d_{(1)}}\otimes_{K}1)
(\phi_{d_{(2)}}\otimes_{K}1)(1\otimes_{K}X) \\
 & = & {\displaystyle\sum_{(d)}}(Y\otimes_{K}1)
(\psi_{d_{(1)}}\phi_{d_{(2)}}\otimes_{K}1)(1\otimes_{K}X) \\
 & = & \varepsilon(d)Z.
\end{eqnarray*}
Thus $Z$ has entries in $H:=(A\otimes_{K}A)^{D}$. Similarly we have that
the entries in $Z^{-1}$ are also in $H$ and hence $Z\in GL_{n}(H)$.
Then,
\[ 1\otimes_{K}X=(X\otimes_{K}1)Z,\ 1\otimes_{K}Y=Z^{-1}(Y\otimes_{K}1)
\in GL_{n}(A\cdot H). \]
This implies $A\otimes_{K}A=A\cdot H$. Therefore $(L/K,A,H)$ is a PV
extension.
\end{proof}

\begin{remark}
\label{rem-compatibility}
Keep the notation just as above. 

(i) Write $Z=(z_{ij})$, $Z^{-1}=(w_{ij})$.
Then $A\otimes_{K}A=A[z_{ij},\det(z_{ij})^{-1}]$ and $H
=k[z_{ij},\det(z_{ij})^{-1}]=k[z_{ij},w_{ij}]$.
Taking $\mu,\theta$ as in Proposition \ref{prop-PV-Hopf-algebra}, we see
\[ \begin{array}{lll}
\theta(X) & = & \mu^{-1}(1\otimes_{K}X)=\mu^{-1}((X\otimes_{K}1)Z)
=\mu^{-1}\left({\displaystyle\sum_{s=1}^{n}}(x_{is}\otimes_{K}1)z_{sj}
\right)_{i,j} \\
 & = & \left({\displaystyle\sum_{s=1}^{n}}x_{is}\otimes_{k}z_{sj}\right)_{i,j}
=(X\otimes_{k}1)(1\otimes_{k}Z).
\end{array} \]
This is often written like
\begin{equation}
\label{comodule-structure}
\theta(X)=X\otimes_{k}Z.
\end{equation}
It follows that the Hopf algebra structure of $H$ is given by
\[ \Delta(Z)=Z\otimes_{k}Z,\qquad \vep(Z)=I,\quad S(Z)=Z^{-1}, \]
here $I$ denotes the identity matrix; see \cite[(3.2), Corollary]{Waterhouse}.
We have a Hopf algebra surjection,
\[ k[GL_{n}]=k[T_{ij},\det(T_{ij})^{-1}]\twoheadrightarrow H,\quad
T_{ij}\mapsto z_{ij}, \]
which gives a closed embedding $\bG(L/K)\rightarrow \bGL_{n}$ of
affine $k$-group schemes.

(ii) Suppose that $D=R[\tau,\tau^{-1}]$, the group algebra of the free
abelian group of rank $1$, and $K$ is a field; $K$ is then an inversive
difference field.
A difference system $\tau\by=B\by$ with $B\in GL_{n}(K)$
arises uniquely from a $K\# D$-module
of $K$-dimension $n$, together with its $K$-basis. We see from
(\ref{in-fact-GL-primitive}) that the $X$ in (\ref{making-GL-primitive}) is
a fundamental matrix \cite[Definition 1.4]{vanderPut-Singer1997}
for the difference system arising from the $V$ and the $\bv$ above,
and so that $A$ is the Picard-Vessiot ring
\cite[Definition 1.5]{vanderPut-Singer1997} for the system.
It will follow from Theorems
\ref{thm-characterization}, \ref{thm-existence} that if $k$ ($=K^{D}$)
is algebraically closed, a Picard-Vessiot ring for any difference system as
above uniquely exists, and is given by such an $A$ as above.
\end{remark}

\begin{corollary}
\label{finitely-generated-PV-extensions}
Let $(L/K,A,H)$ be a PV extension of AS $D$-module algebras. The following
are equivalent:
\begin{enumerate}
\renewcommand{\labelenumi}{\rm(\alph{enumi})}
\item $L/K$ is finitely generated;
\item $L$ is the total quotient ring of a finitely generated $K$-subalgebra
in $L$;
\item $A$ is finitely generated as a $K$-algebra;
\item $H$ is finitely generated as a $k$-algebra.
\end{enumerate}
\end{corollary}
\begin{proof}
((a) $\Rightarrow$ (c) $\Rightarrow$ (b) $\Rightarrow$ (a)) and
((a) $\Rightarrow$ (d)) follow by the proof of
Theorem \ref{thm-characterization}.
If $H$ is a finitely generated $k$-algebra, then we have an ascending
chain condition for Hopf ideals of $H$. Hence we have an ascending chain
condition for intermediate AS $D$-module algebras of $L/K$, which implies (a).
\end{proof}

\begin{corollary}
\label{PV-splitting}
Let $K\subset L$ be an inclusion of AS $D$-module algebras such that
$K^{D}=L^{D}=:k$. Then $L/K$ is a PV extension iff it is a minimal
splitting algebra for such a $K\# D$-module $V$ that is a directed union,
$V=\bigcup_{\lambda}V_{\lambda}$, of $K\# D$-submodules $V_{\lambda}$ of
finite $K$-rank.
\end{corollary}
\begin{proof}
This follows in the same way as \cite[Corollary 3.5]{Takeuchi1989}.
We include the proof for convenience.

(``Only if" part.) Let $(L/K,A,H)$ be a PV extension. Then $H$ is a
directed union of Hopf subalgebras which are finitely generated $k$-algebras
(see \cite[(3.3)]{Waterhouse}).
It follows by Proposition \ref{prop-Galois-correspondence3} and 
Corollary \ref{finitely-generated-PV-extensions} that
$L$ is a directed union, say $L=\bigcup_{\lambda}L_{\lambda}$,
of AS $D$-module subalgebras which are finitely generated PV extensions
over $K$. By Theorem \ref{thm-characterization}, each $L_{\lambda}/K$ is
a minimal splitting algebra for a $K\#D$-module $V_{\lambda}$ of finite
$K$-rank. Then $L/K$ is a minimal splitting algebra for the direct sum
$V=\bigoplus_{\lambda}V_{\lambda}$.

(``If" part.) Suppose that $L/K$ is a minimal splitting algebra for a
$K\#D$-module $V=\bigcup_{\lambda}V_{\lambda}$, a directed union of
$K\#D$-modules $V_{\lambda}$ of finite $K$-rank. Since every $V_{\lambda}$
splits in $L/K$, each $L_{\lambda}:=K\lan V_{\lambda}\ran$ is a minimal
splitting algebra for $V_{\lambda}$ and is a finitely generated PV extension
over $K$ by Theorem \ref{thm-characterization}.
By Lemma \ref{lem-splitting-algebras3}, the union
$\bigcup_{\lambda}L_{\lambda}$ is a directed union of AS $D$-module subalgebras
of $L$. Thus $\bigcup_{\lambda}L_{\lambda}$ is an AS $D$-module subalgebra
of $L$ by Lemma \ref{AS-subalgebra}.
For every $f\in\Hom_{K\#D}(V,L)$, we have $f(V)
=\bigcup_{\lambda}f(V_{\lambda})\subset\bigcup_{\lambda}L_{\lambda}$.
Hence $L=K\lan V\ran=\bigcup_{\lambda}L_{\lambda}$. Let $A_{\lambda}$ 
(resp.\ $H_{\lambda}$) be the principal algebra (resp.\ the Hopf algebra)
for $L_{\lambda}/K$. Then one sees that $(L/K,\bigcup_{\lambda}A_{\lambda},
\bigcup_{\lambda}H_{\lambda})$ is a PV extension.
\end{proof}

\begin{theorem}
\label{thm-existence}
Let $K$ be an AS $D$-module algebra such that the field $K^{D}$ of
$D$-invariants is algebraically closed. Let $V$ be a $K\# D$-module of
finite $K$-rank. Then there exists an AS $D$-module algebra $L$ including $K$
such that $K^{D}=L^{D}$, and $L/K$ is a (necessarily finitely generated)
minimal splitting algebra for $V$. Such an algebra is unique up to $D$-linear
isomorphism of $K$-algebras.
\end{theorem}

To prove this, we need the following:
\begin{lemma}
\label{Levelt}
Let $K$ be an AS $D$-module algebra. Let $A$ be a simple $D$-module algebra,
and let $L=Q(A)$ be the total quotient ring of $A$; by Lemma
\ref{lem-localization}, $L$ is uniquely a $D$-module algebra.
If $A$ is finitely generated as a $K$-algebra, then $L^{D}/K^{D}$
is an algebraic extension of fields.
\end{lemma}
\begin{proof}
We follow Levelt \cite[Appendix]{Levelt} for this proof.
If $x\in L^{D}$, then $(A:x)=\{a\in A\;|\;ax\in A\}$ is a $D$-stable ideal.
Since this contains a non-zero divisor, we have that $(A:x)=A$, and so
$A^{D}=L^{D}$.

If $A$ is finitely generated, then it is noetherian.
By Proposition \ref{prop-NS}, we may suppose that $K$ is a field (and $A$ is a
domain). If $P\subset A$ is a maximal ideal, then the field $A^{D}$ is
included in the field $A/P$, which is algebraic over $K$. Therefore if
$x\in A^{D}$, it is algebraic over $K$.
Let $\varphi(T)=T^{n}+c_{1}T^{n-1}+\dotsb +c_{n}$ denote the minimal
polynomial of $x$ over $K$. Since for any $d\in D$,
$\varepsilon(d)T^{n}+(dc_{1})T^{n-1}+\dotsb +dc_{n}$ has $x$ as a root,
each $c_{i}\in K^{D}$ by the minimality of $\varphi(T)$. Thus $x$ is algebraic
over $K^{D}$.
\end{proof}

\begin{proof}[{Proof of Theorem \ref{thm-existence}}]
{\it Existence}; this is proved by modifying the proof of
\cite[Theorem 4.5]{Takeuchi1989}, as follows.
Let $v_{1},\dotsc,v_{r}$ be a $K$-basis for $V$. For $d\in D$, write
\[ dv_{i}=\sum_{s=1}^{r}c_{is}(d)v_{s} \]
with $c_{is}(d)\in K$. Define a $D$-module algebra structure on $K[X_{ij}]$,
the polynomial $K$-algebra in $r^{2}$ indeterminates, by
\[ d(X_{ij})=\sum_{s=1}^{r}c_{is}(d)X_{sj}\qquad (d\in D). \]
For each $g\in G$, we see that $(c_{ij}(g))$ is an invertible matrix:
$(c_{ij}(g))^{-1}=(gc_{ij}(g^{-1}))$.
Thus $\det(c_{ij}(g))$ is invertible in $K$ for each $g\in G$, and
the $D$-module algebra structure of $K[X_{ij}]$ is uniquely extended to
$F=K[X_{ij},\det(X_{ij})^{-1}]$ by Lemma \ref{lem-localization}. 
Let $I$ be a maximal $D$-stable ideal of $F$, and put $A=F/I$. Since $K$ is
simple, $I\cap K=0$. Hence $A$ is a noetherian simple $D$-module algebra
including $K$. Let $L$ be the total quotient ring of $A$; this is an AS
$D$-module algebra by Proposition \ref{prop-NS} and Lemma
\ref{lem-localization}.
By Lemma \ref{Levelt}, we have $L^{D}=K^{D}$. Let $x_{ij}$ denote the image
of $X_{ij}$ in $A$, and define $K$-linear maps $f_{j} : V\rightarrow L$
($j=1,\dotsc,r$) by $f_{j}(v_{i})=x_{ij}$. Then these maps are in
$\Hom_{K\# D}(V,L)$, and are linearly independent over $L^{D}$, since
$(x_{ij})_{i,j}\in GL_{r}(L)$. Therefore, $L/K$ is a minimal splitting algebra
for $V$ by Lemma \ref{prop-splitting-algebras2} (c).

{\it Uniqueness}; also this proof is essentially the same as the proof given in
\cite[Theorem 4.6]{Takeuchi1989}. Let $L_{1}/K$ and $L_{2}/K$ are
two minimal splitting algebra for $V$ such that $L_{1}^{D}=L_{2}^{D}=K^{D}=k$.
By Theorem \ref{thm-characterization},
$L_{1}/K$ and $L_{2}/K$ are finitely generated PV extensions.
Let $A_{i}$ be the principal algebra for $L_{i}/K$ ($i=1,2$) respectively.
Put $A=A_{1}\otimes_{K}A_{2}$ and let $I$ be a maximal $D$-stable ideal of $A$.
We see $A_{1},A_{2}$ are noetherian simple $D$-module algebras
which are finitely generated $K$-algebras by the proof of Theorem
\ref{thm-characterization} (d) $\Rightarrow$ (a), and by
Proposition \ref{prop-principal-algebra} (i).
Thus $A_{i}\cap I=0$ ($i=1,2$) and hence $A_{1},A_{2}$ are identified with
$D$-module subalgebras of $A/I$. Let $L$ be the total quotient ring of $A/I$.
By Lemma \ref{lem-localization} and Proposition \ref{prop-NS},
$L$ is an AS $D$-module algebra since $A/I$ is a noetherian simple
$D$-module algebra. Furthermore, since $A/I$ is
a finitely generated $K$-algebra, we have $L^{D}=K^{D}=k$
by Lemma \ref{Levelt}. Let $\chi_{i} : L_{i}\hookrightarrow L$ ($i=1,2$)
denote the induced inclusions of $D$-module algebras over $K$.
The injective $k$-linear maps $\Hom_{K\# D}(V,L_{i})\rightarrow
\Hom_{K\# D}(V,L)$ ($i=1,2$) are precisely isomorphisms by Proposition
\ref{prop-splitting-algebras2}. Therefore $f(V)\subset
\chi_{1}(L_{1})\cap\chi_{2}(L_{2})$ for all $f\in\Hom_{K\# D}(V,L)$.
Since $L_{i}$ are generated over $K$ by the image of all $f\in
\Hom_{K\# D}(V,L_{i})$,
we have $\chi_{1}(L_{1})=\chi_{2}(L_{2})=K\langle V\rangle$ in $L$.
Thus we have $\chi_{2}^{-1}\circ\chi_{1} : L_{1}\xrightarrow{\sim}L_{2}$,
a $D$-module algebra isomorphism over $K$.
\end{proof}

Let $K$ be an AS $D$-module algebra. We have the $K^{D}$-abelian symmetric
tensor category $({}_{K\# D}\cM,\otimes_{K},K)$. Let $V$ be an object in
${}_{K\# D}\cM$ of finite $K$-rank. Then the $K$-linear dual $V^{\ast}:=
\Hom_{K}(V,K)$ is a dual object under the $D$-conjugation; see
(\ref{eq-internalHom}). Thus the tensor full subcategory
${}_{K\# D}\cM_{\mathrm{fin}}$ consisting of the finite $K$-rank objects
is rigid. Let $\{\{V\}\}$ denote the abelian, rigid tensor full subcategory
of ${}_{K\# D}\cM$ generated by $V$, that is, the smallest full subcategory
containing $V$ that is closed under subquotients, finite direct sums,
tensor products and duals. Thus an object in $\{\{V\}\}$ is precisely
a subquotient of some finite direct sum $W_{1}\oplus\dotsb\oplus W_{r}$,
where each $W_{i}$ is the tensor product of some copies of $V,V^{\ast}$;
see \cite[Theorem 2.33]{vanderPut-Singer2003} also for comparing with
the following.

\begin{theorem}
\label{equivalence-to-representations}
Let $(L/K,A,H)$ be a finitely generated PV extension of AS $D$-module
algebras. By Theorem \ref{thm-characterization}, we have such a
$K\# D$-module $V$ of finite $K$-rank for which $L/K$ is a minimal splitting
algebra.

{\rm (i)} Let $W\in\{\{V\}\}$. Regard the $A\otimes_{K}W$ as a right
$H$-comodule with the structure induced by $A$. Then $(A\otimes_{K}W)^{D}$
is an $H$-subcomodule with $k$-dimension $\rk_{K}(W)$.

{\rm (ii)} $W\mapsto (A\otimes_{K}W)^{D}$ gives a $k$-linear equivalence
\[ \{\{V\}\}\approx\cM_{\mathrm{fin}}^{H} \]
of symmetric tensor categories, where $\cM_{\mathrm{fin}}^{H}
=(\cM_{\mathrm{fin}}^{H},\otimes_{k},k)$ denotes the rigid symmetric tensor
category of finite-dimensional right $H$-comodules; notice that this is
isomorphic to the category $\mathrm{Rep}_{\bG(L/K)}$ of the same kind,
consisting of finite-dimensional linear representations of the PV group scheme
$\bG(L/K)=\Spec H$.
\end{theorem}
\begin{proof}
Put $D_{k}=D\otimes_{R}k$, a cocommutative Hopf algebra over $k$, and consider
$D_{k}$ as a right $H$-comodule algebra with the trivial structure map
$d\mapsto d\otimes 1$.
Regard naturally $A$ as an algebra in the symmetric tensor category
$({}_{D_{k}}\cM^{H},\otimes_{k},k)$ of right $(H,D_{k}^{\mathrm{op}})$-Hopf
modules (see \cite[\S 8.5]{Montgomery}); its objects are $D_{k}$-modules $N$
which has a $D_{k}$-linear, right $H$-comodule structure
$\rho_{N} : N\rightarrow N\otimes_{k}H$. We then have the
symmetric tensor category ${}_{A}({}_{D_{k}}\cM^{H})$ of $A$-modules in
${}_{D_{k}}\cM^{H}$, which is denoted by $({}_{A\# D}\cM^{H},\otimes_{A},A)$;
this is $k$-abelian. Define $k$-linear functors
\[ \cM^{H}\substack{\xrightarrow{\Theta_{1}} \\ \xleftarrow[\Xi_{1}]{}}
{}_{A\# D}\cM^{H}\substack{\xrightarrow{\Theta_{2}} \\ \xleftarrow[\Xi_{2}]{}}
{}_{K\# D}\cM \]
by 
\begin{eqnarray*}
\Theta_{1}(U) & = & A\otimes_{k}U;\quad \mbox{$H$ coacts codiagonally,} \\
\Xi_{1}(N) & = & N^{D}, \\
\Theta_{2}(N) & = & N^{\coinv H}\quad
(=\{n\in N\;|\;\rho_{N}(n)=n\otimes_{k}1\}), \\
\Xi_{2}(W) & = & A\otimes_{K}W;\quad \mbox{$H$ coacts on $A$.}
\end{eqnarray*}
We see that $\Theta_{1}$ and $\Xi_{2}$ are symmetric tensor functors with
the obvious tensor structures. Moreover by \cite{Schneider},
$\Theta_{2}$ and $\Xi_{2}$ are quasi-inverses of each other, since $A/K$ is
$H$-Galois by Proposition \ref{prop-PV-Hopf-algebra} (iii). Since
$A^{D}=k$, $\Xi_{1}\circ\Theta_{1}$ is isomorphic to the identity functor.
Suppose $N\in {}_{A\# D}\cM^{H}$. Since $A$ is simple by Corollary
\ref{prop-principal-algebra} (i), we see from Corollary
\ref{cor-simple-module-algebras}
that the morphism in ${}_{A\# D}\cM^{H}$
\[ \mu_{N} : \Theta_{1}\circ\Xi_{1}(N)=A\otimes_{k}N^{D}\rightarrow N \]
is an injection. Let $\cN$ denote the full subcategory of
${}_{A\# D}\cM^{H}$ consisting of those $N$ for which $\mu_{N}$ is an
isomorphism. Since each $\Theta_{1}(U)$ is in $\cN$, $\Theta_{1}$ gives an
equivalence
\[ \cM^{H}\approx\cN. \]
Necessarily, $\cN$ is closed under tensor products, and this is an equivalence
of symmetric tensor categories.

Since $A\otimes_{K}V\simeq A^{n}$ ($n=\rk_{K}(V)$) in ${}_{A\# D}\cM$,
$\Xi_{2}(V)=A\otimes_{K}V\in\cN$. We see that $\Theta_{1}$ is exact, and $\cN$
is closed under subquotients. Therefore for (ii), it suffices to prove that
\[ \tilde{V} := \Xi_{1}\circ\Xi_{2}(V)=(A\otimes_{K}V)^{D} \]
generates $\cM_{\mathrm{fin}}^{H}$. Let $v_{1},\dotsc,v_{n}$ be a $K$-free
basis of $V$, and define $X,\bv$ as in (\ref{making-GL-primitive}). We see
from (\ref{in-fact-GL-primitive}) that the entries in $\tilde{\bv}:=X^{-1}
\otimes_{K}\bv$ ($\in (A\otimes_{K}V)^{n}$) are $D$-invariant,
and hence form a $k$-basis in $\tilde{V}$. By (\ref{comodule-structure}),
the $H$-comodule structure $\rho_{\tilde{V}} : \tilde{V}\rightarrow
\tilde{V}\otimes_{k}H$ on $\tilde{V}$ is given by
\[ \rho_{\tilde{V}}({}^{t}\tilde{\bv})
={}^{t}\tilde{\bv}\otimes_{k}{}^{t}Z^{-1}, \]
where ${}^{t}$ denotes the transpose of matrices.
This means that the coefficient $k$-space of $\tilde{V}$ is the subcoalgebra
in $H$ spanned by the entries $w_{ij}$ in ${}^{t}Z^{-1}$. Since $w_{ij}$
together with the entries $S(w_{ij})$ in $Z$ generate the $k$-algebra $H$
(see the proof of Theorem \ref{thm-characterization} (d) $\Rightarrow$
(a)), $\tilde{V}$ generates $\cM_{\mathrm{fin}}^{H}$;
see \cite[(3.5)]{Waterhouse}. This proves part (ii).

If $W\in\{\{V\}\}$, then $\Xi_{2}(W)\in\cN$, and so
\[ \dim_{k}(A\otimes_{K}W)^{D}=\rk_{A}(A\otimes_{K}W)=\rk_{K}(W). \]
This proves part (i).
\end{proof}

\section{Liouville extensions}
\label{sec-Liouville}

Finally we define the notion of Liouville extensions and show the
solvability theorem. As is described in Introduction, we should
define {\em Liouville group schemes} and study how strong the definition
is.

\subsection{Liouville group schemes}
\label{secLGS}

\begin{definition}
\label{defLGS}
Let $\bG$ be an algebraic affine group scheme over a field $k$.

(1) We say $\bG$ is {\em ($k$-)Liouville} (cf.\ \cite[p.\ 374]{Kolchin1973})
iff there exists a normal chain of closed subgroup schemes
\begin{equation}
\label{LNC}
\bG=\bG_{0}\vtr\bG_{1}\vtr\dotsb\vtr\bG_{r}=\{1\}
\end{equation}
such that each $\bG_{i-1}/\bG_{i}$ ($i=1,\dotsc,r$) is at least one of the
following types: finite etale, a closed subgroup scheme of $\Ga$,
or a closed subgroup scheme of $\Gm$.
In this case, we call (\ref{LNC}) a {\em Liouville normal
chain} (LNC).

(2) In (\ref{LNC}), if each $\bG_{i-1}/\bG_{i}$ is merely a closed subgroup
scheme of $\Ga$ or a closed subgroup scheme of $\Gm$, then we call it
a {\em restricted Liouville normal chain} (RLNC).
\end{definition}

We use the following abbreviation of some types on group schemes
which arise as factor group schemes in an LNC:
we say $\bG$ is of {\em $\Ga$-type} (resp., {\em $\Gm$-type}) iff it is a
closed subgroup scheme of $\Ga$ (resp., $\Gm$), and a group scheme of
{\em RL-type} (resp., {\em L-type}) means that it is of $\Ga$-type or
$\Gm$-type (resp., RL-type or finite etale).

\begin{lemma}
\label{solvability1}
{\rm (1)} If $\bG$ is Liouville (resp., has an RLNC), then every closed
subgroup scheme of $\bG$ is Liouville (resp., has an RLNC).
Especially $\bG$ is Liouville iff the connected component $\bG^{\circ}$
is Liouville.

{\rm (2)} Let $\bH$ be a normal closed subgroup scheme of $\bG$. Then
$\bG$ is Liouville (resp., has an RLNC) iff both $\bH$ and $\bG/\bH$ are
Liouville (resp., have an RLNC).

{\rm (3)} If $\bG$ is connected Liouville, then $\bG$ is solvable.
\end{lemma}
\begin{proof}
First we take an LNC (resp., an RLNC):
$\bG=\bG_{0}\vtr\bG_{1}\vtr\dotsb\vtr\bG_{r}=\{1\}$
in each proof of (1), ``only if" part of (2), and (3).

(1) Let $\bH$ be a closed subgroup scheme of $\bG$ and put $\bH_{i}:=
\bH\cap\bG_{i}$ ($i=0,\dotsc,r$). Then we have $\bH_{0}=\bH$ and
$\bH_{i}=\bH_{i-1}\cap\bG_{i}=\Ker(\bH_{i-1}\rightarrow\bG_{i-1}/\bG_{i})$
for $i=1,\dotsc,r$. It follows that $\bH_{i-1}\vtr\bH_{i}$ and
$\bH_{i-1}/\bH_{i}$ is a closed subgroup scheme of $\bG_{i-1}/\bG_{i}$
for $i=1,\dotsc,r$. Therefore $\bH=\bH_{0}\vtr\bH_{1}\vtr\dotsb\vtr\bH_{r}
=\{1\}$ is an LNC (resp., an RLNC).

(2) (``Only if" part.) $\bH$ is Liouville (resp., has an RLNC) by (1).
Put $\bF_{i}:=\bG_{i}/\bH\cap\bG_{i}$ ($i=0,\dotsc,r$).
Each $k[\bF_{i-1}/\bF_{i}]$ is identified with a Hopf subalgebra of
$k[\bG_{i-1}/\bG_{i}]$, since the quotient morphism $\bG_{i-1}\rightarrow
\bF_{i-1}\rightarrow\bF_{i-1}/\bF_{i}$ goes through $\bG_{i-1}\rightarrow
\bG_{i-1}/\bG_{i}$; consider the following commutative diagram:
\[ \begin{CD}
\{1\} @>>> \bG_{i} @>>> \bG_{i-1} @>>> \bG_{i-1}/\bG_{i} @>>> \{1\} \quad
\mbox{(exact)} \\
@. @VV{\mbox{\scriptsize quotient}}V  @VV{\mbox{\scriptsize quotient}}V
@VV{{}^{\exists}\mbox{\scriptsize quotient}}V  @. \\
\{1\} @>>> \bF_{i} @>>> \bF_{i-1} @>>> \bF_{i-1}/\bF_{i} @>>> \{1\} \quad
\mbox{(exact)}.
\end{CD} \]
Thus each $\bF_{i-1}/\bF_{i}$ is of L-type (resp., RL-type) for $i=1,\dotsc,r$.
Therefore $\bG/\bH=\bF_{0}\vtr\bF_{1}\vtr\dotsb\vtr\bF_{r}=\{1\}$ is an LNC
(resp., an RLNC).

(''If" part.) Let $\bG/\bH=\bF_{0}\vtr\bF_{1}\vtr\dotsb\vtr\bF_{r}=\{1\}$ be
an LNC (resp., an RLNC) and $(0)=I_{0}\subset I_{1}\subset\dotsb\subset I_{r}$
the corresponding sequence of Hopf ideals of $k[\bG/\bH]$.
Each $I_{i}/I_{i-1}$ is the normal Hopf ideal of $k[\bF_{i-1}]
=k[\bG/\bH]/I_{i-1}$ corresponding to the Hopf subalgebra
$k[\bF_{i-1}/\bF_{i}]$, and thus $I_{i}/I_{i-1}=k[\bF_{i-1}]\cdot
k[\bF_{i-1}/\bF_{i}]^{+}$. If we put $I_{i}':=k[\bG]\cdot I_{i}$
($i=0,\dotsc,r$), then each $I_{i}'$ becomes a Hopf ideal of $k[\bG]$.
Let $\bG_{i}$ be the closed subgroup scheme of $\bG$ which corresponds to
$I_{i}'$. Since $k[\bG]$ is a faithfully flat $k[\bG/\bH]$-module, we have
an inclusion
\[ k[\bF_{i-1}]=k[\bG/\bH]/I_{i-1}\hookrightarrow k[\bG]\otimes_{k[\bG/\bH]}
(k[\bG/\bH]/I_{i-1})\simeq k[\bG]/I_{i-1}'=k[\bG_{i-1}] \]
for each $i=1,\dotsc,r+1$. Considering $I_{i}/I_{i-1}=k[\bF_{i-1}]\cdot
k[\bF_{i-1}/\bF_{i}]^{+}$ through this inclusion, we have that
$I_{i}'/I_{i-1}'=k[\bG_{i-1}]\cdot k[\bF_{i-1}/\bF_{i}]^{+}$.
Hence we have a normal chain
$\bG=\bG_{0}\vtr\bG_{1}\vtr\dotsb\vtr\bG_{r}=\bH$ such that
$\bG_{i-1}/\bG_{i}\simeq\bF_{i-1}/\bF_{i}$ ($i=1,\dotsc,r$).
Therefore $\bG$ is Liouville (resp., has an RLNC).

(3) We use induction on the least length $r$ of LNC. The case $r=0$ is
clear. Let $r>0$. Since $\bG$ is connected, $\bG/\bG_{1}$ is also connected.
Then $\bG/\bG_{1}$ is of RL-type and hence abelian. Therefore $\sD\bG$
(see \cite[(10.1)]{Waterhouse}) is a connected closed subgroup scheme of
$\bG_{1}$. By (1) and its proof, $\sD\bG$ is connected Liouville and has
an LNC with length $\leq r-1$. Then $\sD\bG$ is solvable by inductive
assumption, concluding the proof.
\end{proof}

The converse of (3) above does not hold in general:
\begin{example}
\label{notLiouville}
(1) A nontrivial anisotropic torus $\bT$ is connected solvable but not
Liouville since both $\Hom(\bT,\Gm)$ and $\Hom(\bT,\Ga)$ are trivial.

(2) Let $k$ be the prime field of $\ch(k)=2$ and $H=k[x]/\lan x^{4}+x^{2}+x
\ran$ with $x$ primitive. Then $H$ is a commutative Hopf algebra and
$\bG=\Spec H$ is abelian, finite etale, and unipotent. The Cartier dual
$\bG^{\ast}$ is finite connected of multiplicative type and hence solvable.
Since $H^{\ast}$ does not have any nontrivial grouplike, $\Hom(\bG^{\ast},\Gm)$
is trivial. Therefore $\bG^{\ast}$ is not Liouville.
\end{example}

\begin{proposition}
\label{RLNC}
Let $\bG$ be a connected algebraic affine group scheme over a field $k$.
Then $\bG$ is Liouville iff $\bG$ has an RLNC.
\end{proposition}
\begin{proof}
The ``if" part is clear. For the ``only if" part,
we use induction on the least length $r$ of LNC $\bG=\bG_{0}\vtr\bG_{1}\vtr
\dotsb\vtr\bG_{r}=\{1\}$. The case $r=0$ is clear. Let $r>0$ and assume
$\bG_{1}^{\circ}$ has an RLNC. By the argument in the proof of 
Lemma \ref{solvability1} (3), we have $\bG\vtr\bG_{1}^{\circ}$ and
$\bG/\bG_{1}^{\circ}$ is abelian. Thus the proof can be reduced to the case
that $\bG$ is connected abelian by Lemma \ref{solvability1} (2).

Let $\bG$ be connected abelian and put $H=k[\bG]$. Let $H_{\mathrm{u}}$
($=H^{1}$) be the irreducible component of $H$ which contains $1$ 
and $H_{\mathrm{s}}=H/HH_{\mathrm{u}}^{+}$. Then we have the exact sequence
\begin{equation}
\label{Jordan}
H_{\mathrm{u}}\rightarrowtail H\twoheadrightarrow H_{\mathrm{s}}.
\end{equation}
Let $\bar{k}$ denote the algebraic closure of $k$. It is known that
$H_{\mathrm{u}}\otimes_{k}\bar{k}$ is also the irreducible component of
$H\otimes_{k}\bar{k}$ containing $1$. (Let $\wedge_{k}$
(resp.\ $\wedge_{\bar{k}}$) denotes the wedge product of $k$-subspaces
(resp.\ $\bar{k}$-subspaces) of $H$ (resp.\ $H\otimes_{k}\bar{k}$)
in the sense of \cite[Ch.\ IX]{Sweedler1969}. Then
$(H\otimes_{k}\bar{k})^{1}=\bigcup_{n=0}^{\infty}\wedge^{n}_{\bar{k}}\bar{k}
=\bigcup_{n=0}^{\infty}(\wedge_{k}k)\otimes_{k}\bar{k}
=H_{\mathrm{u}}\otimes_{k}\bar{k}$.)
The exact sequence (\ref{Jordan}) splits over $\bar{k}$ (the Jordan
decomposition of $\bG_{\bar{k}}$ \cite[(9.5)]{Waterhouse}),
and $\Gs:=\Spec H_{\mathrm{s}}$ is connected of multiplicative type
since $(\Gs)_{\bar{k}}$ is connected diagonalizable.
Put $\Gu:=\Spec H_{\mathrm{u}}$ ($=\bG/\Gs$); this is unipotent.
We see $\Gu$ has an RLNC whose all factor group schemes
are of $\Ga$-type (see \cite[Ch.\ 16, Ex.\ 5]{Waterhouse}).
Then it suffices to show that $\Gs$ has an RLNC. 
Let $\bT$ be a maximal torus of $\Gs$.
$\bT$ includes no nontrivial anisotropic subtorus since it is Liouville.
Hence, by \cite[(7.4)]{Waterhouse}, we see $\bT$ is a split torus and has
an RLNC. Put $\bH=\Gs/\bT$; this is finite connected, Liouville,
and of multiplicative type. Let $\bH\vtr\bH_{1}\vtr\dotsb\vtr\bH_{r}=\{1\}$
be an LNC. We see $\bH/\bH_{1}$ is of $\Gm$-type. Since $\bH$ is finite
connected, $k[\bH]$ is a local algebra of finite dimension. Then its quotient
$k[\bH_{1}]$ is also a local algebra of finite dimension and hence $\bH_{1}$
is connected. By inductive assumption, $\bH_{1}$ has an RLNC. Therefore $\bH$
also has an RLNC, concluding the proof.
\end{proof}

\begin{proposition}
\label{solvability2}
Let $k$ be an algebraically closed field and $\bG$ an algebraic
affine group scheme over $k$. Then $\bG$ is Liouville iff
$\bG^{\circ}$ is solvable.
\end{proposition}
\begin{proof}
In fact we have proved the ``only if" part in Lemma \ref{solvability1} over
an arbitrary field. For the ``if" part, we use induction on the least $m$
such that $\sD^{m}\bG^{\circ}=\{1\}$. The case $m=0$ is clear. Let $m>0$ and
assume that $\sD\bG^{\circ}$ has an RLNC. By Lemma \ref{solvability1} (2), it
suffices to show that $\bG^{\circ}/\sD\bG^{\circ}$ has an RLNC. Thus
the proof can be reduced to the case that $\bG$ is (connected) abelian.

Let $\bG$ be abelian and take the Jordan decomposition $\bG=\Gs\times\Gu$.
$\Gu$ has an RLNC. Since $k$ is algebraically closed, $\Gs$ is diagonalizable
and hence has an RLNC. Therefore $\bG$ has an RLNC.
\end{proof}

As is seen in the following example, we see that the triangulability is
certainly stronger than the condition to have an RLNC, even if $k$ is
algebraically closed.
\begin{example}[{\cite[Ch.\ 10, Ex.\ 3]{Waterhouse}}]
Let $k$ be a field with $\ch(k)=2$ and $\bG$ the closed subgroup scheme
of $\mathbf{SL}_{2}$ over $k$ defined by the Hopf algebra
\[ k[\bG]=k[X_{11},X_{12},X_{21},X_{22}]/\lan X_{11}X_{22}+X_{12}X_{21}+1,
X_{11}^{2}+1,X_{22}^{2}+1,X_{12}^{2},X_{21}^{2}\ran. \]
Put $I=\lan X_{11}X_{22}+X_{12}X_{21}+1,X_{11}^{2}+1,X_{22}^{2}+1,X_{12}^{2},
X_{21}^{2}\ran\subset k[X_{11},X_{12},X_{21},X_{22}]$.
Since $\sqrt{I}=\lan X_{11}+1,X_{22}+1,X_{12},X_{21}\ran$, we have
$k[X_{11},X_{12},X_{21},X_{22}]/\sqrt{I}=k$ and hence $\bG$ is finite
connected. Let $H_{1}=k[\bG]/k[\bG]X_{21}$ and $\bG_{1}=\Spec H_{1}$.
Then we have exact sequences
\[ \begin{array}{rccll}
k[\boldsymbol{\alpha}_{2}]=k[X]/\lan X^{2}\ran & \rightarrowtail & k[\bG] &
\twoheadrightarrow & H_{1} \\
X & \mapsto & X_{21}X_{22} & &
\end{array} \]
and
\[ \begin{array}{rccll}
k[\boldsymbol{\alpha}_{2}] & \rightarrowtail & H_{1} & \twoheadrightarrow &
H_{1}/H_{1}X_{12}\simeq k[\mu_{2}]. \\
X & \mapsto & X_{11}X_{12} & & 
\end{array} \]
Thus $\bG\vtr\bG_{1}\vtr\mu_{2}\vtr\{1\}$ is an RLNC. Since 
\[ kX_{11}\oplus kX_{12}\oplus kX_{21}\oplus kX_{22}\subset k[\bG] \]
is a $4$-dimensional simple subcoalgebra, $k[\bG]$ is not pointed and
hence $\bG$ is not triangulable.
\end{example}

It is known that $\bG$ is unipotent iff $\bG$ has an RLNC whose all factor
group schemes are of $\Ga$-type. We say that $\bG$ is {\em $\Gm$-composite}
iff $\bG$ has an RLNC whose all factor group schemes are of $\Gm$-type.
When $k$ is algebraically closed and $\bG$ comes from the affine algebraic
group $\bG(k)$ (in the sense of \cite[(4.5)]{Waterhouse}), $\bG$ is
$\Gm$-composite iff $\bG(k)$ is solvable and ``quasicompact" in Kolchin's
terminology, which implies that each element of $\bG(k)$ is
semisimple \cite[\S 6, Theorem 2]{Kolchin1948}.
In general it is difficult to characterize the condition to be $\Gm$-composite.
As is seen above, not all group schemes of multiplicative type are
$\Gm$-composite. On the other hand, non-diagonalizable group schemes can be
$\Gm$-composite; see the following example.
\begin{example}
(1) Let $k$ be the prime field with $\ch(k)=p>0$ and take the commutative
Hopf algebra $H=k[x,y]/\lan x^{p}-x,y^{p}-x-y\ran$ with $x,y$ primitive.
One sees $\bG=\Spec H$ is abelian, finite etale, and unipotent. Hence the
Cartier dual $\bG^{\ast}$ is of multiplicative type and connected.
We have the RLNC of $\bG$:
\[ k[x]/\lan x^{p}-x\ran \rightarrowtail H \twoheadrightarrow
k[y]/\lan y^{p}-y\ran. \]
By dualizing this we see that $\bG^{\ast}$ is $\Gm$-composite:
\[ k[\mu_{p}] \twoheadleftarrow H^{\ast} \leftarrowtail
k[\mu_{p}]. \]
The grouplikes of $H^{\ast}$ is given by
\[ \Coalg_{k}(k,H^{\ast})\simeq\Alg_{k}(H,k)=\{(a,b)\in k^{2}\,|\,
a^{p}-a=0,\quad b^{p}-a-b=0\}. \]
Thus we have $\bG^{\ast}$ is not diagonalizable since $p^{2}=\dim_{k}H^{\ast}
\neq p =\#\Alg_{k}(H,k)$.

(2) \cite[\S 2, Remark 2]{Kolchin1948b}. Let $k$ be an infinite field of
arbitrary characteristic. Let $\bG$ be the closed subgroup scheme of $\bGL_{2}$
over $k$ defined by the Hopf algebra
\[ k[\bG]=k[X_{11},X_{12},X_{21},X_{22},1/(X_{11}X_{22}-X_{12}X_{21})]/
\lan X_{11}X_{12},X_{12}X_{22},X_{11}X_{21},X_{21}X_{22}\ran. \]
Since $\bG(k)$ can not be simultaneously diagonalized, $\bG$ is not
diagonalizable (see \cite[(4.6)]{Waterhouse}).
We see $\pi_{0}\bG=\bbz/2\bbz$. Indeed, $\pi_{0}(k[\bG])=ke_{0}\oplus ke_{1}$
where
\[ e_{0}=\frac{X_{11}X_{22}}{X_{11}X_{22}-X_{12}X_{21}},\quad
e_{1}=-\frac{X_{12}X_{21}}{X_{11}X_{22}-X_{12}X_{21}}. \]
Let $H_{1}=k[\bG]/k[\bG]X_{12}X_{21}$, $H_{2}=H_{1}/H_{1}(X_{11}X_{22}-1)$,
and $\bG_{i}=\Spec H_{i}$ ($i=1,2$). Then we have exact
sequences
\[ k(\bbz/2\bbz)^{\ast} \rightarrowtail k[\bG] \twoheadrightarrow H_{1} \]
and
\[ \begin{array}{rccll}
k[\Gm]=k[X,X^{-1}] & \rightarrowtail & H_{1} & \twoheadrightarrow &
H_{2}\simeq k[\Gm]. \\
X & \mapsto & X_{11}X_{22} & & 
\end{array} \]
Thus $\bG\vtr\bG_{1}\vtr\bG_{2}\vtr\{1\}$ is an LNC. If $\ch(k)\neq 2$, then
$\bbz/2\bbz\simeq\mu_{2}$ and hence $\bG$ is $\Gm$-composite. 
\end{example}

\subsection{Finite etale extensions}

In what follows, we always assume that $L/K$ is an extension of AS $D$-module
algebras such that $L^{D}=K^{D}=:k$.

\begin{definition}
\label{def-finite-etale-extensions}
We say that $L/K$ is a {\em finite etale extension} iff $L$ is a separable
$K$-algebra in the sense of \cite{DeMeyer-Ingraham}, i.e.\ $L$ is a projective
$L\otimes_{K}L$-module.
\end{definition}

For a commutative $K$-algebra $A$, let $\pi_{0}(A)$ denote the union of
all separable $K$-subalgebras.
If we take a maximal ideal $P$ of $L$ and put $L'=L/P$ and $K'=K/P\cap K$,
then the following are equivalent:
\begin{itemize}
\item $L/K$ is a finite etale extension.
\item $L$ is a finitely generated $K$-algebra and $\pi_{0}(L)=L$.
\item $L'/K'$ is a finite separable field extension.
\end{itemize}
Take a maximal ideal $\gp$ of $K$.
We say a finite etale extension $L/K$ is {\em copied} (resp., {\em anticopied})
iff $L'=K'$ (resp., $\gp L$ is a maximal ideal of $L$); this condition is
independent of the choice of $P$ (resp., $\gp$).

\begin{lemma}
\label{etale}
Let $(L/K,A,H)$ be a finitely generated PV extension.

{\rm (1)} $(\pi_{0}(A)/K,\pi_{0}(A),\pi_{0}(H))$ is also a PV extension and
hence $\pi_{0}(A)$ is the intermediate AS $D$-module algebra
which corresponds to $\bG(L/K)^{\circ}$.

{\rm (2)} $\pi_{0}(A)=\pi_{0}(L)$.

Especially $L/K$ is a finite etale extension iff $\bG(L/K)$ is finite etale.
\end{lemma}
\begin{proof}
(1) The $\mu$-isomorphism $A\otimes_{k}H\xrightarrow{\sim}A\otimes_{K}A$
restricts to an algebra isomorphism $\pi_{0}(A)\otimes_{k}\pi_{0}(H)
\xrightarrow{\sim}\pi_{0}(A)\otimes_{K}\pi_{0}(A)$.
Hence $\theta^{-1}(A\otimes_{k}\pi_{0}(H))=\pi_{0}(A)$;
this implies that $(\pi_{0}(A)/K,\pi_{0}(A),\pi_{0}(H))$ is a PV extension
by Proposition \ref{prop-Galois-correspondence3}.

(2) By Proposition \ref{prop-interlacing2}, we may assume that $K$ is a field
and $A$ is an integral domain. Put $M=\pi_{0}(A)$. It suffices to show that $M$
is separably closed in $L$. Write $H_{0}=H/H\pi_{0}(H)^{+}$. Then we have
an algebra isomorphism $L\otimes_{k}H_{0}\xrightarrow{\sim}L\otimes_{M}A$
since $(L/M,A,H_{0})$ is a PV extension by Proposition
\ref{prop-Galois-correspondence2}.
By \cite[(6.6)]{Waterhouse}, $(L\otimes_{k}H_{0})/\mathrm{nilradical}$ is
an integral domain. It follows that $(L\otimes_{M}L)/\mathrm{nilradical}$ is
also an integral domain since it is a localization of
$(L\otimes_{M}A)/\mathrm{nilradical}$ by a multiplicative subset not containing
zero. Thus $L\otimes_{M}L$ has no nontrivial idempotent.
This implies that $M$ is separably closed in $L$. Indeed, let $S$ be an
intermediate field of $L/M$ such that $S/M$ is finite separable.
If $e\in S\otimes_{M}S$ ($\subset L\otimes_{M}L$) is a separability
idempotent for $S$ (see \cite[p.\ 40]{DeMeyer-Ingraham}),
then necessarily $e=1\otimes 1$.
Let $J:=\Ker(\mathrm{mult} : S\otimes_{M}S\rightarrow S)$. Then $0=Je=J$.
Hence $a\otimes 1-1\otimes a=0$ in $S\otimes_{M}S$ for all $a\in S$.
This implies $S=M$.
\end{proof}

\begin{corollary}
Let $L/K$ be a finitely generated PV extension.
Then $\bG(L/K)$ is connected iff $\pi_{0}(L)=K$.
\end{corollary}

Let $\sG$ be a finite group and put $H=(k\sG)^{\ast}$. Consider $D_{k}
=D\otimes_{R}k$ as a right $H$-comodule algebra with the trivial structure
map $d\mapsto d\otimes 1$.
We say $L/K$ is a {\em $\sG$-extension} iff
\begin{enumerate}
\renewcommand{\labelenumi}{(\roman{enumi})}
\item $L$ is an algebra in the symmetric tensor category
$({}_{D_{k}}\mathcal{M}^{H},\otimes_{k},k)$, and
\item $L/K$ is a right $H$-Galois extension.
\end{enumerate}
Here ${}_{D_{k}}\mathcal{M}^{H}$ denotes the category of right
$(H,D_{k}^{\mathrm{op}})$-Hopf modules as in the proof of Theorem
\ref{equivalence-to-representations}.
We easily see that $L/K$ is a $\sG$-extension iff
$(L/K,L,(k\sG)^{\ast})$ is a PV extension. If $L/K$ is an anticopied
$\sG$-extension, then $L'/K'$ is a Galois extension
of fields in ordinary sense such that $\Gal(L'/K')=\sG$. Conversely,
when $L/K$ is a finite Galois extension of fields, $L/K$ is
a $\Gal(L/K)$-extension iff every element of $\Gal(L/K)$ is $D$-linear.

\subsection{$\Ga$-primitive extensions and $\Gm$-primitive extensions}

\begin{definition}
\label{def-primitive-extensions}
(1) An $x\in L$ is called {\em $\Ga$-primitive} over $K$ iff $d(x)\in K$
for all $d\in D^{+}$. In this case, we say that
$K\lan x\ran/K$ is a {\em $\Ga$-primitive extension}.

(2) An $x\in L$ is called {\em $\Gm$-primitive} over $K$ iff
$x$ is a non-zero divisor of $L$ and $d(x)x^{-1}\in K$ for all $d\in D$.
In this case, we say that $K\lan x\ran/K$ is a {\em $\Gm$-primitive extension}.
\end{definition}

As in \cite[(2.5a), (2.5b)]{Takeuchi1989}, we have the following lemmas:
\begin{lemma}
\label{primitive1}
{\rm (1)} Let $K\lan x\ran/K$ be a $\Ga$-primitive extension. Put $A=K[x]$
and $l=1\otimes_{K}x-x\otimes_{K}1\in (A\otimes_{K}A)^{D}$. Then
$(K\lan x\ran/K,A,k[l])$ is a PV extension with $l$ primitive and the
PV group scheme $\bG(K\lan x\ran/K)$ of $\Ga$-type.

{\rm (2)} Let $K\lan x\ran/K$ be a $\Gm$-primitive extension.
Put $A=K[x,x^{-1}]$ and $g=x^{-1}\otimes_{K}x\in (A\otimes_{K}A)^{D}$. Then
$(K\lan x\ran/K,A,k[g,g^{-1}])$ is a PV extension with $g$ grouplike
and the PV group scheme $\bG(K\lan x\ran/K)$ of $\Gm$-type.
\end{lemma}
\begin{proof}
(1) $x$ is $\Ga$-primitive iff there exists a $\varphi\in\Hom_{R}(D,K)$
such that $d(x)=\varepsilon(d)x+\varphi(d)$ for all $d\in D$. Then
\[ X=\left(\begin{array}{cc} 1 & 1 \\ x & x+1 \end{array}\right)\in
GL_{2}(K\lan x\ran) \]
is $GL_{2}$-primitive over $K$. In fact,
\[ dX=\left(\begin{array}{cc} \varepsilon(d) & 0 \\ \varphi(d) &
\varepsilon(d) \end{array}\right)X\qquad (d\in D). \]
Recalling the proof of Theorem \ref{thm-characterization} 
(d) $\Rightarrow$ (a), we see
\[ Z=(X^{-1}\otimes_{K}1)(1\otimes_{K}X)=\left(\begin{array}{cc} 1-l & -l \\
l & 1+l \end{array}\right), \]
which concludes the proof.

(2) This is equivalent to saying that $x$ is $GL_{1}$-primitive over $K$.
\end{proof}

\begin{lemma}
\label{primitive2}
{\rm (1)} If $l\in (L\otimes_{K}L)^{D}$ and if $l$ is primitive
in the $L$-coring $L\otimes_{K}L$, then there exists
an $x\in L$ such that $l=1\otimes_{K}x-x\otimes_{K}1$ and $x$ is
$\Ga$-primitive over $K$.

{\rm (2)} If $g\in (L\otimes_{K}L)^{D}$ and if $g$ is grouplike
in $L\otimes_{K}L$, then there exists a non-zero divisor $x\in L$
such that $g=x^{-1}\otimes_{K}x$ and $x$ is $\Gm$-primitive over $K$.
\end{lemma}
\begin{proof}
(1) Primitive elements in the $L$-coring $L\otimes_{K}L$ are precisely
$1$-cocycles in the Amitsur complex:
\[ \begin{array}{l}
0\rightarrow L\xrightarrow{\delta_{0}}L\otimes_{K}L\xrightarrow{\delta_{1}}
L\otimes_{K}L\otimes_{K}L\xrightarrow{\delta_{2}}\dotsb, \vspace{1ex} \\
\delta_{0}(x)=1\otimes_{K}x-x\otimes_{K}1, \vspace{1ex} \\
\delta_{1}({\displaystyle\sum}x_{i}\otimes_{K}y_{i})=
{\displaystyle\sum}1\otimes_{K}x_{i}\otimes_{K}y_{i}
{}-{\displaystyle\sum}x_{i}\otimes_{K}1\otimes_{K}y_{i}
{}+{\displaystyle\sum}x_{i}\otimes_{K}y_{i}\otimes_{K}1,\quad \dotsc,
\end{array} \]
whose $n$-th cohomology is $H^{n}(L/K,\Ga)$.
But $H^{1}(L/K,\Ga)=0$ since $L/K$ is a faithfully flat extension
(see \cite[Ch.\ 17, Ex.\ 10]{Waterhouse}).
Then $l\in\Ker\delta_{1}=\im\delta_{0}$ and hence there exists some
$x\in L$ such that $l=1\otimes_{K}x-x\otimes_{K}1$. Since $dl=\varepsilon(d)l$
for all $d\in D$, we have $(dx)\otimes_{K}1=1\otimes_{K}(dx)$ for all $d\in
D^{+}$. This implies $dx\in K$ for all $d\in D^{+}$.

(2) Grouplike elements in $L\otimes_{K}L$ are precisely $1$-cocycles in the
complex:
\[ \begin{array}{l}
\{1\}\rightarrow\Gm(L)\xrightarrow{\delta_{0}}\Gm(L\otimes_{K}L)
\xrightarrow{\delta_{1}}\Gm(L\otimes_{K}L\otimes_{K}L)\xrightarrow{\delta_{2}}
\dotsb, \vspace{1ex} \\
\delta_{0}(x)=(1\otimes_{K}x)(x\otimes_{K}1)^{-1}=x^{-1}\otimes_{K}x,
\vspace{1ex} \\
\delta_{1}({\displaystyle\sum}x_{i}\otimes_{K}y_{i})=
({\displaystyle\sum}1\otimes_{K}x_{i}\otimes_{K}y_{i})
({\displaystyle\sum}x_{i}\otimes_{K}1\otimes_{K}y_{i})^{-1}
({\displaystyle\sum}x_{i}\otimes_{K}y_{i}\otimes_{K}1),\quad \dotsc,
\end{array} \]
whose $n$-th cohomology is $H^{n}(L/K,\Gm)$.
But $H^{1}(L/K,\Gm)=\Pic(L/K)\subset\Pic(K)=\{1\}$ since $K$ is a finite
product of fields. Then $g\in\Ker\delta_{1}=\im\delta_{0}$ and hence
there exists some $x\in\Gm(L)$ such that $g=x^{-1}\otimes_{K}x$.
Since $dg=\varepsilon(d)g$, we have
\[ 1\otimes_{K}dx=d(1\otimes_{K}x)=d((x\otimes_{K}1)g)=d(x)x^{-1}\otimes_{K}x
\]
for all $d\in D$. By multiplying $1\otimes_{K}x^{-1}$, we have
$1\otimes_{K}d(x)x^{-1}=d(x)x^{-1}\otimes_{K}1$ for all $d\in D$, which 
implies $d(x)x^{-1}\in K$ for all $d\in D$.
\end{proof}

\begin{proposition}
\label{primitive3}
$L/K$ is a $\Ga$-primitive (resp., $\Gm$-primitive) extension iff
$L/K$ is a PV extension and $\bG(L/K)$ is of $\Ga$-type (resp., $\Gm$-type).
\end{proposition}
\begin{proof}
(``Only if" part.) This has been proved in Lemma \ref{primitive1}.

(``If" part.) Let $k[l]$ (resp., $k[g,g^{-1}]$) be the Hopf algebra for $L/K$.
By Lemma \ref{primitive2}, there exists the corresponding $x\in L$. 
Then $K\lan x\ran$ is an intermediate AS $D$-module algebra of $L/K$ such
that $K\lan x\ran/K$ is a PV extension. Since the Hopf algebras of
$K\lan x\ran/K$ and $L/K$ coincide, we have $L=K\lan x\ran$.
\end{proof}

\subsection{The solvability theorem}

\begin{definition}
\label{defLE}
Let $F/K$ be a finitely generated extension of AS $D$-module algebras.
We call $F/K$ a {\em Liouville extension} iff $F^{D}=K^{D}=k$ and
there exists a sequence of AS $D$-module algebras
\begin{equation}
\label{LC}
K=F_{0}\subset F_{1}\subset\dotsb\subset F_{r}=F
\end{equation}
such that each $F_{i}/F_{i-1}$ ($i=1,\dotsc,r$) is at least one of the
following types: $\Ga$-primitive extension, $\Gm$-primitive extension,
or finite etale extension. In this case, the sequence (\ref{LC}) is called
a {\em Liouville chain}. 
Moreover, $F/K$ is called a {\em Liouville extension of type ($j$)}
($j=1,\dotsc,10$) iff $F/K$ has a Liouville chain (\ref{LC}) such that
each extension $F_{i}/F_{i-1}$ ($i=1,\dotsc,r$) is
\begin{enumerate}
\renewcommand{\labelenumi}{\rm(\arabic{enumi})}
\item $\Ga$-primitive, $\Gm$-primitive, or finite etale,
\item $\Ga$-primitive or $\Gm$-primitive,
\item $\Gm$-primitive or finite etale,
\item $\Ga$-primitive or finite etale,
\item $\Ga$-primitive or a $\sG$-extension for a finite solvable group $\sG$,
\item $\Gm$-primitive,
\item $\Ga$-primitive,
\item finite etale,
\item a $\sG$-extension for a finite solvable group $\sG$,
\item trivial (i.e.\ $F_{i}=F_{i-1}$),
\end{enumerate}
respectively. Here we are taking priority of the compatibility with
\cite[\S 24]{Kolchin1948}.
We observe an anticopied $\sG$-extension for
a finite solvable group $\sG$ is identified with a Galois extension by
radicals and is also a Liouville extension of type (6).
\end{definition}

To show the solvability theorem, we need the following lemma
(cf.\ \cite[\S 21]{Kolchin1948}).
\begin{lemma}
\label{newelement}
Let $L/K$ be a finitely generated PV extension and $F$ an AS $D$-module algebra
including $L$ such that $F^{D}=K^{D}=k$. Take one $t\in F$.
Then $L\lan t\ran/K\lan t\ran$ is a finitely generated PV extension and
$\bG(L\lan t\ran/K\lan t\ran)\simeq\bG(L/K\lan t\ran\cap L)$.
\end{lemma}
\begin{proof}
By Theorem \ref{thm-characterization},
there exists a $GL_{n}$-primitive $X=(x_{ij})
\in GL_{n}(L)$ over $K$ such that $L=K\lan x_{ij}\ran$. Since $L\lan t\ran
=K\lan t,x_{ij}\ran$, we have that $L\lan t\ran/K\lan t\ran$ is a finitely
generated PV extension. Write $M=K\lan t\ran\cap L$, $Z=(X^{-1}\otimes_{M}1)
(1\otimes_{M}X)=(z_{ij})$, and $Z^{-1}=(w_{ij})$. Then $H=k[z_{ij},w_{ij}]$
becomes the Hopf algebra for $L/M$. Similarly by writing
$Z'=(X^{-1}\otimes_{K\lan t\ran}1)(1\otimes_{K\lan t\ran}X)=(z_{ij}')$,
and $(Z')^{-1}=(w_{ij}')$, we obtain the Hopf algebra $H'=k[z_{ij}',w_{ij}']$
for $L\lan t\ran/K\lan t\ran$. It follows that there exists a surjective 
Hopf algebra map $\varphi : H\twoheadrightarrow H'$, $z_{ij}\mapsto z_{ij}'$.
This implies that $\bG(L\lan t\ran/K\lan t\ran)$ is a closed subgroup
scheme of $\bG(L/M)$. Let $I=\Ker\varphi$ be the corresponding Hopf ideal.

$\varphi$ is the restriction (to $H$) of the natural map $\tilde{\varphi} :
L\otimes_{M}L\rightarrow L\lan t\ran\otimes_{K\lan t\ran}L\lan t\ran$.
Since the coideal $I\cdot(L\otimes_{M}L)$ of $L\otimes_{M}L$, which corresponds
to $I$, is included in $\Ker\tilde{\varphi}$,
we have $\{a\in L\,|\,a\otimes_{M}1-1\otimes_{M}a\in I\cdot(L\otimes_{M}L)\}
\subset\{a\in L\,|\,a\otimes_{M}1-1\otimes_{M}a\in\Ker\tilde{\varphi}\}
= L\cap K\lan t\ran =  M$.
This implies that the intermediate AS $D$-module algebra of $L/M$ which
corresponds to $I$ is $M$. Thus $I=0$.
\end{proof}

\begin{theorem}
\label{mainTh}
Let $L/K$ be a finitely generated PV extension.
Then the following are equivalent:
\begin{enumerate}
\renewcommand{\labelenumi}{\rm(\alph{enumi})}
\item $L/K$ is a Liouville extension.
\item There exists a Liouville extension $F/K$ such that $L\subset F$.
\item $\bG(L/K)$ is Liouville.
\end{enumerate}
When $k$ is algebraically closed, these are equivalent to:
\begin{enumerate}
\renewcommand{\labelenumi}{\rm(d)}
\item $\bG(L/K)^{\circ}$ is solvable.
\end{enumerate}
\end{theorem}
\begin{proof}
((a) $\Rightarrow$ (b)) This is clear.

((b) $\Rightarrow$ (c)) Take a Liouville chain of $F/K$:
\[ K=F_{0}\subset F_{1}\subset\dotsb\subset F_{r}=F. \]
We use induction on $r$. The case $r=0$ is obvious. Let $r>0$. 
Since there are finite $t_{1},\dotsc,t_{s}\in F$ such that
$F_{1}=K\lan t_{1},\dotsc,t_{s}\ran$, we have that $L\lan t_{1},\dotsc,t_{s}
\ran/F_{1}$ is a finitely generated PV extension and
$\bG(L\lan t_{1},\dotsc,t_{s}\ran/F_{1})\simeq\bG(L/F_{1}\cap L)$ by
Lemma \ref{newelement}. By the inductive assumption, $\bG(L/F_{1}\cap L)$
is Liouville.

If $F_{1}/K$ is a finite etale extension, then $F_{1}\cap L\subset\pi_{0}(L)$.
Hence we have $\bG(L/F_{1}\cap L)\supset\bG(L/\pi_{0}(L))=\bG(L/K)^{\circ}$
(Lemma \ref{etale}).
Thus, $\bG(L/F_{1}\cap L)^{\circ}=\bG(L/K)^{\circ}$ and hence (c) holds
by Lemma \ref{solvability1} (1).

If $F_{1}/K$ is a $\Ga$-primitive extension, then there exists a
$\Ga$-primitive $x\in F_{1}$ such that $F_{1}=K\lan x\ran$.
Write $L_{1}=F_{1}\cap L$. One sees that $L_{1}/K$ is also a $\Ga$-primitive
extension (see \cite[(2.9a)]{Takeuchi1989}).
Hence $\bG(L/K)\vtr\bG(L/L_{1})$ and $\bG(L/K)/\bG(L/L_{1})=\bG(L_{1}/K)$
is of $\Ga$-type. Therefore (c) holds.

If $F_{1}/K$ is a $\Gm$-primitive extension, then there exists a
$\Gm$-primitive $x\in F_{1}$ such that $F_{1}=K\lan x\ran$.
Write $L_{1}=F_{1}\cap L$. One sees that $L_{1}/K$ is also a $\Gm$-primitive
extension (see \cite[(2.9b)]{Takeuchi1989}).
Then we obtain (c) in the same way to the above.

((c) $\Rightarrow$ (a)) Let $\bG(L/K)=\bG_{0}\vtr\bG_{1}\vtr\dotsb\vtr\bG_{r}
=\{1\}$ be an LNC and $L_{i}$ ($i=0,\dotsc,r$) the intermediate AS
$D$-module algebra which corresponds to $\bG_{i}$. Then by Lemma \ref{etale}
and by Proposition \ref{primitive3}, $K=L_{0}\subset L_{1}\subset\dotsb
\subset L_{r}=L$ is a Liouville chain.
\end{proof}

By Proposition \ref{RLNC}, we have the following.

\begin{corollary}
Let $L/K$ be a finitely generated PV extension.
If $L/K$ is a Liouville extension, then there exists a Liouville chain
\[ K=L_{0}\subset\pi_{0}(L)=L_{1}\subset L_{2}\subset\dotsc\subset L_{r}=L \]
such that each $L_{i}/L_{i-1}$ ($i=2,\dotsc,r$) is $\Gm$-primitive or
$\Ga$-primitive extension.
\end{corollary}

\begin{corollary}
\label{precisely}
Let $L/K$ be a finitely generated PV extension. Then $L/K$ is (included in)
a Liouville extension of type ($j$) ($j=1,\dotsc,10$) iff
\begin{enumerate}
\renewcommand{\labelenumi}{\rm(\arabic{enumi})}
\item $\bG(L/K)$ is Liouville,
\item $\bG(L/K)$ has an RLNC,
\item $\bG(L/K)^{\circ}$ is $\Gm$-composite,
\item $\bG(L/K)^{\circ}$ is unipotent,
\item $\pi_{0}\bG(L/K)$ is finite constant and solvable,
and $\bG(L/K)^{\circ}$ is unipotent,
\item $\bG(L/K)$ is $\Gm$-composite,
\item $\bG(L/K)$ is unipotent,
\item $\bG(L/K)$ is finite etale,
\item $\bG(L/K)$ is finite constant and solvable,
\item $\bG(L/K)$ is trivial,
\end{enumerate}
respectively.
\end{corollary}

This corollary can become more explicit when $K$ is a perfect field and $k$ is
algebraically closed.
In such a case, if $(L/K,A,H)$ is a finitely generated PV extension,
then $A\otimes_{K}A$ is reduced (see \cite[Ch.\ 6, Ex.\ 2]{Waterhouse}), and so
$H$ is reduced. Thus $\bG(L/K)$ corresponds to the affine algebraic
group $\bG(L/K)(k)=\AutDlinKalg(L)$ in the sense of \cite[(4.5)]{Waterhouse}.
There we can do the following replacement on condition about $\bG(L/K)$:
\begin{enumerate}
\renewcommand{\labelenumi}{\rm(\arabic{enumi})}
\item $\bG(L/K)^{\circ}$ is solvable ($\Leftrightarrow$ $\bG(L/K)^{\circ}$ is
triangulable),
\item $\bG(L/K)$ is solvable,
\item $\bG(L/K)^{\circ}$ is diagonalizable,
\item $\bG(L/K)^{\circ}$ is unipotent,
\item $\bG(L/K)$ is solvable and $\bG(L/K)^{\circ}$ is unipotent,
\item $\bG(L/K)(k)$ is solvable and quasicompact (in Kolchin's sense),
\item $\bG(L/K)$ is unipotent,
\item $\bG(L/K)$ is finite constant,
\item $\bG(L/K)$ is finite constant and solvable,
\item $\bG(L/K)$ is trivial.
\end{enumerate}


\bibliographystyle{amsalpha}

\begin{thebibliography}{99}

\bibitem{Intro-vanderPut-Singer} M.\ van der Put, M.F.\ Singer,
\textit{``Galois theory of difference equations"},
Lecture Notes in Math.\ 1666, Springer, 1997.

\bibitem{Intro-Takeuchi1989} M.\ Takeuchi,
\textit{A Hopf algebraic approach to the Picard-Vessiot theory},
J.\ Algebra 122 (1989), 481--509.

\end{thebibliography}

\begin{thebibliography}{99}

\bibitem{AFK} K.\ Amano, M.\ Fujigami, T.\ Kogiso,
\textit{Construction of irreducible relative invariant of the prehomogeneous
vector space $(SL_{5}\times GL_{4},\Lambda^{2}(\bbc^{5})\otimes\bbc^{4})$},
Linear Algebra Appl.\ 355 (2002), 215--222.

\bibitem{Gyoja1} A.\ Gyoja,
\textit{Construction of invariants},
Tsukuba J.\ Math.\ 14 (1990), 437--457.   

\bibitem{Kable} A.C.\ Kable,
\textit{The concomitants of a prehomogeneous vector space},
J.\ Algebra 271 (2004), 295--311.

\bibitem{Kable-Yukie1} A.C.\ Kable, A.\ Yukie,
\textit{On the space of quadruples of quinary alternating forms},
J.\ Pure Appl.\ Algebra 186 (2004), 277--295.

\bibitem{Kable-Yukie2} A.C.\ Kable, A.\ Yukie,
\textit{A construction of quintic rings},
Nagoya Math.\ J.\ 173 (2004), 163--203.

\bibitem{Ochiai} H.\ Ochiai,
\textit{Quotients of some prehomogeneous vector spaces},
J.\ Algebra 192 (1997), 61--73.

\bibitem{Sato-Kimura} M.\ Sato, T.\ Kimura,
\textit{A classification of irreducible prehomogeneous vector spaces and their
relative invariants},
Nagoya Math.\ J.\ 65 (1977), 1--155.

\bibitem{Wright-Yukie} D.J.\ Wright, A.\ Yukie,
\textit{Prehomogeneous vector spaces and field extensions},
Invent.\ math.\ 110 (1992), 283--314.

\end{thebibliography}

\begin{thebibliography}{99}


\bibitem{Artin} E.\ Artin, translated by M.\ Butler,
\textit{``The gamma function"},
Holt, Rinehart and Winston, 1964.

\bibitem{Bernshtein} I.N.\ Bernshtein,
\textit{The analytic continuation of generalized functions with respect
to a parameter},
Funct.\ Anal.\ Appl.\ 6 (1972), 273--285.

\bibitem{Fujigami} M.\ Fujigami,
\textit{A generalization of the Bohr-Mollerup theorem and the $\Gamma$-factors
of local functional equations} (Japanese),
S\=urikaisekikenky\=usho K\=oky\=uroku, Kyoto Univ., 1238 (2001),
20--29.

\bibitem{Gelfand-Graev-Retakh} I.M.\ Gel'fand, M.I.\ Graev, V.S.\ Retakh,
\textit{General hypergeometric systems of equations and series of
hypergeometric type},
Russian Math.\ Surveys 47:4 (1992), 1--88.

\bibitem{Gyoja} A.\ Gyoja,
\textit{Theory of prehomogeneous vector spaces without regularity condition},
Publ.\ Res.\ Inst.\ Math.\ Sci., Kyoto Univ., 27 (1991), 861--922.

\bibitem{Igusa1986} J.\ Igusa,
\textit{On functional equations of complex powers},
Invent.\ math.\ 85 (1986), 1--29

\bibitem{Igusa2000} J.\ Igusa,
\textit{``An introduction to the theory of local zeta functions"},
Studies in Advanced Mathematics 14, AMS/IP, 2000.

\bibitem{Kashiwara} M.\ Kashiwara,
\textit{$B$-functions and holonomic systems---rationality of roots of
$b$-functions},
Invent.\ math.\ 38 (1976), 33--53.

\bibitem{Kimura} T.\ Kimura, translated by M.\ Nagura and T.\ Niitani,
\textit{``Introduction to prehomogeneous vector spaces"},
Translations of Mathematical Monographs 215, AMS, 2003.

\bibitem{Loeser-Sabbah} F.\ Loeser, C.\ Sabbah,
\textit{Equations aux diff\'erences finies et d\'eterminants d'int\'egrales
de fonctions multiformes},
Comment.\ Math.\ Helvetici 66 (1991), 458--503.

\bibitem{Mostow} G.D.\ Mostow,
\textit{Self-adjoint groups},
Ann.\ of Math.\ 62 (1955), 44--55.

\bibitem{Ore} O.\ Ore,
\textit{Sur la forme des fonctions hyperg\'eom\'etriques de plusieurs
variables},
J.\ Math.\ Pures et Appl.\ 9 (1930), 311--326.


\bibitem{FSato1982} F.\ Sato,
\textit{Zeta functions in several variables associated with prehomogeneous
vector spaces I: Functional equations},
T\^ohoku Math.\ J.\ (2) 34 (1982), 437--483.


\bibitem{Sato-Shintani-Muro} M.\ Sato, note by T.\ Shintani,
translated by M.\ Muro,
\textit{Theory of prehomogeneous vector spaces (algebraic part)---the English
translation of Sato's lecture from Shintani's note},
Nagoya Math.\ J.\ 120 (1990), 1--34.

\end{thebibliography}

\begin{thebibliography}{99}

\bibitem{Liouville} K.\ Amano,
\textit{Liouville extensions of artinian simple module algebras},
Comm.\ Algebra 34 (2006), 1811--1823.

\bibitem{Amano-Masuoka} K.\ Amano, A.\ Masuoka,
\textit{Picard-Vessiot extensions of artinian simple module algebras},
J.\ Algebra 285 (2005), 743--767.

\bibitem{Andre} Y.\ Andr\'e,
\textit{Diff\'erentielles non commutatives et th\'eorie de Galois
diff\'erentielle ou aux diff\'erences},
Ann.\ Sci.\ \'Ecole Norm.\ Sup.\ (4) 34 (2001), 685--739.

\bibitem{Bialynicki-Birula} A.\ Bia{\l}ynicki-Birula,
\textit{On Galois theory of fields with operators},
Amer.\ J.\ Math.\ 84 (1962), 89--109.

\bibitem{Deligne} P.\ Deligne,
\textit{Cat\'egories tannakiennes},
In P.\ Cartier et al.\ (editors), ``Grothendieck Festschrift", Vol.\ 2,
Progress in Math.\ 87, Birkh\"auser, 1990, pp.\ 111--195.

\bibitem{Deligne-Milne} P.\ Deligne, J.S.\ Milne,
\textit{Tannakian categories},
In P.\ Deligne et al.\ (editors),
``Hodge cycles, motives and Shimura varieties",
Lecture Notes in Math.\ 900, Springer, 1982, pp.\ 101--228.

\bibitem{DeMeyer-Ingraham} F.\ DeMeyer, E.\ Ingraham,
\textit{``Separable algebras over commutative rings"},
Lecture Notes in Math.\ 181, Springer, 1971.

\bibitem{Franke} C.H.\ Franke,
\textit{Picard-Vessiot theory of linear homogeneous difference equations},
Trans.\ Amer.\ Math.\ Soc.\ 108 (1963), 491--515.


\bibitem{Hendriks-Singer} P.A.\ Hendriks, M.F.\ Singer,
\textit{Solving difference equations in finite terms},
J.\ Symbolic Computation 27 (1999), 239--259.

\bibitem{Heyneman-Sweedler} R.G.\ Heyneman, M.E.\ Sweedler,
\textit{Affine Hopf algebras II},
J.\ Algebra 16 (1970), 271--297.

\bibitem{Kolchin1948} E.R.\ Kolchin,
\textit{Algebraic matric groups and the Picard-Vessiot theory of homogeneous
linear ordinary differential equations},
Ann.\ of Math.\ 49 (1948), 1--42.

\bibitem{Kolchin1948b} E.R.\ Kolchin,
\textit{On certain concepts in the theory of algebraic matric groups},
Ann.\ of Math.\ 49 (1948), 774--789.

\bibitem{Kolchin1973} E.R.\ Kolchin,
\textit{``Differential algebra and algebraic groups"},
Pure and Applied Mathematics 54, Academic Press, New York, 1973.

\bibitem{Levelt} A.H.M.\ Levelt,
\textit{Differential Galois theory and tensor products},
Indag.\ Mathem., N.S., 1(4) (1990), 439--450.

\bibitem{Masuoka-Yanai} A.\ Masuoka, T.\ Yanai,
\textit{Hopf module duality applied to X-outer Galois theory},
J.\ Algebra 265 (2003), 229--246.

\bibitem{Montgomery} S.\ Montgomery,
\textit{``Hopf algebras and their actions on rings"},
CBMS Reg.\ Conf.\ Series in Math.\ 82, AMS, 1993.

\bibitem{Nichols} W.D.\ Nichols,
\textit{Quotients of Hopf algebras},
Comm.\ Algebra 6 (1978), 1789--1800.

\bibitem{Okugawa} K.\ Okugawa,
\textit{Basic properties of differential fields of an arbitrary characteristic
and the Picard-Vessiot theory},
J.\ Math.\ Kyoto Univ.\ 2--3 (1963), 295--322.

\bibitem{Petkovsek} M.\ Petkov\v{s}ek,
\textit{Hypergeometric solutions of linear recurrences with polynomial
coefficients},
J.\ Symbolic Computation 14 (1992), 243--264.

\bibitem{vanderPut-Singer1997} M.\ van der Put, M.F.\ Singer,
\textit{``Galois theory of difference equations"},
Lecture Notes in Math.\ 1666, Springer, 1997.

\bibitem{vanderPut-Singer2003} M.\ van der Put, M.F.\ Singer,
\textit{``Galois theory of linear differential equations"},
Grundlehren Math.\ Wiss.\ 328, Springer, 2003.

\bibitem{Schneider} H.-J.\ Schneider,
\textit{Principal homogeneous spaces for arbitrary Hopf algebras},
Israel J.\ Math.\ 72 (1990), 167--195.

\bibitem{Sweedler1969} M.E.\ Sweedler,
\textit{``Hopf algebras"},
Benjamin, New York, 1969.

\bibitem{Sweedler1975} M.E.\ Sweedler,
\textit{The predual theorem to the Jacobson-Bourbaki theorem},
Trans.\ Amer.\ Math.\ Soc.\ 213 (1975), 391--406.

\bibitem{Takeuchi1972} M.\ Takeuchi,
\textit{A correspondence between Hopf ideals and sub-Hopf algebras},
Manuscripta Math.\ 7 (1972), 251--270.

\bibitem{Takeuchi1979} M.\ Takeuchi,
\textit{Relative Hopf modules---equivalences and freeness criteria},
J.\ Algebra 60 (1979), 452--471.

\bibitem{Takeuchi1989} M.\ Takeuchi,
\textit{A Hopf algebraic approach to the Picard-Vessiot theory},
J.\ Algebra 122 (1989), 481--509.

\bibitem{Tyc-Wisniewski} A.\ Tyc, P.\ Wi\'sniewski,
\textit{The Lasker-Noether theorem for commutative and noetherian module
algebras over a pointed Hopf algebra},
J.\ Algebra 267 (2003), 58--95.

\bibitem{Waterhouse} W.C.\ Waterhouse,
\textit{``Introduction to affine group schemes"},
Grad.\ Texts in Math.\ 66, Springer, 1979.

\end{thebibliography}

\end{document}